\documentclass[11pt]{article}
\usepackage{amssymb,amsmath,amsthm}
\usepackage{upgreek}
\usepackage{yfonts}
\usepackage{bbm}
\usepackage{mathrsfs}
\usepackage{stmaryrd}
\usepackage{graphicx}
\usepackage{color}
\newtheorem{remark}{Remark}
\newtheorem{definition}{Definition}

\newtheorem{proposition}{Proposition}
\newtheorem{lemma}{Lemma}
\newtheorem{corollary}{Corollary}
\newtheorem{theorem}{Theorem}

\newcommand{\Z}{\mathbb{Z}}

\def\rest{\hskip 1pt{\hbox to 10.8pt{\hfill
\vrule height 7pt width 0.4pt depth 0pt\hbox{\vrule height 0.4pt width 7.6pt depth 0pt}\hfill}}}

\newcommand{\eps}{\varepsilon}

\def \rc {{\rm c }}
\def \rd {{\rm d}}

\def \rba {{\rm \bf a}}

\def \rbe {{\rm \bf e}}
\def \rbo{{\rm \bf o}}
\def \rbc {{\rm \bf c}}
\def \rbd {{\rm \bf d}}
\def \rbh {{\rm \bf h}}

\def \rbu {{\rm \bf u}}
\def \rbv {{\rm \bf v}}
\def \rbq{{\rm \bf q}}
\def \rby{{\rm \bf y}}
\def \rbn{{\rm \bf n}}

\def \rbk{{\rm \bf k}}

\def \rbI {{\rm \bf I}}
\def \rbQ {{\rm \bf Q}}
\def \rbQ {{\rm \bf Q}}
\def \rbS{{\rm \bf S}}
\def \rbC {{\rm \bf C}}
\def \rbA {{\rm \bf A}}

\def \rbD{{\rm \bf D}}
\def \rbV{{\rm \bf V}}
\def \rbU{{\rm \bf U}}
\def \ureg {{\rm u_{\rm reg}}}

\def \uir{{\rm u_{\rm irr}}}

\def \vregna {{\rm v_{{\rm reg}, n}^a}}

\def \virna{{{\rm v}_{{\rm irr},  n}^a}}
\def \mbWL {{\mathbb W_{\rL}}}
\def\B{{\mathbb B}}
\def\E{{\rm E}}

\def\R{{\mathbb R}}
\def\N{{\mathbb N}}
\def\S{{\mathbb S}}

\def \rL{{\rm L}}

\def \be{\vec {\bf e}}

\def \upsigmam{{\upsigma}_{\rm main}}

\DeclareMathAlphabet{\mathpzc}{OT1}{pzc}{m}{it}

\def \bv{ {\mathbbm v }}

\def \bm{ {\mathbbm m }}

\def \rbN{{ \rm \bf  N}}
\def \rbf{{\rm \bf f}}

\def \rK {{\rm K}}
\def \rbu {{\rm \bf u}}
\def \rbO {{\rm \bf O}}
\def \rbR{{\rm \bf R}}
\def \rbJ{{\rm \bf J}}
\def \rbx {{\rm \bf x}}
\def \rbK {{\rm \bf K}}
\def \rbM {{\rm \bf M}}
\def \HL {{H^1_{\rm per}(\Lambda_\rL)}}
\def \WL {{\mathbb W_{\rm L}}}
\def \bw{{\mathbbmss w}}
\def \sech {{\rm  sech\,  }}
\def \dist{{\rm \bf dist} }
\def \Tr{{\rm \bf Tr}}

\def\QED{\hbox{${\vcenter{\vbox{
   \hrule height 0.4pt\hbox{\vrule width 0.4pt height 6pt
   \kern5pt\vrule width 0.4pt}\hrule height 0.4pt}}}$}\vspace{7pt}}

\begin{document}

\author{Fabrice BETHUEL\thanks{Sorbonne Universit\'e, CNRS, Universit\'e de Paris, Laboratoire Jacques-Louis Lions (LJLL), F-75005 Paris  
},  
Dongheng LIU  \thanks{Sorbonne Universit\'e, CNRS, Universit\'e de Paris, Laboratoire Jacques-Louis Lions (LJLL), F-75005 Paris,}
and 
Ramon OLIVER-BONAFOUX \thanks{Dipartimento di Informatica, Universit\`a di Verona. Strada le Grazie 15, 37134 Verona, Italy. E-mail: \texttt{ramon.oliverbonafoux@univr.it}}
}
\title{Pseudo-profiles for vectorial Allen-Cahn systems}
\date{}
\maketitle
\begin{abstract}
The purpose of the present paper is to present a new type of solutions to the elliptic  vectorial Allen-Cahn system of  equations. Its main feature is that in the classical  rescaling its asymptotics do not follow the rule of balance of energies, which is crucial in the analysis of the asymptotics for the scalar problem.

\end{abstract}

\section {Introduction}
\label{introduction}
\subsection{ The  Allen-Cahn  equation on $\R^N$}
The purpose of the present paper is to construct special \emph{entire} solutions  $u: \R^N\to \R^k$,  for $N \in \N^\star$,  to the vectorial  Allen-Cahn equation
\begin{equation}
 \label{elip0}
  -\Delta u(x)=-\nabla_u V(u(x))  {\rm \ for  \ } x \in \R^N. 
  \end{equation}
The function $V$ appearing  in \eqref{elip0},  termed the \emph{potential}, denotes   a smooth scalar   function on   $\R^k$, where $k \in \N^\star$ is a given integer,  and 
$\displaystyle{\nabla_u V(\rbu)=\left(\frac{\partial V}{\partial u_1}(\rbu),\ldots  \frac{\partial V}{\partial u_k}(\rbu) \right)}$, for any $\rbu=(u_1, \ldots, u_k)$. 
Equation \eqref{elip0} is  hence a \emph{system of $k$ scalar partial differential equations}  for each of the components of the map $u=(u_1, \ldots , u_k)$. The system of Equations \eqref{elip0}   corresponds to the  Euler-Lagrange equation   of the  energy functional  ${\E}$ which is defined for a function $u$ and any subdomain $\Omega \subset \R^N$  by the formula
    \begin{equation}
 \label{glfunctional0}
 \E (u, \Omega)= \int_\Omega e(u(x))\rd x, {\rm \ where \ } e (u(x))=\frac{\vert \nabla u (x)\vert ^2}{2}+ V(u(x)), {\rm \ for \ } x \in \R^N
 \end{equation}
  and $\E(u)=\E(u, \R^N)$.
We notice  also that, if $u$ is a solution to \eqref{elip0} on $\R^{N}$, then the map  $\tilde u: \R^{N+1}{\color{black}\to}\R^k$ defined by
\begin{equation}
\label{plusun}
\tilde u(\rbx, x_{N+1})=u(\rbx), {\rm \ for \ any \ } \rbx\in \R^N,{\rm \ and \ any \ } x_{N+1} \in \R,
\end{equation}
is a solution to \eqref{elip0} in dimension $N+1$. Hence solutions obtained in dimension $N$ allow us to construct  related solutions in any dimension $m\geq N$.

\smallskip
  We  assume throughout that the potential $V$ is bounded below, so that we may impose, without loss of generality and changing possibly $V$ by a suitable  additive  constant,  
     \begin{equation}
   \label{infimitude}
   \inf V=0.
   \end{equation}
 We introduce the  set $\Sigma$ of minimizers of $V$,  sometimes called the \emph{vacuum manifold}, that is the  subset of $\R^k$ defined by
     $$\Sigma\equiv \{ y \in \R^k, V(y)=0 \}.$$
We assume also throughout that the set $\Sigma$ is finite and has at least two minimizers.   A typical example is given in the scalar case by the 
\emph{Allen-Cahn} potential
\begin{equation}
\label{scalarAC}
V_{\rm AC}(u)=\frac{(1-u^2)^2}{4}, {\rm \ for \ }  u \in \R, 
\end{equation}
for which $\Sigma=\Sigma_{\rm AC}=\{+1, -1\}$ has exactly two elements. 
Other, more technical,  assumptions on $V$ required by our statements will be introduced in due place. 

\smallskip
Given $\eps>0$, we introduce the  scaled  function  $u_\eps$  defined by
\begin{equation}
\label{ueps}
u_\eps(x)= u( \frac{x}{\eps}), {\rm \ for  \  any \ } x\in {\color{black}\R^N}, 
\end{equation}
 so that, if $u$ is a solution to \eqref{elip0} then $u_\eps$ is a solution to
  \begin{equation}
 \label{elipes}
  -\Delta u_\eps=-\eps^{-2}\nabla_u V(u_\eps)  {\rm \ on  \  }  {\color{black}\R^N},  
  \end{equation}
 The asymptotic behavior as $\eps \to 0$ of  families of solutions $(u_\eps)_{\eps>0}$ of  \eqref{elipes}   has been widely investigated, in particular in the scalar case $k=1$.
 Equation \eqref{elipes}   corresponds to the scaled   Euler-Lagrange equation   of the  energy functional  ${\E_\eps}$ which is defined for a function on any subdomain $\Omega \subset \R^N$  by the formula
    \begin{equation}
 \label{glfunctional}
 \E_\eps(u, \Omega)= \int_\Omega e_\eps(u(x))\rd x, {\rm \ where \ } e_\eps (u(x))=\frac{\eps\vert \nabla u (x)\vert ^2}{2}+\frac{1}{\eps} V(u(x)), {\rm \ for \ } x \in \R^N.
 \end{equation}
 The  asymptotic analysis of families  $(u_\eps)_{\eps>0}$ of  \eqref{elipes} involves several measures, in particular
\begin{equation}
 \upnu_\eps= e_\eps(u_\eps)\rd x, \, \upzeta_\eps \equiv \frac{  V(u_{\eps})}{\eps}  {\rm d} x,  {\rm \ and \ }
 \upmu_{\eps, i, j}=\eps {u_{\eps}}_{x_i}\cdot  {u_{\eps}}_{x_j} \rd x, {\rm \ for \ } i, j=1, \ldots, N.
\end{equation}

   \subsection{ The one-dimensional case: profiles}
  
  \noindent
  For $N=1$,    equation \eqref{elip0} reduces   to the ordinary differential equation
\begin{equation}
	\label{heterocline}
	\frac{{\rm d}^2 w}{{\rm d} s
		^2}=\nabla_w V(w) {\rm \ on \ } \R.
\end{equation} 
Finite energy solutions to \eqref{heterocline} necessarily connect at $\pm \infty$ two  minimizers {$\upsigma^-$} and {$\upsigma^+$}: They are called 
{\emph{profiles}} or {\emph{heteroclinic connections}}, if {$\upsigma^-\not = \upsigma^+$}. Multiplying \eqref{heterocline} by 
$\dot w$, we are led to the conservation law 
 \begin{equation}
 \label{conserve}
 \frac{\rd}{\rd s} \left(V(w)-\frac{\vert \dot w\vert^2}{2}\right)=0,  
 \end{equation}
so that for  solutions to \eqref{heterocline}  one derives  the identity
 \begin{equation}
 \label{first-order}
  \vert \dot w\vert=\sqrt{2V(w)} {\rm \ on \ } \R, 
   \end{equation}
    referred to as \emph{the balance of energies}. Notice also the scale invariance
    \begin{equation}
     \E_\eps(w_\eps, \R)= \E(w, \R), {\rm  \ for \ } \eps>0, {\rm \ with \   \E=\E_1}.
    \end{equation}
   In the  \emph{scalar} case, the first order equation \eqref{first-order} reduces to  $w_s=\pm \sqrt{2V(w)}$, which can be    integrated by separation of variables, so  that profiles  connect only \emph{nearby minimizers} $\upsigma^-$ and $\upsigma^+$ of the potential, and  are essentially unique,  up to translations and symmetries. For instance, in the case of the   \emph{Allen-Cahn potential} \eqref{scalarAC}, the solution  is given,  up to translation and symmetry, by  
\begin{equation}
\label{kink}
w(s)= \tanh \left(\frac{s} {\sqrt{2}} \right), {\rm \ for \ } s \in \R.
\end{equation}
The situation is  different in the \emph{vectorial case}, since relation \eqref{first-order} is less constraining: Under additional assumptions on the potential $V$, one may find several profiles connecting two minimizers of the potential. The search for such solutions is an active field of research (see e.g. \cite{ Alibetu,  Zuster, monsant, bonafoux}).

\smallskip
  Given $N >1$, a solution  $w$ to \eqref{heterocline}  gives rise to a solution $u$ to \eqref{elip0} on $\R^N$, as mentioned earlier, setting 
  $u(x_1, \ldots, x_N)=w(x_N)$. For $\eps>0$, we then consider  the map $u_\eps$ provided by \eqref{ueps}. 
We observe that the energy is locally bounded independently of $\eps>0$, that is $\E_\eps(u_\eps, \B(R)) \leq C(R)$, for any $R>0$ and $\eps >0$, and that, in the sense of measures  
\begin{equation}
\label{lilimit}
\left\{
\begin{aligned}
&\underset{\eps \to 0}\lim \upmu_{\eps, i, i}=0, {\rm \ for  \ } i=1, \ldots, N-1 \\
& \underset{\eps \to 0} \lim \upnu_\eps= \rc_0 \rd \lambda, {\rm \  where \ }  \rd  \lambda=\mathcal H^{N-1}\rest  P_N, \\
& 2\underset{\eps \to 0} \lim \, \upzeta_\eps=\underset{\eps \to 0} \lim \, \upmu_{\eps, N, N}=\underset{\eps \to 0} \lim \upnu_\eps
=\rc_0 \rd \lambda, 
\end{aligned}
\right.
\end{equation}
where $P_N$ is the hyperplane $P_N=\{(x_1, \dots,  {\color{black}x_N})\in \R^N, x_N=0\}$ and $\rc_0=\E(w, \R)$ is the energy of the profile. The third line in \eqref{lilimit} is a consequence of balance of energies  and will be referred to as the \emph{asymptotic balance of energies}.  Setting 
$\upzeta_\star=\underset{\eps \to 0} \lim \, \upzeta_\eps$, and $\upnu_\star=\underset{\eps \to 0} \lim \upnu_\eps$,  the asymptotic balance of energies is written 
hence as 
\begin{equation}
\label{asymptobal}
2\upzeta_\star=\upnu_\star.
\end{equation}
   \subsection{Profiles in the scalar case}

\noindent
In the scalar case, the theory developed in the higher-dimensional case $N\geq 1$, shows that, for a larger part, the asymptotic situation described in the example \eqref{lilimit} corresponds, at least vaguely,  to  the  general situation. Indeed, 
 the results in  \cite{hutchtone}  embody  the intuitive idea that locally, the equation reduces to a one-dimensional problem. More precisely,  it asserts that there exists some $N-1$ dimensional stationary varifold $\mathfrak S_\star$ ({\color{black}a weak} notion of {\color{black}minimal} surface), such that  up to some subsequence, 
   the expected situation reduces,  \emph{locally near some point}  $x_0\in \mathfrak S_\star$,  to the case
\begin{equation}
\label{quoique}
u_\eps(x)\underset{x \to x_0} \simeq w_\eps (x_N),  {\rm \ with \ }  x=(x_1,\ldots,  x_N) \in {\color{black}\R^N}, 
\end{equation} 
where the coordinates are chosen so that the tangent to $\mathfrak S_\star$ at $x_0$ has equation $x_N=0$, and where $w_\eps$ stands for {\color{black}an} $\eps$-scaled version of the solution to the one-dimensional problem \eqref{heterocline}. The possibility of gluing of several such one-dimensional solutions is not excluded, but we will not discuss this here.   Ultimately, the results in  \cite{ilmanen, hutchtone} provide a rather simple picture of the solutions. They involve a minimal surface, the solution may be represented  as  one-dimensional profiles glued to the surface in the transversal direction, so that one is  tempted to write the correspondence
\begin{equation}
\label{sim}
 {\rm \ solutions \ to \   \eqref{elipes} }  \sim {\rm \ minimal \ surface} + {\rm glued \ profiles}.
 \end{equation}
 The general structure of solutions is hence fairly well understood.

 \subsection{The  vectorial case : pseudo-profiles and solutions on cylinders}
\label{theprogram}

The picture  in the   vectorial case is \emph{more complex}. Firstly, as seen, the set of one-dimensional profiles is much larger, it could perhaps be even infinite.  Besides this, in a case where the set of minimizing solutions contains at least two elements and  the potential $V$ satisfies some suitable conditions, there exists a solution on $\R^2$, concentrating near $P_2$,  which  \emph{cannot be reduced to a single   profile},  but performs a gradual transition between two minimizing solutions. This was first proven  in \cite{alabrogui} for potentials satisfying a symmetry assumption{\color{black}.} It was later proven in \cite{schatzman} (see also \cite{monsant, smyrnelis} for more recent proofs through
diﬀerent approaches) that such a symmetry assumption can be removed provided that the
minimizing heteroclinic solutions are \emph{(spectrally) non-degenerate up to translations}, an assumption  which
holds for generic potentials when $k=2$. More precisely, in \cite{alabrogui} (see also \cite{schatzman}), it is proved under suitable assumptions on $V$, 
that  there exists a \emph{locally minimizing} solution $u: \R^2 \to \R^k$  of 
  $\displaystyle{-\Delta u+\nabla_u V(u)=0}$  such that 
  \begin{equation*}
  \left\{
  \begin{aligned}
  u(x_1, x_2) & \to w_1(x_2) {\rm \ as \ } x_1  \to -\infty \\
  u(x_1, x_2) & \to w_2(x_2) {\rm \ as \ } x_1  \to +\infty,  \\ 
  \end{aligned}
\right.
  \end{equation*} 
  where $w_1$ and $w_2$ are two distinct minimizing solutions of \eqref{heterocline}.
  For $\eps>0$,   the {scaled map $u_\eps$} is defined by 
 $\displaystyle{u_\eps(\rbx)=u\left(\frac{\rbx}{\eps}\right), \rbx  \in \R^2}$.
The map   $u_\eps $ is {locally minimizing for $\E_\eps$} and 
concentrates on the line  {\bf $\mathfrak S_\star=\{(x_1, 0), x_1 \in \R\}$}, and 
\begin{equation*}
\left\{
\begin{aligned}
u_\eps(x_1, x_2) & \simeq  w_1(\frac{x_2}{\eps}) {\rm \ for  \ } x_1<0 . \\
  u_\eps(x_1, x_2) & \simeq w_2(\frac{x_2}{\eps}) {\rm \ as \ } x_1>0. \\ 
\end{aligned}
\right. 
\end{equation*}
hence, \eqref{sim}  is verified only on  the \emph{half-lines}, not on the whole line.  Notice however that relations \eqref{lilimit} are satisfied by the family 
$(u_\eps)_{\eps>0}$, hence the balance of energies holds, at least in the asymptotic sense. 

\smallskip
In the present paper we present a different construction of solutions to $\displaystyle{-\Delta u+\nabla_u V(u){\color{black}=0}}$ on $\R^2$ which we refer to as \emph{pseudo-profiles}. An important feature and main motivation for this construction is the fact  that the  family of scaled maps $(u_\eps)_{\eps>0}$ do not satisfy \emph{the  balance of energies in the asymptotic sense}, hence providing a (\emph{new}) major difference between the scalar and the vectorial case. Notice that the assumptions on the potential  $V$ that we are going to choose are quite close to the ones introduced in \cite{schatzman}. Indeed, we do not impose any symmetry assumption on $V$, but instead   we rely in some places on the spectral properties of $w_1$ and $w_2$. We introduce the cylinder

 $\displaystyle{\Lambda_\rL=[-\frac{\rL}{2}, \frac{\rL}{2}] \times \R}$ with $L>0$, and impose periodic boundary conditions, that is, we consider the problem: Find $u \in C^2(\Lambda_\rL)$ such that $u$ is a solution to
 \begin{equation}
 \label{themainthing} 
 \left\{
 \begin{aligned}
 & -\Delta u +  \nabla_u V(u)=0  {\rm \ on \ } \Lambda_\rL, \\
   & u(-\frac \rL 2, x_2)=u(\frac \rL 2, x_2), {\rm \ for \ any \ }  x_2 \in \R {\rm \ and \ } \\
   & \frac{\partial u}{\partial x_1}(-\frac{\rL}{2}, x_2)= \frac{\partial u}{\partial x_1 }(\frac \rL 2, x_2), {\rm \ for \ any \ }  x_2 \in \R.
  \end{aligned}
\right. 
  \end{equation}
An alternate and equivalent approach is to study solutions on the tube $\displaystyle{\R\slash\rL\Z\times \R}$.  Given a solution $u \in C^2(\Lambda_\rL)$   which satisfies \eqref{themainthing},  we extend the solution $u$  as a map $\rbu$ on $\R^2$, by periodicity, setting
\begin{equation}
\label{periodisation}
\rbu(x_1+ \rL N, x_2) = u(x_1, x_2), {\rm \ for \ any \ } x_1\in [-\frac{\rL}{2}, \frac{\rL}{2}], N\in \Z, x_2\in \R, 
\end{equation}
 so that $u$ is now a solution of \eqref{elip0} on  the whole of $\R^2$. Finally, we consider the scaled map
defined   on $\R^2$ by,  for $x=(x_1, x_2) $  
\begingroup\color{black}
\begin{equation}
\label{glasnost}
\rbu_\eps(x)=\rbu\left(\frac{x}{\eps}\right),\qquad x\in\R^2.
\end{equation}
\endgroup
 If $u$ has finite energy on $\Lambda_L$,  that is $\E(u, \Lambda_L)< \infty$, then  by scaling we have 
 \begin{equation}
 \E_\eps (u_\eps, \Lambda_{\eps L})=\eps\E(u, \Lambda_L).
  \end{equation}
  It follows that
  \begin{equation}
  \label{ouistiti}
  \left\{
  \begin{aligned}
  \eps  \left  \vert  \frac{\partial u_\eps}{\partial x_1}  \right  \vert^2  \rightharpoonup \upmu_{\star, 1,1}&=\rc_{1, 1} \mathcal H^1\rest D, {\rm \ with \ } D=\left\{(x_1, 0), x_1\in \R \right\},  \\
{\rm \ where \  }  \rc_{1, 1}&=\frac{1}{\rL}\int_{\Lambda_L} \left( \frac{\partial u}{\partial x_1}\right)^2 \rd x.
\end{aligned}
\right.
 \end{equation}
 In particular, if
 \begin{equation}
  \label{bla}
 \frac{\partial u}{\partial x_1}\not  =  0  {\rm \ on \ } \Lambda_\rL, {\rm \ so \ that \ } \rc_{1, 1}\not=0, 
  \end{equation}
then the  solution $u$  is \emph{not a one-dimensional profile}.   Indeed, if $u$ is of the form $u(x_1, x_2)=w(x_2)$,  for any $(x_1, x_2) \in \Lambda_\rL$,  then necessarily  $\frac{\partial u}{\partial x_1}= 0$.

\begin{remark}{\rm  In \cite{bethu}, it was established that, for families $(u_\eps)_{\eps>0}$ of solutions to \eqref{elipes} with suitable energy bounds and limiting interface contained in the horizontal line $D$, one has the identities (where the measures are defined as in \eqref{lilimit} and up to a possible subsequence)
\begin{equation}
\label{discrepvec11}
2\upzeta_\star=\upmu_{\star, 2,  2}-\upmu_{\star, 1, 1} =\upnu_\star- \upmu_{\star, 1, 1}.
\end{equation}
In particular, if the family $(u_\eps)_{\eps>0}$  is given by \eqref{glasnost}, with $u$ being a solution to \eqref{themainthing}, then
\begin{equation}
\label{graal}
2\upzeta_\star=\upnu_\star-{\color{black}\rc_{1, 1}} \mathcal H^1\rest D,
{\rm \ with \ }
\rc_{1, 1}=\frac{1}{\rL}\int_{\Lambda_L} \left( \frac{\partial u}{\partial x_1}\right)^2 \rd x.
\end{equation}
}
\end{remark}  
 
The remainder of the paper is devoted to the existence, under suitable assumptions on $V$, of a solution to \eqref{themainthing} satisfying \eqref{bla}.

\subsection{Assumptions on the potential}
 Our main result, which provides a non constant solution to \eqref{themainthing} and \eqref{bla},  requires assumptions on the potential $V$ which are parallel to {\color{black}those} in \cite{schatzman}.
 As {\color{black}mentioned}, we assume that the infimum of $V$ is zero,  the set of minimizers is non-empty, finite, that is 
    \begin{equation}
   \label {H1}
    {\color{black}(\mathbf{H}_{\mathbf{1}})} 
    \left\{
    \begin{aligned}
 & \inf V=0,   V \in C^3(\R^k)\\
 & V(x) \to + \infty {\rm  \ as \ } \vert x \vert \to + \infty,   {\rm \ and \ } \\
 &\Sigma\equiv \{ y \in \R^k, V(y)=0 \}=\{\upsigma_1, ..., \upsigma_{\mathfrak q_\star}\},\ \mathfrak q_\star \geq 2, \ \upsigma_i \in \R^k, \ \forall i=1,...,\mathfrak q_\star, 
  \end{aligned}
  \right.
   \end{equation}
   for some $\mathfrak q_\star \in \N^\star$.  We impose the following  additional conditions on the potential $V$. We first impose the behavior at infinity:
   
   \smallskip
\noindent
  ${\color{black}(\mathbf{H}_{\mathbf{2}})}$ {\it There exist  constants  $\upalpha_\infty >0$ and  $R_\infty >0$  such that  
  \begin{equation}
  \label{condinfty}
  \left\{
   \begin{aligned}
&   y\cdot\nabla V( y )\geq \upalpha_\infty \vert y \vert ^2, \ \hbox {if }  \vert y \vert >R_\infty {\rm \  and \ }  \\
 &V(x) \to +\infty {\rm \ as  \  }  \vert x \vert \to + \infty, \\
 &\vert V(x)\vert + \vert  \nabla V(x)  \vert + \vert  \nabla^2V(x)  \vert \leq C_ m\left( \vert x \vert^m+1\right), {\rm \ for \ some \ } C_m>0, {\rm \ and  \ } m \geq 3.
 \end{aligned}
 \right.
 \end{equation}   
 }

 Next, we specify the behavior near the points of $\Sigma$:
   
    \smallskip
    \noindent
${\color{black}(\mathbf{H}_{\mathbf{3}})}$   {\it  The matrix $\nabla^2V(\upsigma_i)$ is positive definite at each point $\upsigma_i$ of $\Sigma$, in other words, if $\lambda_i^-$ denotes its smallest eigenvalue, then  $\lambda_i^->0$. We denote by $\lambda_i^+$ its largest eigenvalue. }

\medskip 
 The two next  assumptions are of different nature, and deal with  properties of \emph{minimizing one-dimensional profiles}, that is solutions to the ODE \eqref{heterocline}.  Given two distinct elements $\upsigma^-$ and $\upsigma^+$ in $\Sigma$, we set 
 \begin{equation}
 \label{Xsigma}
 \mathbb X(\upsigma^-, \upsigma^+)=\left\{ v\in H^1_{\rm loc}(\R, \R^k), {\rm \ such \ that \ } \E(v, \R)<+\infty, {\rm \ and \ } v(s) \underset {s\to \pm{\color{black}\infty}} \to \upsigma^\pm \right\}.
\end{equation}
 Clearly $\mathbb X(\upsigma^-, \upsigma^+)$ is non-empty.   Set 
 \begin{equation}
 \label{mini1D}
  \bm(\upsigma^-, \upsigma^+)=\inf \left\{ \E(v, \R), v \in \mathbb X(\upsigma^-, \upsigma^+)\right\}. 
 \end{equation}
 If the infimum is achieved in \eqref{mini1D}, then a minimizer is a  solution to \eqref{heterocline} on $\R$.
Notice that we have invariance by translation: For $a\in \R$,  and $f : \R \to \R^k$, set $\tau_a(f)=f(x-a)$.  By invariance of the energy by translations, we have  for any 
 $v \in \mathbb X(\upsigma^-, \upsigma^+)$, and any $a \in \R$, $\tau_a(v) \in  \mathbb X(\upsigma^-, \upsigma^+)$ and 
 $$\E(\tau_a (v), \R)=\E(v, \R), {\rm \ for \ any \ } a \in \R. $$
  We are therefore led to  introduce  the equivalence class, for $w\in  \mathbb X(\upsigma^-, \upsigma^+)$
\begin{equation}
\label{equivclass}
\mathcal C(w)=\{z \in \mathbbmss X(\upsigma^-, \upsigma^+), {\rm \ such \ that \ } \exists a\in \R, {\rm \ such \ that \ } 
z(x)=w(x-a), \forall x \in \R\}.
\end{equation}
Since the equation \eqref{heterocline} is invariant by translations, if $w$ is a solution to \eqref{heterocline}, so is any element in 
$\mathcal C(w)$, and, if moreover, $w$ is minimizing, so is any element in the equivalence class $\mathcal C(w)$.
  
 \medskip
 \noindent
  ${\color{black}(\mathbf{H}_{\mathbf{4}})}$  There exist two elements $\upsigma_0^+$ and $\upsigma_0^-$ in $\Sigma$  and  two distinct minimizers $w_1$ and $w_2$ in $\mathbb X(\upsigma_0^-, \upsigma_0^+)$ of $\E$,  such that $w_2 \not\in \mathcal  C(w_1)$. 
 Moreover,  there is no other minimizing solution outside of $\mathcal C(w_1)\cup \mathcal C(w_2)$.
  We assume additionally  that 
  \begin{equation}
  \label{strictosenso}
  \bm(\upsigma_0^-, \upsigma_0^+)<\underset{k=1}{\overset{\mathfrak p}\sum}\bm(\upsigma_k^{-}, \upsigma_k^+)
  \end{equation}
  for any $\mathfrak p\geq 2$ in $\N^\star$, and for $k=1, \ldots, \mathfrak p$, $\upsigma_k^\pm\in \Sigma$ satisfying the compatibility condition
  \begin{equation}
  \label{quadrilatere}
  \left\{
  \begin{aligned}
  &\upsigma_1^-=\upsigma_0^- {\rm \ and \ } \upsigma_{\mathfrak p}^+=\upsigma_0^+, {\rm \ and \ } \\
 &\upsigma_k^+=\upsigma_{k+1}^-{\color{black},}\, {\rm \ for \ } 
 k=1,  \ldots, \mathfrak p-1{\color{black},}\\
&{\color{black}\upsigma_k^-\not=\upsigma_k^+,\quad k=1,\ldots,\mathfrak p.}
   \end{aligned}
 \right.
  \end{equation}
\medskip 

\begin{remark}{\rm  Condition \eqref{strictosenso} is clearly satisfied if the number of wells $\mathfrak q_\star$ is $2$, or if 
   $$\bm(\upsigma_0^-, \upsigma_0^+)
   =\inf\{\bm(\upsigma^-, \upsigma^+),  (\upsigma^-, \upsigma^+)\in \Sigma^2, \ \upsigma^-\not =\upsigma^+\}.$$
The explicit example  provided in the Appendix falls in the case $\mathfrak q_\star=2$.
}
\end{remark}

\begin{remark}
{\rm
It is well-known that, under {\color{black}condition \eqref{quadrilatere}}, we have the  classical, triangle-type inequality  
 \begin{equation}
 \label{v3}
  \bm(\upsigma_0^-, \upsigma_0^+)=\bm (\upsigma_1^-, \upsigma_{\mathfrak p}^+) \leq
  \underset{k=1}{\overset {\mathfrak p} \sum }\bm(\upsigma_k^-, \upsigma_k^+), 
 \end{equation}
so that the strict inequality in \eqref{strictosenso} is crucial in our arguments.
}
\end{remark}

\begin{remark}
\label{remimportante}
{\rm   For given $\mathfrak p\geq 1$,
 consider $(\upsigma_1, \ldots,  \upsigma_{\mathfrak p} )\in \Sigma^{\mathfrak p}$, and define
 $\upsigma_1^-=\upsigma_0^-$, \ $\upsigma_{\mathfrak p}^+=\upsigma_0^+$, 
 $\upsigma_k^-=\upsigma_{k-1}$,  for $k=2, \ldots, \mathfrak p$,  and $\upsigma_k^+=\upsigma_{k}$, for $k=1, \ldots, \mathfrak p-1$, 
 so that {\color{black}the endpoint and matching conditions in \eqref{quadrilatere} hold.}
 We then set
 \begingroup\color{black}
\begin{equation}
\label{deltastar}
\left\{
\begin{aligned}
&\rbS(\upsigma_1,\ldots,\upsigma_{\mathfrak p})
 =\sum_{k=1}^{\mathfrak p}\bm(\upsigma_k^-,\upsigma_k^+),
 \qquad(\upsigma_1,\ldots,\upsigma_{\mathfrak p})\in\Sigma^{\mathfrak p},\\
&\updelta_\star=
 \inf_{\mathfrak p\geq2}
 \inf_{\substack{(\upsigma_1,\ldots,\upsigma_{\mathfrak p})\in\Sigma^{\mathfrak p}\\
 \upsigma_k^-\not=\upsigma_k^+,\ 1\leq k\leq\mathfrak p}}
 \rbS(\upsigma_1,\ldots,\upsigma_{\mathfrak p})
 -\bm(\upsigma_0^-,\upsigma_0^+).
\end{aligned}
\right.
\end{equation}
\endgroup
{\color{black}Since $\Sigma$ is finite and every non-trivial edge has positive minimal energy, only finitely many admissible chains have total energy at most $\bm(\upsigma_0^-,\upsigma_0^+)+1$. Applying \eqref{strictosenso} to these chains gives}
\begin{equation}
\label{deltastarplus}
\delta_\star>0.
\end{equation}
}
\end{remark}

 \begin{remark}{\rm   {\color{black}It is well-known that (see, e.g., inequality \eqref{weirdos} in Section \ref{one-d})}
 \begin{equation}
  \bm(\upsigma^-, \upsigma^+)=0 {\rm \ if \ and \ only \ if \ } \upsigma^+=\upsigma^-{\color{black}.} 
 \end{equation}
}
 \end{remark}

 Our last assumption involves  the second derivative of the energy $\E$  at the one-dimensional minimizing profiles $w_1$ and $w_2$ and more precisely,  their  spectral properties.  Notice that given any $v \in \mathbb X(\upsigma^-, \upsigma^+)$, we may expand the energy at second order, 
 $\rbh=(h_1, \ldots, h_k) \in H^1(\R, \R^k)$, as 
 \begin{equation}
 \label{svilupesque}
 \E(v+\rbh)=\E(v)+\langle \rd \E(v), \rbh\rangle + \frac{1}{2}\rd^2 \E(v)(\rbh, \rbh) +
 \underset{\Vert \rbh \Vert_{H^1}\to 0} \rbo(\Vert h \Vert_{H^1}^2)
 \end{equation}
 with 
 \begin{equation}
 \label{derdederdeux}
 \rd^2E(v)(\rbh, \rbh)=\int_\R \vert \rbh'(s)\vert^2\rd s +\int_\R \rbh(s) \nabla^2 V(v(s)) \rbh(s) \rd s.
 \end{equation}
 If $v$ is a minimizing solution, then, in view of \eqref{svilupesque}, we have  $\rd^2\E(v)(\rbh, \rbh)\geq 0$, for any $\rbh \in H^1(\R)$.
 If $\rbh \in H^2(\R, \R^k)$, identity {\color{black}\eqref{derdederdeux}} can be rephrased as
 \begin{equation}
 \rd^2E(v) (\rbh, \rbh)=\langle \rbA(v) (\rbh), \rbh\rangle_{L^2(\R)}, {\rm \ where \ } 
 \rbA(v)(h)=-\frac{\rd^2 \rbh}{\rd s^2}+\nabla^2 V(v)\rbh, 
 \end{equation}
 so that $\rbA(v)$ has no negative eigenvalue, if $v$ is a minimizing solution. However, the kernel {\color{black}is} never reduced to zero. Indeed,
 \begin{equation}
 \label{coree}
  \rbA(v) v'=0, {\rm \ for \ any \  finite \ energy \ solution \ } v   {\rm \ to \ } \eqref{heterocline}.
 \end{equation}
  In order to establish \eqref{coree} it suffices to check that  $v'\in H^1(\R)$  and to {\color{black}differentiate} \eqref{heterocline}. 
 
 \smallskip
 In this context,  our last assumption is the following:

  \medskip
 \noindent
  ${\color{black}(\mathbf{H}_{\mathbf{5}})}$ The kernel of $\rbA(w_i)$ is one-dimensional, for $i=1, 2$,  and we have 
 \begin{equation}
 \label{noyau}
  \ker \rbA(w_i)=\R {\color{black}w_i'},  {\rm \ for \ }  i=1, 2.
  \end{equation}
    
  \begin{remark}{\rm   As noticed in \cite{schatzman}, it turns out that   Assumption ${\color{black}(\mathbf{H}_{\mathbf{5}})}$ is equivalent to the existence of  some number $\upalpha_i >0$ such that, for any $\rbh \in H^1(\R)$, we have
     \begin{equation}
  \label{lespectre0}
   \rd^2\E(w_i)(\rbh, \rbh)\geq {\color{black}\upalpha_i} \Vert \rbh \Vert_{H^1(\R)}^2, {\rm  \ provided \ } \langle \rbh, w_i'\rangle_{L^2(\R)}=0, 
 \end{equation} 
 see Proposition \ref{spectralou}.

 }
 \end{remark}
 \begin{remark}
 {\rm Notice that assumptions ${\color{black}(\mathbf{H}_{\mathbf{1}})}$, ${\color{black}(\mathbf{H}_{\mathbf{2}})}$, ${\color{black}(\mathbf{H}_{\mathbf{3}})}$, ${\color{black}(\mathbf{H}_{\mathbf{4}})}$  and ${\color{black}(\mathbf{H}_{\mathbf{5}})}$, are quite  standard in the literature on this subject, and we rely in many places on existing {\color{black}results}, in particular {\color{black}provided} by  \cite{schatzman, alabrogui, bonafoux, bethu}.
  In particular, let us emphasize Assumption ${\color{black}(\mathbf{H}_{\mathbf{5}})}$ was considered before in \cite{schatzman}, see also \cite{monsant}.  Moreover, as shown in \cite{schatzman}, Theorem 4.3 and Remark
4.4,  Assumption ${\color{black}(\mathbf{H}_{\mathbf{5}})}$ is {\color{black}generic}, at least in the case $k=2$.
  }
 \end{remark}

 \medskip
 In the Appendix, we produce an \emph{explicit example} of a potential $\rbV$ (for $k=2$) which fulfills all the required assumptions, that is 
 ${\color{black}(\mathbf{H}_{\mathbf{1}})}$, ${\color{black}(\mathbf{H}_{\mathbf{2}})}$, ${\color{black}(\mathbf{H}_{\mathbf{3}})}$, ${\color{black}(\mathbf{H}_{\mathbf{4}})}$  and ${\color{black}(\mathbf{H}_{\mathbf{5}})}$. In \cite {schatzman}, a {\color{black}genericity} result was given.

 \subsection{Statement of the main result}
 The main aim of this paper is to establish the following:
 
 \begin{theorem}
 \label{maintheo}
 Assume that $V$ satisfies assumptions ${\color{black}(\mathbf{H}_{\mathbf{1}})}$, ${\color{black}(\mathbf{H}_{\mathbf{2}})}$, ${\color{black}(\mathbf{H}_{\mathbf{3}})}$ and
 ${\color{black}(\mathbf{H}_{\mathbf{4}})}$ and ${\color{black}(\mathbf{H}_{\mathbf{5}})}$. 
 Then, there exists $\rL_\star>0$, such that, if $\rL\geq \rL_\star$, then there  exists a solution $u \in C^2(\Lambda_\rL)\cap \mathbb W_\rL (\upsigma_0^-, \upsigma_0^+)$  of \eqref{themainthing}  satisfying \eqref{bla}.
  \end{theorem}  
   
  As a consequence, the program   announced in Section \ref{theprogram} is fulfilled, in the case where the potential $V$ satisfies assumptions 
  ${\color{black}(\mathbf{H}_{\mathbf{1}})}$, ${\color{black}(\mathbf{H}_{\mathbf{2}})}$, ${\color{black}(\mathbf{H}_{\mathbf{3}})}$, ${\color{black}(\mathbf{H}_{\mathbf{4}})}$  and ${\color{black}(\mathbf{H}_{\mathbf{5}})}$. Indeed,  following the constructions \eqref{periodisation} and \eqref{glasnost}, we obtain a  sequence of maps $(\rbu_\eps)_{\eps>0}$  satisfying \eqref{ouistiti} and \eqref{graal}, hence contradicting the asymptotic balance of energies, contrary to the scalar case.
    
    \subsection{Some elements  of the proof of Theorem \ref{maintheo}}
    Many of our arguments are close or parallel to those used in \cite{rabinowitz, bonafoux} {\color{black}in} the one-dimensional case.  In particular, an important part of our argument is to transpose them to the context of our cylinders.
    
    \subsubsection{Functional framework}
The construction of the solution  relies on   \emph{variational} arguments. The functional framework is provided by the set $\mathbb W_\rL$, for given $\rL>0$, which is defined by  
\begin{equation}
\label{functionalset}
\mathbb W_\rL=\left\{v\in H^1_{{\color{black}\rm loc}}(\Lambda_\rL, \R^k),     \E(v, \Lambda_\rL)<+\infty, {\rm \ s.t \ }  v(-\frac{\rL}{2}, x_2)=v(\frac{\rL}{2}, x_2), \forall x_2 \in \R\right\}, 
\end{equation}
where the gradient of $v$ is defined in the distributional sense.  As for the one-dimensional case, we may define asymptotic values as $x_2 \to \pm \infty$.   For that purpose set, for $\rba \in \R^k$, and $v \in \mathbb W$, 
\begin{equation}
\label{rbI}
\rbI^\pm (v, x_2, \rba)\equiv \int_{[-\frac \rL 2, +\frac \rL 2] \times \{s\in\R:\ \pm(s-x_2)\geq0\}} \vert v-\rba \vert^2\rd x.
\end{equation}

\begin{lemma} 
\label{asympto2D0}
Assume that $V$ satisfies  ${\color{black}(\mathbf{H}_{\mathbf{1}})}$, ${\color{black}(\mathbf{H}_{\mathbf{2}})}$ and  ${\color{black}(\mathbf{H}_{\mathbf{3}})}$.
 Let $v$ be in $\mathbb W_\rL$. Then there exist elements $\upsigma^-$ and $\upsigma^+$ in $\Sigma$ such that 
\begin{equation}
\label{asympto2D}
\rbI^\pm(v, x_2, \upsigma^\pm)\equiv \int_{[-\frac \rL 2, +\frac \rL 2] \times \{s\in\R:\ \pm(s-x_2)\geq0\}} \vert v-\upsigma^\pm \vert^2\rd x \to 0 {\rm \ as \ } x_2\to\pm \infty. 
\end{equation}
\end{lemma}
In view of Lemma \ref{asympto2D0}, given any elements  $\upsigma^-$ and $\upsigma^+$ in $\Sigma$, we define the set
\begin{equation}
\label{mathbbW}
\mathbb W_\rL(\upsigma^-, \upsigma^+)=\left\{ v \in \mathbb W_\rL, {\rm \ such \ that  \  }  {\color{black}\underset {x_2 \to \pm \infty}\lim \rbI^\pm (v, x_2,\upsigma^\pm)=0}\right \}.
\end{equation}

\begin{lemma} 
\label{affinage}
Assume that $V$ satisfies  ${\color{black}(\mathbf{H}_{\mathbf{1}})}$, ${\color{black}(\mathbf{H}_{\mathbf{2}})}$ and  ${\color{black}(\mathbf{H}_{\mathbf{3}})}$. Given $\upsigma^-$ and $\upsigma^+$ in 
$\Sigma$, the set $\mathbb W_\rL(\upsigma^-, \upsigma^+)$ defined in \eqref{mathbbW} is an affine space. More precisely, we have the decomposition
\begin{equation}
\label{hilbertus}
\left\{
\begin{aligned}
\mathbb W_\rL(\upsigma^-, \upsigma^+)&=\{\bv_0\} + H_{\rm per}^1(\Lambda_\rL), {\rm \ where \ }\\
H_{\rm per}^1(\Lambda_\rL)&=\left\{u \in H^1(\Lambda_\rL), {\rm \ such \ that \ }  u({\color{black}-}\frac {\rL}{2}, x_2)=u(\frac {\rL}{2}, x_2),\forall {\color{black}x_2} \in \R \right\}, 
\end{aligned}
\right.
\end{equation}
and where $\bv_0: \Lambda_\rL\to \R^k $ is of the form
$\bv_0(x_1, x_2)=v_0(x_2)$, for $(x_1, x_2) \in \Lambda_\rL$, {\color{black}where} $v_0\in C^1(\R)$ is any function satisfying
\begin{equation}
\label{bv0}
\left\{
\begin{aligned}
v_0(s)&=\upsigma^+ {\rm \ for \ } s\geq 1\\
v_0(s)&=\upsigma^- {\rm \ for \ } s\leq -1.
\end{aligned}
\right. 
\end{equation}
\end{lemma}

The set $\mathbb W_\rL$ can therefore be decomposed as a disconnected and {\color{black}\emph{finite union}} of  \emph{affine spaces} 
 \begin{equation}
 \label{disconnected}
 \mathbb W_\rL=\underset{(\upsigma^+, \upsigma^-)\in \Sigma^2} \cup \mathbb W_\rL(\upsigma^-, \upsigma^+).
 \end{equation}

\subsubsection{Critical points of the energy}
Since  $\mathbb W_\rL \left(\upsigma^-, \upsigma^+\right)$ is an affine space, we may define critical points. 

\begin{proposition} 
  \label{differentiability}
   The functional $\E$ is $C^1$ on the affine space $\mbWL(\upsigma^-, \upsigma^+) = \{\bv_0\}+ H^1_{\rm per}(\Lambda_\rL)$, whose tangent space is 
   $H^1_{\rm per}(\Lambda_\rL)$. Moreover, we have
  \begin{equation}
  \label{weakasse}
  \rd \E(u) (w)=\int_{\Lambda_\rL } \nabla u\cdot \nabla w+  \nabla_u V (u){\color{black}\cdot} w, {\rm \ for \ any \ } w \in H^1_{\rm per}(\Lambda_\rL).
  \end{equation}
  \end{proposition}

  \begin{lemma} The critical points in $\mbWL(\upsigma^-, \upsigma^+)$ of $\E$ are exactly the  weak solutions to \eqref{themainthing}.
  \end{lemma}
  \begin{proof} If we take $w=\varphi$ with $\varphi \in \mathcal D(\mathring {\Lambda}_\rL)$, then \eqref{weaki}  yields
  \begin{equation}
  \label{diffe}
  \langle -\Delta u+  \nabla_u V (u), \varphi \rangle_{\mathcal D'(\mathring {\Lambda}_\rL), \mathcal D(\mathring {\Lambda}_\rL)}=0, \ \forall \varphi \in
  {\color{black}\mathcal D}(\mathring {\Lambda}_\rL).
  \end{equation}
Hence the result.
  \end{proof}

\begin{proposition} 
\label{regularita}
Assume that $V$ satisfies  ${\color{black}(\mathbf{H}_{\mathbf{1}})}$, ${\color{black}(\mathbf{H}_{\mathbf{2}})}$, and ${\color{black}(\mathbf{H}_{\mathbf{3}})}$. Any critical point of $\E$  belongs to $C^2(\Lambda_\rL)$ and verifies the system of partial differential equations
\eqref{themainthing}. Moreover the periodised map $\rbu$ defined in \eqref{periodisation} satisfies the system of partial differential equations
\begin{equation}
\label{Rdeux}
-\Delta \rbu+\nabla_u V(\rbu)=0 {\rm \ on \ } \R^2.
\end{equation}
\end{proposition}

  Among critical points, we may single out  the \emph{minimizing solutions}. Consider, for $i=1, 2$,  the two maps $\bw_i: \Lambda_\rL \to \R^k $  defined by 
\begin{equation}
\label{bwi}
\bw_i(x_1, x_2)=w_i(x_2),  \forall \, x_1\in [-\frac\rL 2, \frac \rL 2], \forall\,  x_2 \in \R, {\rm \ for \ } i=1, 2, 
\end{equation}
where $w_1$ and $w_2$ are as in Assumption ${\color{black}(\mathbf{H}_{\mathbf{4}})}$,   so that $\bw_i \in \mathbb W(\upsigma_0^-, \upsigma_0^+)$, for $i=1, 2$, and correspond to \emph{one-dimensional profiles}. 
Set $\bm:=\bm(\upsigma_0^-,\upsigma_0^+)$. We first have 
\begin{lemma}  We have $\bw_i \in \mathbb W_\rL(\upsigma_0^-, \upsigma_0^+)$, for $i=1, 2$  and the identities
\begin{equation}
\underset {v\in \mathbb W(\upsigma_0^-, \upsigma_0^+)} \inf \E(v, \Lambda_\rL)=\E(\bw_i, \Lambda_\rL)=\rL \bm, {\rm \ for \ }  i=1, 2. 
\end{equation}
Any minimizing solution $\bw$ in $\mathbb W(\upsigma_0^-, \upsigma_0^+)$  is of the form
$$
\bw(x_1, x_2)=w(x_2),  \forall x_1\in [-\frac\rL 2, \frac \rL 2], x_2 \in \R,
$$
{\color{black}where} $w: \R \to \R^k$ is such that $\E(w)=\bm$.
\end{lemma}
 
Since we seek solutions which are not one-dimensional profiles, we are hence led to the \emph{search for non-minimizing solutions}.

\subsubsection{The mountain-pass geometry}
 The core idea  in the proof is  to apply  the \emph{mountain-pass argument}  to the functional $\E$ on the set $\mathbb W(\upsigma_0^-, \upsigma_0^ +)$. For that purpose, we introduce the set of continuous paths joining the profiles $\bw_1$ and $\bw_2$, namely 
 \begin{equation}
 \label{paths}
 \mathcal P=\left \{p \in C^0\left({\color{black}[0, 1]}, \mathbb W_\rL(\upsigma_0^-, \upsigma_0^+)\right), {\rm \ such \ that \ } p(0)={\color{black}\bw_1} {\rm \ and \ } p(1)={\color{black}\bw_2}\right\}, 
 \end{equation}
 where $\bw_i$ is defined, for $i=1, 2$,  in \eqref{bwi}. We consider the number 
   \begin{equation}
  \label{mountainpass}
  \rbc_\rL= \underset{p \in \mathcal P}\inf  \left( \underset{s \in [0, 1]} \sup \E(p(s))  \right).
\end{equation}
The aim is to  show that $\rbc_\rL$ is a  \emph{generalized critical value of $\E$}, as defined next.     

\begin{definition} Let $c \in \R$. We say that $c$ is a \emph{generalized critical value of}  a $C^1$ functional $\E$  defined on a Banach space $X$ if and only if there exists a sequence $(u_n)_{n \in \N}$ such that 
\begin{equation}
\label{palais}
\left\{
\begin{aligned}
&\E(u_n) \to c, {\rm \ as \ } n \to + \infty, \\
&\rd\E(u_n) \to 0 {\rm  \ in  \ } X^\star.
\end{aligned}
\right.
\end{equation}
\end{definition}
 A sequence satisfying \eqref{palais} is called a \emph{Palais-Smale sequence at level c}. We say that problem  \eqref{mountainpass}  has the mountain-pass geometry if
  \begin{equation}
  \label{mountaingeo}
  \rbc_\rL  > \E(\bw_1)=\E(\bw_2)=\bm \rL.
  \end{equation}
  The mountain-pass Theorem then states that $\rbc_\rL$ is a generalized critical value of $\E$ provided \eqref{mountaingeo} holds.

  \subsubsection{Properties of $\rc_\rL$}
 We have the following upper-bound on $\rc_\rL$.
\begin{proposition}  
\label{Collevel}
 Assume that the potential $V$ satisfies assumptions ${\color{black}(\mathbf{H}_{\mathbf{1}})}$, ${\color{black}(\mathbf{H}_{\mathbf{2}})}$,  and   ${\color{black}(\mathbf{H}_{\mathbf{3}})}$.
There exists a constant $\rbC_{\rm pass}$ independent of $\rL$  such that, for any $\rL>0$
\begin{equation}
\label{passlevel}
\rbc_\rL \leq \bm\rL+ \rbC_{\rm pass}.
\end{equation}
\end{proposition}

We have also a lower bound:

\begin{proposition}
\label{lowerbound} 
Assume that the potential $V$ satisfies assumptions ${\color{black}(\mathbf{H}_{\mathbf{1}})}$, ${\color{black}(\mathbf{H}_{\mathbf{2}})}$,   ${\color{black}(\mathbf{H}_{\mathbf{3}})}$ and ${\color{black}(\mathbf{H}_{\mathbf{4}})}$. {\color{black}Let $\rL\geq1$.} Then there exists a constant $\rbC_{\rm sink}>0$ {\color{black}independent of} $\rL$ such that
\begin{equation}
\label{sinky}
\rbc_\rL \geq \bm\rL+ \rbC_{\rm sink}>\E(\mathbbmss w_i, \Lambda_\rL), {\rm \ for \ } i=1, 2.
\end{equation}
\end{proposition}

We are then in position to apply the mountain-pass theorem (see e.g. \cite{rabinowitz, bonafoux}), namely:
\begin{proposition}
  \label{genecritic}
   Assume that the potential $V$ satisfies assumptions ${\color{black}(\mathbf{H}_{\mathbf{1}})}$, ${\color{black}(\mathbf{H}_{\mathbf{2}})}$,   ${\color{black}(\mathbf{H}_{\mathbf{3}})}$ and ${\color{black}(\mathbf{H}_{\mathbf{4}})}$. {\color{black}Let $\rL\geq1$.} Then $\rbc_\rL$ is a generalized critical value of $\E$, that is there exists a  Palais-Smale sequence $(u_n)_{n \in \N}$ at level $\rbc_\rL$. 
  \end{proposition}

  \subsubsection{Asymptotic analysis  of Palais-Smale sequences}
In view of Proposition \ref{genecritic}, we are led to study Palais-Smale sequences for $\E$.
  
  \begin{proposition}
  \label{chainofsolutions}
  Assume that the potential $V$ satisfies assumptions ${\color{black}(\mathbf{H}_{\mathbf{1}})}$, ${\color{black}(\mathbf{H}_{\mathbf{2}})}$ and  ${\color{black}(\mathbf{H}_{\mathbf{3}})}$, and let $\rL>1$. Let  $\upsigma^-$ and  $\upsigma^+$ be {\color{black}two distinct elements} in $\Sigma$, and let 
$(u_n)_{n \in \N}$ in $\mathbb W_\rL(\upsigma^-, \upsigma^+)$ be a sequence satisfying \eqref{palais}.   Passing possibly to a subsequence, still {\color{black}denoted as} $(u_n)_{n \in \N}$,  there exists an integer $\mathfrak p_\star \in \N$, such that the following holds:  
\begin{itemize}
\item There exist  $\mathfrak p_\star$ sequences  of  real numbers  
$\displaystyle{(\rba_{n, k})_{n \in \N}}$   such that 
$\displaystyle{\rba_{n, 1}<\rba_{n, 2}<\ldots< \rba_{n, \mathfrak p_\star}}$ and  such that 
$\displaystyle{d_{n, k}\equiv \rba_{n, k+1}-\rba_{n, k}}$ tend to $+ \infty$, as $n \to + \infty$, for $k=1, \ldots, {\color{black}\mathfrak p_\star-1}$. 
\item  There exist  $\mathfrak p_\star$ non constant solutions $\rbu_k$, for $k={\color{black}1}, \ldots, \mathfrak p_\star$,  to \eqref{themainthing}, {\color{black}with endpoints satisfying \eqref{losange},} 
  \end{itemize}
  such that, setting $\rbe_2=(0, 1)$, we have the decomposition
  \begin{equation}
  \label{chainofthought}
  u_n=\underset {k=1}{\overset {\mathfrak p_\star} \sum} \rbu_k\left(\cdot-\rba_{n, k}\rbe_2\right){\color{black}-\sum_{k=1}^{\mathfrak p_\star-1}\upsigma_k^+}+ g_n, {\rm \ for \ any \ } n \in \N, 
  \end{equation}
  where {\color{black}the sum over intermediate wells is empty when $\mathfrak p_\star=1$,} $g_n \in \HL$ and $g_n \to 0$ in $H^1(\Lambda_\rL)$ as $n \to +\infty$.
 We have the identity 
 \begin{equation}
 \label{theidentity}
\underset {n \to +\infty}  \lim \E(u_n)=\underset {k=1}{\overset {\mathfrak p_\star} \sum}\,  \E(\rbu_k).
\end{equation}
  \end{proposition}
  \begin{remark}{\rm Similar  or related  results appear  in many places in the study of non linear PDE's with lack of compactness. For instance, it is  very close to the results for the corresponding ODE. We refer to \cite{rabinowitz, bonafoux} for references.
  }
  \end{remark}
  
{\color{black}The profiles have matching endpoints, as proved independently below:}
\begingroup\color{black}
\begin{equation}
\label{losange}
\left\{
\begin{aligned}
\rbu_k&\in\mathbb W_\rL(\upsigma_k^-,\upsigma_k^+),\quad\upsigma_k^\pm\in\Sigma,
\quad k=1,\ldots,\mathfrak p_\star,\\
\upsigma_k^+&=\upsigma_{k+1}^-,\quad k=1,\ldots,\mathfrak p_\star-1,\\
\upsigma_1^-&=\upsigma^-,\qquad\upsigma_{\mathfrak p_\star}^+=\upsigma^+.
\end{aligned}\right.
\end{equation}
\endgroup

 We have the following
 
 \begin{proposition}
 \label{nicht} 
    Assume that the potential $V$ satisfies assumptions ${\color{black}(\mathbf{H}_{\mathbf{1}})}$, ${\color{black}(\mathbf{H}_{\mathbf{2}})}$, ${\color{black}(\mathbf{H}_{\mathbf{3}})}$ and ${\color{black}(\mathbf{H}_{\mathbf{4}})}$. There exists a constant $\rL_\star>1$, depending only on $V$, such that, for $\rL\geq\rL_\star$, the decomposition in Proposition \ref{chainofsolutions} of any Palais-Smale sequence in $\mathbb W_\rL(\upsigma_0^-,\upsigma_0^+)$ at level $\rbc_\rL$ satisfies 
 $$
 \upsigma_k^-\not = \upsigma_k^+, {\rm \ for \ } k=1, \ldots, \mathfrak p_\star.
 $$
 \end{proposition}

 \smallskip
  
    We will also rely on the following:
    
    \begin{proposition}
    \label{lwb0}
  Assume that the potential $V$ satisfies assumptions ${\color{black}(\mathbf{H}_{\mathbf{1}})}$, ${\color{black}(\mathbf{H}_{\mathbf{2}})}$, ${\color{black}(\mathbf{H}_{\mathbf{3}})}$, and let $\rL>1$. There exists a constant $\upeta_{\rm min}$, depending only on $V$, such that any  non-constant solution $\rbu$ to \eqref{themainthing} satisfies the lower bound
    \begin{equation}
    \label{lwb1}
    \E(\rbu, \Lambda_\rL) \geq \upeta_{\rm min} .
    \end{equation}
    \end{proposition}
   
 \subsection{Some properties of one-dimensional  profiles satisfying ${\color{black}(\mathbf{H}_{\mathbf{5}})}$}
   Proposition \ref{genecritic} and Proposition \ref{chainofsolutions}  have provided $\mathfrak p_\star$ solutions of \eqref{themainthing}. We have to check that at least one of these solutions is  not a one-dimensional profile. To that aim, we rely on properties of one-dimensional minimizing solutions, in particular on the spectral property ${\color{black}(\mathbf{H}_{\mathbf{5}})}$.

\begin{proposition}
\label{gapounet}
 Assume that ${\color{black}(\mathbf{H}_{\mathbf{1}})}$, ${\color{black}(\mathbf{H}_{\mathbf{2}})}$, ${\color{black}(\mathbf{H}_{\mathbf{3}})}$, ${\color{black}(\mathbf{H}_{\mathbf{4}})}$, and ${\color{black}(\mathbf{H}_{\mathbf{5}})}$ hold. Then, there exists some constant $\uprho_\star>0$, such that, if $\rbq\in \mathbb X(\upsigma_0^-, \upsigma_0^+)$ is a solution to \eqref{heterocline}, satisfying  the energy bound
 \begin{equation}
 \label{gapos}
 \E(\rbq, \R) \leq \bm(\upsigma_0^-, \upsigma_0^+)+\uprho_\star,
 \end{equation}
 then   we have $\rbq(\cdot)=w_i(\cdot-a)$, for some $i\in \{1, 2\}$  and   some $a \in \R$.
\end{proposition}
  \subsection{Outline of the  completion of the proof of Theorem \ref{maintheo} }
   \label{qed}
  At this stage, we are given an integer $\mathfrak p_\star$, and a family of non constant  solutions  $\rbu_k$ to \eqref{themainthing} such that 
  \begin{equation}
  \label{mainidentity}
  \bm\rL  + \rbC_{\rm  sink} \leq  \underset{k=1}{\overset {\mathfrak p_\star} \sum } \E(\rbu_k)=\rbc_\rL \leq \bm \rL+ \rbC_{\rm pass}.
  \end{equation}
  The proof of the  Theorem is deduced from the relations \eqref{mainidentity}. {\color{black}It} requires two important ingredients. The first one allows us to exclude the possibility of homoclinic solutions,  and is stated in Proposition {\color{black}\ref{mono}}. The second one relies on a spectral property ${\color{black}(\mathbf{H}_{\mathbf{5}})}$, and actually deals with one-dimensional profiles. It states that non-minimizing solutions to \eqref{heterocline} cannot be arbitrarily close to  the orbits (by the action of translations) of the minimizing solutions $w_1$ and $w_2$. 
\subsection{Outline of the paper}

  We gather in the next three sections some background material, mostly some well-known facts that we adapt to our problem and our {\color{black}assumptions.} More specifically,  the next section, Section \ref{functframe} gathers a number of elementary facts on the functional {\color{black}framework, in} particular on the affine space structure (well-known, at least in the one-dimensional case).   We then present in {\color{black}Sections} \ref{propsol} and 4 some general properties of solutions and approximate solutions, which will then be used  in Section \ref{Passy}, in the description of the behavior of general Palais-Smale sequences, in particular the proof of Proposition \ref{chainofsolutions}.  Section \ref{one-d} is devoted to properties of one-dimensional solutions.  The actual proof of Theorem \ref{maintheo} starts in Section \ref{geometry}, where the mountain-pass geometry is established (i.e. the proofs of Propositions \ref{Collevel} 	and \ref{lowerbound}).  The proof of Theorem \ref{maintheo} is then completed in Section \ref{final}.  In the Appendix, we provide an example of an explicit potential $\rbV$ satisfying all our assumptions 
  ${\color{black}(\mathbf{H}_{\mathbf{1}})}$, ${\color{black}(\mathbf{H}_{\mathbf{2}})}$, ${\color{black}(\mathbf{H}_{\mathbf{3}})}$, ${\color{black}(\mathbf{H}_{\mathbf{4}})}$, and ${\color{black}(\mathbf{H}_{\mathbf{5}})}$.
    \tableofcontents

   \numberwithin{theorem}{section} \numberwithin{lemma}{section}
\numberwithin{proposition}{section} \numberwithin{remark}{section}
\numberwithin{corollary}{section}
\numberwithin{equation}{section}

  \section{The functional framework}
  \label{functframe}
  In this section, we  provide the proof of Lemma \ref{affinage} establishing the affine structure of $\mathbb W(\upsigma^-, \upsigma^+)$ and  also  specify the regularity properties of the energy functional $\E$.
  
\subsection{Properties of the potential $V$}  
We introduce parameters controlling the behavior of $V$ in the vicinity of the set 
$\Sigma$ as well as far from this set.  First, we introduce some  $\upmu_0>0$ sufficiently small so that
    \begin{equation}
    \label{kiv}
    \left\{
    \begin{aligned}
 \B^k(\upsigma_i, 2\upmu_0)\cap \B^k(\upsigma_j, 2\upmu_0) &= \emptyset,
 {\rm \ for \  all  \ } i\neq j  {\rm \  in \ }  \{1,\cdots,\mathfrak q_\star\} {\rm \ and \  such  \   that  \ }  \\
 \frac{1}{2}\lambda_i^-{\rm Id} \leq  
 \nabla^2 V (y) &\leq 2\lambda_i^+{\rm Id},    \  \ 
  {\rm \ for \  all  \ } i\in \{1,\cdots,\mathfrak q_\star\} {\rm \  and  \ }  y \in \B(\upsigma_i,2\upmu_0).
 \end{aligned}
  \right. 
\end{equation}
The existence of $\upmu_0$ follows from assumption  ${\color{black}(\mathbf{H}_{\mathbf{3}})}$. As a first consequence, we have:
\begin{proposition} [see e.g. \cite{bethu}]
\label{potto}
 For any $i=1, \ldots, \mathfrak q_\star$  and any   $y \in \B^k(\upsigma_i, 2\upmu_0)$, we have the local bound
\begin{equation}
\label{glutz}
\left\{
\begin{aligned}
\frac{1}{4}\lambda_i^-\vert y-\upsigma_i \vert^2 & \leq V(y) \leq  \lambda_i^+\vert y-\upsigma_i \vert^2, 
  \\
  \frac{1}{2}\lambda_i^-\vert y-\upsigma_i \vert^2 & \leq  \nabla V(y) \cdot (y-\upsigma_i) \leq  2\lambda_i^+\vert y-\upsigma_i \vert^2.
    \end{aligned}
  \right.
 \end{equation}   
Choosing possibly an even smaller constant $\upmu_0$, we  have
   \begin{equation}
   \label{extrut}
   V(y) \geq  \upalpha_0\equiv{\color{black}\frac{1}{64}} \lambda_0 \upmu_0^2>0 {\rm \ on \ } 
   \R^k  \setminus  \underset{i=1}{\overset {\mathfrak q_\star} \cup} \B^k(\upsigma_i, \frac{\upmu_0}{4}), 
  \end{equation}
  where we  set  
  \begin{equation}
  \label{extrut2}
  {\color{black}\lambda_0=\lambda_0^-:=\min\{\lambda_i^-, i=1, \ldots, \mathfrak q_\star\}.}
  \end{equation}
    \end{proposition}

   For the behavior far from $\Sigma$, we introduce the number
   \begin{equation}
\label{grandr0}
{\rm \bf R}_0=\sup\{\vert \upsigma  \vert, \upsigma \in \Sigma\},
\end{equation}     
 It can be shown, in view of assumption ${\color{black}(\mathbf{H}_{\mathbf{2}})}$,  that there exists a constant $\upbeta_\infty >0$ such that 
\begin{equation}
\label{upbetainfty}
V(y) \geq \upbeta_\infty \vert y \vert^2,   {\rm \ for \ any \ } y {\rm \ such \ that \ }  \vert y \vert  \geq  2 {\rm \bf R}_0.  
\end{equation}

 The following technical Lemma will also be useful.
 
\begin{lemma} 
 \label{technical}
 Assume $V$ satisfies assumptions  ${\color{black}(\mathbf{H}_{\mathbf{1}})}$, ${\color{black}(\mathbf{H}_{\mathbf{2}})}$ and  ${\color{black}(\mathbf{H}_{\mathbf{3}})}$. Then,  there exists a constant $\rK_0>0$  depending only on $V$, such that, for any $\upsigma \in \Sigma$, and any $y \in \R^k$, we have the bounds
 \begin{equation}
 \label{technical0}
 \left\{
 \begin{aligned}
 0\leq V(\upsigma +y)&\leq \rK_0 \left(\vert y \vert^2+ \vert y \vert^m\right),  {\rm \ and  \ } \\
  \vert \nabla V(\upsigma+y)\vert &\leq \rK_0\left(\vert y \vert+ \vert y \vert^m\right).
 \end{aligned}
 \right.
 \end{equation}
 \end{lemma}
 \begin{proof}  Let 
 \begin{equation}
 \label{lambda0lus}
 \lambda_0^+=\lambda_{\rm max}^+={\rm Max} \{\lambda_i^+, i=1, \ldots, \mathfrak q_\star\},
 \end{equation}
  where $\lambda_i^+$ is defined in \eqref{kiv}. We deduce from \eqref{kiv}  
 \begin{equation}
 \label{uno}
 V(y+\sigma_i)\leq \lambda_0^+ \vert y \vert^2, {\rm \ for \ } y \in  \B^k ({\color{black}0}, 2 \upmu_0), {\rm \ for \ any \ }
  i=1, \ldots \mathfrak q_\star. 
 \end{equation}
 On the other hand, by  assumption ${\color{black}(\mathbf{H}_{\mathbf{2}})}$, we have also 
 \begin{equation}
 \label{due}
  V(y) \leq C_m(\vert y \vert^m+1) = \frac{C_m}{\upmu_0^m} \left( \upmu_0^m \vert y \vert^m+ \upmu_0^m\right), \forall   y \in \R^k.
  \end{equation}
 {\color{black}Enlarging $C_m$ if necessary, the same growth bound holds for $V(\upsigma+y)$ and $\vert\nabla V(\upsigma+y)\vert$, uniformly for $\upsigma\in\Sigma$, since $\Sigma$ is finite.}
 We claim that
 \begin{equation}
 \label{claim}
 V(y+\upsigma) \leq\lambda_0^+ \vert y \vert^2+\frac{C_m(\upmu_0^m+1)}{\upmu_0^m}   \vert y \vert^m,  {\rm \ for \ any \ } y \in \R^k.
\end{equation}
 Indeed, if $\vert y \vert \leq \upmu_0$, we apply directly \eqref{uno}. Otherwise $\vert y \vert > \upmu_0$: in this case 
 $ \upmu_0^m \leq \vert y \vert^m$, so that, in view of \eqref{due}
 \begin{equation}
 \label{patate}
V(\upsigma + y) \leq \frac{C_m(\upmu_0^m+1)}{\upmu_0^m}   \vert y \vert^m {\rm \ for \ } \vert y \vert \geq \upmu_0.
\end{equation}
 \eqref{claim} holds  likewise. Setting  
 $\displaystyle{ \rK_0={\rm Max} \left\{{\color{black}2}\lambda_0^+, \frac{C_m(\upmu_0^m+1)}{\upmu_0^m}\right\}}$, the first inequality in \eqref{technical0} follows. For the second inequality in \eqref{technical0},  we observe similarly  that
 it follows from \eqref{kiv} that
 \begin{equation}
 \label{unodue}
\vert  \nabla V(y+\sigma_i)\vert \leq 2\lambda_0^+ \vert y \vert, {\rm \ for \ } y \in  \B^k ({\color{black}0}, 2 \upmu_0), {\rm \ for \ any \ } i=1, \ldots \mathfrak q_\star. 
 \end{equation}
 whereas  by  assumption ${\color{black}(\mathbf{H}_{\mathbf{2}})}$, we have, as above 
$$\vert  \nabla V(\upsigma_i +y)\vert  \leq C_m(\vert y \vert^m+1) = \frac{C_m}{\upmu_0^m} \left( \upmu_0^m \vert y \vert^m+ \upmu_0^m\right), 
$$
 If $\vert y \vert \leq \upmu_0$, we apply directly \eqref{unodue}. Otherwise, we argue as above to assert that 
 \begin{equation}
 \label{luistrio}
 \begin{aligned}
 \vert  \nabla V(y+\sigma_i)\vert &\leq \frac{C_m}{\upmu_0^m} \left( \upmu_0^m \vert y \vert^m + \vert y \vert^m\right)
 =\frac{C_m (\upmu_0^m+1)}{\upmu_0^m} \left(\vert y \vert^m \right), {\rm \ for \ } \vert y \vert \geq \upmu_0. \\
 &\leq \rbK_0  \vert y \vert^m.
\end{aligned}
 \end{equation}
and we conclude similarly.
\end{proof}

\subsection{On the $L^2$ distance to the vacuum set $\Sigma$}  

The main result in this section is:

\begin{lemma} 
\label{uniclo}
Assume that $V$ satisfies  ${\color{black}(\mathbf{H}_{\mathbf{1}})}$, ${\color{black}(\mathbf{H}_{\mathbf{2}})}$ and  ${\color{black}(\mathbf{H}_{\mathbf{3}})}$. Let $v: {\rbQ}=[0, 1]^2 \to \R^k$ be such that $\E(v,\rbQ )< +\infty$.  Then, there exists $\upsigma_{\rbQ}\in \Sigma$ such that  we have the bound
\begin{equation}
\label{unif2D}
\int_\rbQ \vert v(x)-\upsigma_{\rbQ}\vert^2 \rd x \leq \rbC_2  \E(v, \rbQ), 
\end{equation}
where   the constant {\color{black}$\rbC_2>0$} depends only on the constants $\lambda_0$, $\upalpha_0$, ${\rm \bf R_0}$,  and $\upbeta_\infty$.
\end{lemma}
The proof relies on a corresponding related \emph{one-dimensional uniform estimate} proven in \cite{bethu}, Lemma 2.4:

     \begin{lemma} [\cite{bethu}]
     \label{valli}
   Assume that $V$ satisfies  ${\color{black}(\mathbf{H}_{\mathbf{1}})}$, ${\color{black}(\mathbf{H}_{\mathbf{2}})}$ and  ${\color{black}(\mathbf{H}_{\mathbf{3}})}$.  There exists a constant ${\rm C}_{\rm unf}>0$, depending only on $V$,  such that, for any   given  map    $u: [0, 1] \to \R^k$,  there exists an element $\upsigmam \in \Sigma$ such that 
      \begin{equation}
      \label{bornunif}
     \vert u(\ell)-\upsigmam \vert^2    \leq  {\rm C}^2_{\rm unf}\int_{[0,1]} e(u)(s){\rm d} s, {\rm \ for \ } \ell \in [0, 1].  
      \end{equation}
       \end{lemma}

 \begin{remark}
 {\rm In \cite{bethu}, the estimates are established on circles, but the proof is unchanged if we replace circles by intervals. The constant 
${\rm C}^2_{\rm unf}$ may actually  be estimated explicitly.   For that purpose, we introduce the constant. 
Then, we set 
  \begin{equation}
   \label{extruite}
     \upalpha_1=\inf \{ \sqrt{\lambda_0}\frac{\upmu_0^2}{64},\frac{\pi\upmu_0\lambda_0}{16},  \upalpha_0\}.
  \end{equation}
where $\upalpha_0$ and $\lambda_0$ are defined in \eqref{extrut}.
The constant $C_{\rm unf}$ is then given by
\begin{equation}
\label{defcunf}
 {\rm C}_{\rm unf}^2  ={\rm Max} 
 \left\{
 8\sqrt{\lambda_0^{-1}}+4\lambda_0^{-1},
 \left(\frac{8}{\sqrt{\upbeta_\infty}} +
    \frac{10{\rm \bf  R}_0^2 }{\upalpha_1}
    \right), 
 \frac{2}{ \upbeta_\infty  } \left(  1+ \left(\frac{2{\rm \bf R}_0^2 \upbeta_\infty}{\upalpha_1}\right) \right)
 \right\}.
\end{equation}
It depends hence only on the numbers $\upmu_0, \lambda_0, {\rm \bf R_0}$,  and $\upbeta_\infty$.
}
\end{remark}

\begin{proof}[Proof of Lemma \ref{uniclo} completed] We  apply first  Lemma \ref{valli} in the horizontal direction, on each of the intervals $[0, 1]\times\{x_2\}$, for $x_2\in [0,1]$. This yields, for every $x_2\in [0, 1]$,  an element $\upsigma(x_2) \in \Sigma$ such that,
 \begin{equation}
      \label{bornunifx2}
     \vert u(\ell, x_2)-\upsigma(x_2) \vert^2    \leq  {\rm C}^2_{\rm unf}\int_{[0,1]\times \{x_2\}} e(u)(s, x_2){\rm d} s, {\rm \ for \ } \ell \in [0, 1].  
      \end{equation}
      Integrating, for fixed $x_2$,  on the segment $\rbI(x_2)=[0, 1] \times \{x_2\}$, we obtain
      \begin{equation}
      \label{intx10}
      \int_{[0, 1] \times \{x_2\}}  \vert u(\ell, x_2)-\upsigma(x_2) \vert^2  \rd \ell \leq     
      {\rm C}^2_{\rm unf}\int_{[0,1]\times \{x_2\}} e(u)(\ell, x_2){\rm d} \ell.
      \end{equation}
Next, we work in the vertical direction: By a mean value argument, we may find some $s_0\in [0, 1]$ such that 
\begin{equation}
\int_{\{s_0\}\times [0,1]} e(u)(s_0, s){\rm d} s \leq \E(u, \rbQ).
\end{equation}It follows from Lemma \ref{valli} and \eqref{bornunifx2} that there exists some element $\upsigma_0\in \Sigma$ such that 
\begin{equation}
\label{gepeto}
\vert u(s_0, x_2)-\upsigma_0\vert^2 \leq  {\rm C}^2_{\rm unf} \int_{\{s_0\}\times[0,1]} e(u)(s_0, s){\rm d} s \leq {\rm C}^2_{\rm unf}\E(u, \rbQ), 
{\rm \ for \ any \ } x_2\in [0,1].
\end{equation}
Combining \eqref{bornunifx2},  for $\ell=s_0$ with \eqref{gepeto}, we are led, for any $x_2\in [0, 1]$
\begin{equation}
\label{gepeto2}
\begin{aligned}
\vert \upsigma_0-\upsigma(x_2)\vert^2  &\leq 2 \left(\vert u(s_0, x_2)-\upsigma(x_2)\vert^2+\vert \upsigma_0-u(s_0,x_2)\vert^2\right)\\
& \leq 2{\rm C}^2_{\rm unf}\left( \E(u, \rbQ)+\int_{[0,1]\times \{x_2\}} e(u)(s, x_2){\rm d} s\right).
\end{aligned}
\end{equation}
We will also invoke in some places  the obvious bound, $\rbR_0$ being defined in \eqref{grandr0}, 
\begin{equation}
\label{bebeto}
\vert \upsigma-\upsigma_0\vert^2\leq 4\rbR_0^2, {\rm \ for \ any \  } \upsigma \in \Sigma.
\end{equation}
Next, we distinguish several cases

\smallskip
\noindent
{\it Case 1}. $\displaystyle{\E(u, \rbQ) \geq \frac{{\color{black}\upmu_0^2}}{16{\rm C}^2_{\rm unf}}}$. In this case, integrating \eqref{intx10} with respect to $x_2$ on $[0,1]$ we obtain
\begin{equation}
\Vert u-\upsigma(x_2) \Vert_{L^2({\color{black}\rbQ})}^2 
=\int_{[0, 1]^2} \vert u(x_1, x_2)-\upsigma(x_2) \vert^2\rd x_1\rd x_2\leq  {\rm C}^2_{\rm unf}\E(u, \rbQ).
\end{equation}
We deduce, combining with \eqref{bebeto},  that 
\begin{equation}
\label{bebeto2}
\begin{aligned}
\Vert u-\upsigma_0 \Vert_{L^2({\color{black}\rbQ})}^2&\leq 2\left( \Vert u-\upsigma(x_2) \Vert_{L^2({\color{black}\rbQ})}^2+\Vert {\color{black}\upsigma_0}-\upsigma(x_2) \Vert_{L^2({\color{black}\rbQ})}^2\right)\\
&\leq 2{\rm C}^2_{\rm unf}\E(u, \rbQ)+ 8\rbR_0^2 \leq 2{\rm C}^2_{\rm unf}\E(u, \rbQ)+(\frac{128 {\rm C}^2_{\rm unf} {\color{black}\rbR_0^2}}{{\color{black}\upmu_0^2}})\frac{{\color{black}\upmu_0^2}}{16 {\rm C}^2_{\rm unf}}\\
&\leq \left( 2{\rm C}^2_{\rm unf}+\frac{128 {\rm C}^2_{\rm unf} {\color{black}\rbR_0^2}}{{\color{black}\upmu_0^2}}\right)\E(u, \rbQ).
\end{aligned}
\end{equation}
Hence \eqref{unif2D} is verified with $\upsigma_\rbQ=\upsigma_0$, provided we choose $\rbC_2$ so that 
\begin{equation}
\label{grantorino}
\rbC_2 \geq  2{\rm C}^2_{\rm unf}+\frac{128{\rm C}^2_{\rm unf}{\color{black}\rbR_0^2}}{{\color{black}\upmu_0^2}}.
\end{equation}

\medskip
\noindent
{\it Case 2}. $\displaystyle{\E(u, \rbQ) < \frac{{\color{black}\upmu_0^2}}{16{\rm C}^2_{\rm unf}}}$.  
Going back to \eqref{gepeto2}, we are led, for every $x_2\in [0, 1]$,  to 
\begin{equation}
\label{gepeto3}
\begin{aligned}
\vert \upsigma_0-\upsigma(x_2)\vert^2  
&\leq \frac{{\color{black}\upmu_0^2}}{8}+2{\rm C}^2_{\rm unf}\int_{[0,1]\times \{x_2\}} e(u)(s, x_2){\rm d} s=
\frac{{\color{black}\upmu_0^2}}{8}+2{\rm C}^2_{\rm unf} \rbJ_1 (u, x_2),  \\
 {\rm where \ } &\rbJ_1(u, x_2)=\int_{[0,1]\times \{x_2\}} e(u)(s, x_2){\rm d} s.
\end{aligned}
\end{equation}
We distinguish once more two subcases for $x_2$:

\smallskip
\noindent
{\it Subcase 2a}. $\displaystyle{\rbJ_1(u, x_2) \leq \frac{{\color{black}\upmu_0^2}}{16{\rm C}^2_{\rm unf}}}$.  In view of  \eqref{gepeto3}, we obtain 
\begin{equation}
\vert \upsigma_0-\upsigma(x_2)\vert^2  \leq \frac{{\color{black}\upmu_0^2}}{4},
\end{equation}
so that, in view of the definition \eqref{kiv} of $\upmu_0$, we deduce that $\upsigma(x_2)=\upsigma_0$.  Going back to \eqref{intx10}, we are led to
  \begin{equation}
      \label{intx1}
      \int_{[0, 1] \times \{x_2\}}  \vert u(\ell, x_2)-\upsigma_0 \vert^2  \rd \ell \leq     {\rm C}^2_{\rm unf}\int_{[0,1]\times \{x_2\}} e(u)(s, x_2){\rm d} s.
      \end{equation}
      
  \smallskip
\noindent
{\it Subcase 2b}. $\displaystyle{\rbJ_1(u, x_2) \geq \frac{{\color{black}\upmu_0^2}}{16{\rm C}^2_{\rm unf}}}$.  We argue as in case 1:
\begin{equation}
\label{pinoc}
\begin{aligned}
\  \int_{[0, 1] \times \{x_2\}}  \vert u(\ell, x_2)-\upsigma_0 \vert^2  \rd \ell &\leq 2\left(  \int_{[0, 1] \times \{x_2\}}  \vert u(\ell, x_2)-\upsigma(x_2) \vert^2  \rd \ell+
  \vert \upsigma_0-{\color{black}\upsigma(x_2)}\vert^2\right) \\
  &\leq 2\left( {\rm C}^2_{\rm unf}\rbJ_1(u, x_2)+4\rbR_0^2\right)\\
  &=2 {\rm C}^2_{\rm unf}\rbJ_1(u, x_2)+
  \frac{{\color{black}128}\rbR_0^2{\rm C}^2_{\rm unf}}{{\color{black}\upmu_0^2}}\frac{{\color{black}\upmu_0^2}}{16{\rm C}^2_{\rm unf}} \\
  &\leq \left( 2{\rm C}^2_{\rm unf}+ \frac{{\color{black}128}\rbR_0^2{\rm C}^2_{\rm unf}}{{\color{black}\upmu_0^2}}\right)\rbJ_1(u, x_2)
\end{aligned}
\end{equation}
We observe, that in both subcases 2a and 2b, if $\rbC_2$ satisfies condition \eqref{grantorino}, then, we have 
\begin{equation}
\label{granturismo}
 \int_{[0, 1] \times \{x_2\}}  \vert u(\ell, x_2)-\upsigma_0 \vert^2  \rd \ell \leq \rbC_2\rbJ_1(u, {\color{black}x_2}).
\end{equation}
Integrating with respect to $x_2$, we finally obtain
$$
\Vert u-\upsigma_0\Vert^2_{L^2(\rbQ)} \leq \rbC_2 \E(u, \rbQ), 
$$
hence,  inequality \eqref{unif2D} is  established with $\upsigma_\rbQ=\upsigma_0$, both in case 1 and Case 2, so that the proof is complete.
\end{proof}

  Lemma \ref{uniclo} does not necessarily provide uniqueness of $\upsigma_\rbQ$. This is however true when the energy 
$\E(u, \rbQ)$ is sufficiently small.
  
\begin{lemma}
\label{petitement}
Assume that $V$ satisfies  ${\color{black}(\mathbf{H}_{\mathbf{1}})}$, ${\color{black}(\mathbf{H}_{\mathbf{2}})}$ and  ${\color{black}(\mathbf{H}_{\mathbf{3}})}$. There exists a constant ${\color{black}\upeta_0}> 0$ depending only on $V$
such that, if $\E(u, \rbQ)\leq {\color{black}\upeta_0}$ then there exists a unique element $\upsigma_\rbQ\in \Sigma$ verifying inequality \eqref{unif2D}.
\end{lemma}  
\begin{proof} Assume by contradiction that there are two such elements $\upsigma_1$ and $\upsigma_2$. 
We have, {\color{black}thanks} to \eqref{unif2D}
\begin{equation*}
\begin{aligned}
\vert \upsigma_1-\upsigma_2\vert^2&\leq 2\left(\int_\rbQ \vert u(x)-\upsigma_1\vert^2 \rd x +\int_\rbQ \vert u(x)-\upsigma_2\vert^2 \rd x\right) \\
&\leq  4\rbC_2 \E(u, \rbQ) \leq 4 \rbC_2{\color{black}\upeta_0}.
\end{aligned}
\end{equation*}
Set 
\begin{equation}
\label{infdist}
\upgamma_0=\inf\{ \vert \upsigma'- \upsigma \vert^2, {\rm \ for \ } \upsigma'\not= \upsigma, {\rm \ in \  } \Sigma\}.
\end{equation}
 If we choose  $4 \rbC_2{\color{black}\upeta_0}<\upgamma_0$, we obtain
$\vert \upsigma_1-\upsigma_2\vert^2<\upgamma_0$, a contradiction.
\end{proof}

  \subsection{Adjacent and overlapping cubes}
 \label{adjacent-squares}
  
  We consider here a situation with two adjacent  or overlapping  squares, namely, for $0\leq r \leq 1$  the squares $\rbQ^+$ and $\rbQ_r^-$ defined by
  \begin{equation*}
  \left\{
  \begin{aligned}
  &\rbQ_r^-=[-1+r,r]\times [0, 1], \ \rbQ^+= [0,1]^2, {\rm \ so \ that \ } \\
  &\Theta_r\equiv \rbQ_r^-\cap \rbQ^+= [0, r]\times [0, 1]
  {\rm \ and \ }   \rbK_r=\rbQ^+\cup \rbQ_r^-=[-1+r,  1] \times [0, 1], 
  \end{aligned}
  \right.
  \end{equation*}
  so that $\rbQ_r^-$ and $\rbQ^+$ overlap when $0<r\leq1$, and are adjacent for $r=0$.
  
  \begin{lemma}
  \label{adjacent}
  Assume that $V$ satisfies  ${\color{black}(\mathbf{H}_{\mathbf{1}})}$, ${\color{black}(\mathbf{H}_{\mathbf{2}})}$ and  ${\color{black}(\mathbf{H}_{\mathbf{3}})}$.
  Let $\upsigma_r^-$ and $\upsigma^+$ be elements of $\Sigma$ provided by Lemma \ref{uniclo} on $\rbQ_r^-$ and $\rbQ^+$ respectively. There exists a constant $\upeta_1>0$ depending only on the potential $V$ such that, if 
  \begin{equation}
  \label{adjacebt0}
  \E(u, \rbK_r) \leq 2\upeta_1, 
  \end{equation}
 then $\upsigma^+=\upsigma_r^-$.
 \end{lemma}
  \begin{proof} 
   We may impose that ${\color{black}2\upeta_1\leq\upeta_0}$, where $\upeta_0$ is defined in Lemma \ref{petitement}, so that $\upsigma_r^-$ and $\upsigma^+$ are uniquely defined in view of \eqref{adjacebt0}.   We distinguish two cases:
  
  \smallskip
  \noindent
  {\it Case 1: $\displaystyle{r\geq \frac14}$, large overlap}.
  We  have, in view of \eqref{unif2D}, the fact that $\vert \Theta_r\vert=r$  and the triangle inequality
  \begin{equation*}
\begin{aligned}
\frac14 \vert \upsigma_r^-  -\upsigma^+\vert^2\leq r\vert \upsigma_r^--\upsigma^+\vert^2&\leq 2\left(\int_{\Theta_r} \vert u(x)-\upsigma_r^-\vert^2 \rd x +
\int_{\Theta_r} \vert u(x)-\upsigma^+\vert^2 \rd x\right) \\
&\leq  2\left(\int_{\rbQ_r^-} \vert u(x)-\upsigma_r^-\vert^2 \rd x +
\int_{\rbQ^+} \vert u(x)-\upsigma^+\vert^2 \rd x\right) \\
&\leq 2 \rbC_2\left( \E(u, \rbQ_r^-)+\E(u, \rbQ^+) \right)\leq {\color{black}8}\rbC_2 \upeta_1.
\end{aligned}
\end{equation*}
  Next, we choose $\upeta_1$ so that $32\rbC_2 \upeta_1<\upgamma_0$, where $\upgamma_0$ is defined in \eqref{infdist}. {\color{black}It} follows that
  {\color{black}$\vert \upsigma_r^--\upsigma^+\vert^2<\upgamma_0$,} and hence, in view of the definition of $\upgamma_0$ we deduce that $\upsigma^+=\upsigma_r^-$. 

\smallskip
  \noindent
  {\it Case 2: $\displaystyle{0\leq r\leq \frac14}$, small overlap}.     We  consider the cube $\rbQ_{1\slash2}^-=[-1\slash 2, 1\slash2]\times[0, 1]$, and the element $\upsigma_{1\slash2}\in \Sigma$ provided by Lemma \ref{uniclo} for the restriction of $u$ to $\rbQ_{1\slash 2}^-$. We are in position to apply the result of case $1$,  so that we obtain $\upsigma_{1\slash 2}=\upsigma^+$.  Similarly, we are in position to apply case $1$ to the restriction of the  squares $\rbQ_r^-$ and $\rbQ_{1\slash2}^-$, since their overlap is sufficiently large. Hence, we deduce $\upsigma_{1\slash 2}=\upsigma_r^-${\color{black}.} Combining the previous identities, we are led to $\upsigma_r^-=\upsigma^+$.
 \end{proof}
  
  For $a \in \R\cup\{-\infty\}$, $b \in \R \cup \{+\infty\}$  such that $b-a>0$, we  introduce the set 
  \begin{equation}
  \label{Gamma}
  \Gamma_\rL (a,  b)= [-\frac{\rL}{2}, \frac {\rL}{2}] \times[a, b] \subset \Lambda_\rL.
  \end{equation}

   \begin{corollary} 
   \label{corollo}
  Assume that $V$ satisfies  ${\color{black}(\mathbf{H}_{\mathbf{1}})}$, ${\color{black}(\mathbf{H}_{\mathbf{2}})}$ and  ${\color{black}(\mathbf{H}_{\mathbf{3}})}$ and assume that $\rL\geq 1$. Let $u\in H^1_{\rm loc }(\Lambda_\rL)$  be such that  $\E(u, \Lambda_\rL)<+\infty$. Let $a \in \R\cup\{-\infty\}$, $b \in \R \cup \{+\infty\}$  such that $b-a\geq 1$. 
  If we have, for   any square $\rbQ(x_1, x_2)\equiv  [x_1, x_1+1]\times[x_2, x_2+1]\subset \Gamma_\rL(a,b)$, the inequality
    \begin{equation}
    \label{petitou0}
  \E\left(u, \rbQ(x_1, x_2) \right)\leq \upeta_1,    
  \end{equation}                                                 
 then, there exists a unique $\upsigma^*\in \Sigma$ such that 
 \begin{equation}
 \label{majorana1}
 \int_{\rbQ(x_1, x_2)}\vert u-\upsigma^*\vert^2\leq \rbC_2\E\left(u, \rbQ(x_1, x_2)\right).
 \end{equation}
  \end{corollary}
 The proof is an immediate consequence of Lemma \ref{adjacent} and Lemma \ref{uniclo} and is left to the reader.

  \subsection{Local averaged densities}
   We define, for $s \in \R$, the quantity
\begin{equation}
\label{localD}
\rbD_\rL (u, s)=\E(u, \Gamma_\rL (s-1, s+1)),  
\end{equation}
which represents a local \emph{averaged density} of the energy. Arguing as in the proof of  Corollary {\color{black}\ref{corollo}}, we obtain:
\begin{corollary}
 \label{corolli}
 Assume that $V$ satisfies  ${\color{black}(\mathbf{H}_{\mathbf{1}})}$, ${\color{black}(\mathbf{H}_{\mathbf{2}})}$ and  ${\color{black}(\mathbf{H}_{\mathbf{3}})}$ and assume that $\rL\geq 1$.   Let $a \in \R\cup\{-\infty\}$, $b \in \R \cup \{+\infty\}$  such that $b-a\geq 2$. Assume that 
 \begin{equation}
 \label{petitou1}
 \rbD_\rL(u, s) \leq \upeta_1 {\rm \ for \ } s \in [a, b]
 \end{equation}
 then, there exists  $\upsigma^*\in \Sigma$ such that 
 \begin{equation}
 \label{majorana2}
 \int_{\Gamma_\rL (a, b)}\vert u-\upsigma^*\vert^2\leq 16\rbC_2\E\left(u, \Gamma_\rL (a, b)\right).
 \end{equation}
  \end{corollary}
 \begin{proof} {\color{black}First assume that $a$ and $b$ are finite.} We perform a covering of $\Gamma_\rL(a, b)$ by squares  with sides of length $1$. For that purpose, we introduce the integers 
 $\ell =\llbracket \rL\rrbracket$, the largest integer smaller or equal to $\rL$, $m= \llbracket \ b-a\rrbracket$, the largest integer smaller or equal to $b-a$, so that, in view of our assumption $\ell \geq 1$ and $m\geq 2$.  Next, we introduce the points 
 \begin{equation*}
 \left\{
 \begin{aligned}
 &p_i=-\frac {\rL}{2}+i, {\rm \ for \ }  i=0, \ldots, \ell-1, {\rm \ and \ }  {\color{black}p_\ell}=\frac {\rL}{2}-1\\
 &z_j=a+j, {\rm \ for \ }  j=0, \ldots, m-1,\, {\rm \ and \ }  z_m=b-1, {\rm and \ consider \ the  \  cubes \ }\\
 &\rbQ_{i, j}=[p_i, p_i+1]\times [z_j, z_j+1], {\rm \ for \ } (i, j)\in  \rbA= \{0, \ldots,\ell\} \times \{0, \ldots, m\},
\end{aligned}
\right.
\end{equation*}
 so that $\rbQ_{i, j}$ is a square of sides of length $1$ and 
 $$
   \underset{(i, j)\in \rbA}\cup \rbQ_{i, j}=\Gamma_\rL (a, b).
   $$
Since
 {\color{black}
 $$ \rbQ_{i, j}\subset \Gamma_\rL (z_j, z_j+1), {\rm \ we \ have \ }
  \E(u, \rbQ_{i, j}) \leq \rbD_\rL(u,z_j)\leq \upeta_1, $$}
  by assumption \eqref{petitou1}. Hence, we deduce,   from  Corollary \ref{corollo}, that there exists a unique $\upsigma^\star\in \Sigma$, such that
\begin{equation}
\label{gluon}
  \int_{\rbQ_{i, j}}\vert u-\upsigma^*\vert^2\leq \rbC_2\E\left(u, \rbQ_{i, j}\right).
 \end{equation}
 We verify that almost every point in $\Gamma_\rL(a, b)$  belongs to \emph{at most  4 squares } in $\displaystyle{ \underset{(i, j)\in \rbA}\cup \rbQ_{i, j}}$. Applying 
inequality \eqref{gluon}, we deduce \eqref{majorana2}. {\color{black}For infinite endpoints, apply the result to an increasing exhaustion by finite strips. Uniqueness on overlapping squares fixes the same well throughout the exhaustion, and the result follows by monotone convergence.} 
 \end{proof} 

A variant is the following:
\begin{lemma}  
\label{locDL2}
 Assume that $V$ satisfies  ${\color{black}(\mathbf{H}_{\mathbf{1}})}$, ${\color{black}(\mathbf{H}_{\mathbf{2}})}$ and  ${\color{black}(\mathbf{H}_{\mathbf{3}})}$ and assume that $\rL\geq 1$. Let $u\in H^1_{\rm loc }(\Lambda_\rL)$  be such that  $\E(u, \Lambda_\rL)<+\infty$. Let $a \in \R\cup\{-\infty\}$, $b \in \R \cup \{+\infty\}$  such that $b-a\geq 2$. Assume that 
\begin{equation}
\label{minus}
\rbD_\rL (u, s)  \leq \upeta_1, {\rm \ for \ every \ } s \in [a, b]
\end{equation}
Then, there exists a unique $\upsigma^\star \in \Sigma$, such that 
\begin{equation}
\label{locoloco}
 \int_{\Gamma_\rL (s-1, s+1)}\vert u-\upsigma^{\color{black}\star}\vert^2\leq 16\rbC_2\rbD_\rL(u, s), {\rm \ for \ any  \ }  s\in [a+1, b-1].
\end{equation}
\end{lemma}

An immediate consequence of the above results is the following:

\begin{corollary} 
\label{infy}
Assume that $V$ satisfies  ${\color{black}(\mathbf{H}_{\mathbf{1}})}$, ${\color{black}(\mathbf{H}_{\mathbf{2}})}$ and  ${\color{black}(\mathbf{H}_{\mathbf{3}})}$ and assume that $\rL\geq 1$. Let $u\in\mathbb W_\rL$. If 
\begin{equation}
\label{micro}
\E(u, {\color{black}\Lambda_\rL})\leq \upeta_1 {\rm \ or \ } \rbD_\rL (u, s)  \leq \upeta_1, {\rm \ for \  any \ } s \in \R, 
\end{equation}
then there exists $\upsigma^\star \in \Sigma$, such that $u\in \mathbb W_\rL(\upsigma^\star, \upsigma^\star).$
\end{corollary}

As a consequence of Corollary \ref{infy}, we recover a well-known result (see e.g \cite{bonafoux}):
\begin{corollary} 
\label{balky}
Assume that $V$ satisfies  ${\color{black}(\mathbf{H}_{\mathbf{1}})}$, ${\color{black}(\mathbf{H}_{\mathbf{2}})}$ and  ${\color{black}(\mathbf{H}_{\mathbf{3}})}$ and assume that $\rL\geq 1$. Let 
$\upsigma^+$ and $\upsigma^-$ be two \emph{distinct} elements in $\Sigma$. Then
\begin{equation}
\label{macedoine}
\inf\{ \E(u), u \in  \mathbb W_\rL(\upsigma^-, \upsigma^+){\color{black}\}}= \rL \mathfrak m(\upsigma^-, \upsigma^+)\geq \upeta_1.
\end{equation}
\end{corollary}

  \subsection{ The $H^1$ norm versus the energy $\E$}
  An  immediate consequence  of Lemma \ref{locDL2} is the following: 
  
  \begin{corollary}
  \label{EimpliesH}
  Assume that $V$ satisfies  ${\color{black}(\mathbf{H}_{\mathbf{1}})}$, ${\color{black}(\mathbf{H}_{\mathbf{2}})}$ and  ${\color{black}(\mathbf{H}_{\mathbf{3}})}$ and assume that $\rL\geq 1$. Let $u\in H^1_{\rm loc }(\Lambda_\rL)$  be such that  $\E(u, \Lambda_\rL)<+\infty$. Let $a \in \R\cup\{-\infty\}$, $b \in \R \cup \{+\infty\}$ {\color{black}be} such that $b-a\geq 2$. Assume that 
 \eqref{minus} is satisfied, and let $\upsigma^\star\in \Sigma$ be as in Lemma \ref{locDL2}.  Then $v=u-\sigma^\star\in H^1(\Gamma_\rL(a, b))$   and, for some constant $\rbC_{\rm cap}>0$, we have 
\begin{equation}
\label{eimplesH0}
\Vert v \Vert_{H^1\left(\Gamma_\rL(a, b)\right)}^2=\Vert v \Vert_{L^2(\Gamma_\rL(a, b))}^2+ \Vert \nabla u\Vert^2_{L^2(\Gamma_\rL(a, b))}
 \leq \rbC_{\rm cap}\E(u, \Gamma(a, b)).
\end{equation}                                                                                                                                         
  \end{corollary}
\begin{proof}
{\color{black}For finite $a,b$, set $n_0=\lceil(b-a)/2\rceil$ and choose
$s_i=a+1+2i$ for $i=0,\ldots,n_0-2$, and $s_{n_0-1}=b-1$.
Then $s_i\in[a+1,b-1]$, the intervals $[s_i-1,s_i+1]$ cover $[a,b]$, and their multiplicity is at most two almost everywhere.
For an infinite endpoint use the analogous one-sided sequence, and for $[a,b]=\R$ use $s_i=2i$, $i\in\Z$.
Writing $I$ for the resulting index set, we add \eqref{locoloco} to obtain}
\begingroup\color{black}
\begin{equation}
\label{coverage}
\begin{aligned}
\Vert v\Vert_{L^2(\Gamma_\rL(a,b))}^2
&\leq\sum_{i\in I}\Vert v\Vert_{L^2(\Gamma_\rL(s_i-1,s_i+1))}^2\\
&\leq16\rbC_2\sum_{i\in I}\rbD_\rL(u,s_i)
\leq32\rbC_2\E(u,\Gamma_\rL(a,b)).
\end{aligned}
\end{equation}
\endgroup
The {\color{black}result} follows, choosing ${\color{black}\rbC_{\rm cap}=32\rbC_2+2}$,
{\color{black}since $\int_{\Gamma_\rL(a,b)}\vert\nabla u\vert^2\leq2\E(u,\Gamma_\rL(a,b))$.}
\end{proof}

 \begin{lemma}  
 \label{EvsH}
 Assume that $V$ satisfies  ${\color{black}(\mathbf{H}_{\mathbf{1}})}$, ${\color{black}(\mathbf{H}_{\mathbf{2}})}$ and  ${\color{black}(\mathbf{H}_{\mathbf{3}})}$ and assume that $\rL\geq 1$. Let $u\in H^1_{\rm loc }(\Lambda_\rL)$  be such that  $\E(u, \Lambda_\rL)<+\infty$. Let $a \in \R\cup\{-\infty\}$, $b \in \R \cup \{+\infty\}$ {\color{black}be} such that $b-a\geq 2$. Assume that {\color{black}\eqref{minus}} is satisfied, and let $\upsigma^\star\in \Sigma$ be as in Lemma \ref{locDL2}. Then, {\color{black}with $v=u-\upsigma^\star$,} we have
\begin{equation}
\label{lautre}
\E(u, \Gamma_\rL (a, b))\leq 4\rbK_0 \Vert v \Vert_{H^1(\Gamma_\rL(a, b))}^2.
\end{equation}
 \end{lemma}
 \begin{proof}  We have, in view of Lemma \ref{technical}
  $$
  0\leq V(\upsigma +y)\leq \rK_0 \left(\vert y \vert^2+ \vert y \vert^m\right), {\rm  \ for \  any \ } y \in \R^k.
 $$
Using this inequality with $y=u(x)-\upsigma^\star$ and integrating, we  are led, for any $s \in {\color{black}[a+1, b-1]}$, 
\begin{equation}
\label{back}
\int_{\Gamma_\rL (s-1, s+1)} V(u(x)) \rd x \leq  \rK_0 \left(\Vert v \Vert_{L^2(\Gamma_\rL(s-1, s+1))}^2+   \Vert v \Vert_{L^m(\Gamma_\rL(s-1, s+1))}^m\right)
\end{equation}
By the Sobolev embedding ${\color{black}H^1(\Gamma_\rL(s-1,s+1))\hookrightarrow L^m(\Gamma_\rL(s-1,s+1))}$, we obtain
{\color{black}
$$
\Vert v\Vert_{L^m(\Gamma_\rL(s-1,s+1))}^2
\leq C_S\Vert v\Vert_{H^1(\Gamma_\rL(s-1,s+1))}^2
\leq\rbC_4\rbD_\rL(u,s).
$$
Here $C_S$ is independent of $\rL\geq1$, as follows by a covering by unit squares, and $\rbC_4=C_S(16\rbC_2+2)$; the last inequality follows from Lemma \ref{locDL2} and the gradient part of the energy.} Going back to \eqref{back}, we are led to 
\begin{equation*}
\int_{\Gamma_\rL (s-1, s+1)} V(u(x)) \rd x \leq  \rK_0 \left(\Vert v \Vert_{L^2(\Gamma_\rL(s-1, s+1))}^2+  
\rbC_4^{\frac{m}{2}} \rbD_\rL(u, s)^\frac{m}{2}\right)
 \end{equation*}
 and thus, imposing $\rbK_0\geq 4$, and invoking assumption {\color{black}\eqref{minus}}, 
 \begin{equation*}
 \begin{aligned}
 \rbD_\rL(u, s)=\E(u, \Gamma_\rL (s-1, s+1)) &\leq  (\rK_0 +1)\Vert v \Vert_{H^1(\Gamma_\rL(s-1, s+1))}^2+  
 \rbC_4^{\frac{m}{2}}\rK_0 \rbD_\rL(u, s)^{\frac m2} \\
 &\leq (\rK_0 +1)\Vert v \Vert_{H^1(\Gamma_\rL(s-1, s+1))}^2+
  \rbC_4^{\frac{m}{2}}\rK_0  \upeta_1^{\frac{m}{2}-1}{\color{black}\rbD_\rL(u, s)}. 
 \end{aligned}
 \end{equation*}
 Since $m\geq 3$, we have $\displaystyle{\frac{m}{2}-1\geq \frac12}$, so that we may impose an additional  smallness condition  of $\upeta_1$, namely  
 $\displaystyle{{\color{black}\rbK_0}\rbC_4^{\frac{m}{2}} \upeta_1^{\frac{m}{2}-1}\leq \frac 14}$, so that, since we assume $\rbK_0\geq 4$
 $$\rbD_\rL(u, s) \leq 2\rK_0 \Vert v \Vert_{H^1(\Gamma_\rL(s-1, s+1))}^2.$$
Using a covering as in \eqref{coverage}, we finally obtain the desired result \eqref{lautre}. 
\end{proof}

  \subsection{Proof of Lemma \ref{asympto2D0}}
 We first assume $\rL\geq1$. We only deal with the behavior at $+\infty$, the argument at $-\infty$ being similar. Since, by assumption, $\E(u, \Lambda_\rL)<+\infty$, we deduce that $\E(u, \Gamma_\rL(a, +\infty){\color{black})}\to 0$ as $a\to + \infty$. It follows that there exists $a_0$ such that
 \begin{equation*}
  \E\left(u, \Gamma_\rL (a_0-1, +\infty)\right)\leq \upeta_1.
 \end{equation*}
We are hence in position to apply Corollary \ref{corolli}  to assert that there exists an element $\upsigma^+\in \Sigma$ such that 
\begin{equation}
\int_{\Gamma_\rL (a_0, +\infty)} \vert u-\upsigma^+\vert^2< +\infty.
\end{equation}
 Hence $\displaystyle{\int_{\Gamma_\rL (a, +\infty)} \vert u-\upsigma^+\vert^2 \to 0}$ as $a\to + \infty$. For $0<\rL<1$, choose an integer $q$ with $q\rL\geq1$, apply the preceding argument to the periodic extension on $\Lambda_{q\rL}$, and restrict back to $\Lambda_\rL$. The proof is complete.

  \subsection{On the structure of $\mbWL(\upsigma^-, \upsigma^+)$}
    
  \subsubsection{Proof of Lemma \ref{affinage}}
 We divide  the proof of Lemma \ref{affinage}  in Steps.
 
      \smallskip
  \noindent
 {\it Step 1: $\mbWL(\upsigma^-, \upsigma^+) \subset \{\bv_0\}+ H^1_{\rm per}(\Lambda_\rL)$}. Let $u\in\mbWL(\upsigma^-, \upsigma^+)$. It follows from Lemma \ref{asympto2D0} that 
 \begin{equation*}
  \int_{\Gamma_\rL (0, +\infty)} \vert u-\upsigma^+\vert^2< +\infty, 
  \end{equation*}
whereas, by construction we have 
\begin{equation*}
  \int_{\Gamma_\rL (0, +\infty)} \vert \bv_0-\upsigma^+\vert^2< +\infty, 
  \end{equation*}
  Setting  $v=u-\bv_0$, we obtain, combining the two estimates, 
  \begin{equation*}
  \int_{\Gamma_\rL (0, +\infty)} \vert v\vert^2\leq 2 \left( \int_{\Gamma_\rL (0, +\infty)} \vert u-\upsigma^+\vert^2+ \int_{\Gamma_\rL (0, +\infty)} \vert \bv_0-\upsigma^+\vert^2\right)< +\infty, 
    \end{equation*}
The same estimates hold on $\Gamma_\rL(-\infty,0)$ with $\upsigma^-$, 
so that $v \in L^2(\Lambda_\rL)$. Since $\E({\color{black}u}, \Lambda_\rL)<+\infty$ and, since by construction, $\E(\bv_0, \Lambda_\rL)<+\infty$, we have also 
$\nabla v \in L^2(\Lambda_\rL)$, hence $v\in H^1(\Lambda_\rL)$.  Concerning the periodic boundary condition, we have 
\begin{equation*}
{\color{black}u}(-\frac{\rL}{2}, \cdot)={\color{black}u}(\frac{\rL}{2}, \cdot) {\rm \ and \ } \bv_0(-\frac{\rL}{2}, \cdot)=\bv_0(\frac{\rL}{2}, \cdot) {\rm \ in \ the \ sense \ of \ traces,  \ }
\end{equation*}
 so that, by subtracting the two identities, we obtain $\displaystyle{v(-\frac{\rL}{2}, \cdot)=v(\frac{\rL}{2}, \cdot)}$. Hence, $v\in H^1_{\rm per}(\Lambda_\rL)$.

  \bigskip
  \noindent
 {\it Step 2: $ \{\bv_0\}+H^1_{\rm per}(\Lambda_\rL)\subset \mbWL(\upsigma^-, \upsigma^+)$}.  Let $v \in H^1_{\rm per}(\Lambda_\rL)$, and $u=\bv_0+v$. We have to show  that  
 \begin{equation}
 \label{butul}
\mathbb V(u)\equiv  \int_{\Lambda_\rL}V(u(x))\rd x<+\infty
 \end{equation}
 First, we have
  \begin{equation}
 \begin{aligned}
\int_{\Gamma_\rL (1, +\infty)} V(u(x))\rd x&=\int_{\Gamma_\rL (1, +\infty)} V(v(x)+ \bv_0)\rd x  \\
&=\int_{\Gamma_\rL (1, +\infty)} V(v(x)+ \upsigma^+)\rd x \\
&\leq
\rK_0\int_{\Gamma_\rL (1, +\infty)} \vert v(x) \vert^m+ \vert v(x) \vert^2\rd x
&<+\infty, 
\end{aligned}
\end{equation}
 where for the last inequality, we used the embedding $H^1(\Lambda_\rL) \hookrightarrow L^m  (\Lambda_\rL)$. Similarly, we prove 
 $\displaystyle{\int_{\Gamma_\rL (-\infty, -1)} V(u(x))\rd x< \infty}$  whereas
 $\displaystyle{\int_{\Gamma_\rL (-1, 1)} V(u(x))\rd x< \infty}$  follows directly from \eqref{due}.
 Hence the proof is complete.

 \medskip
 We notice the following consequence of Lemma \ref{technical}.
  \begin{corollary}
  \label{corona}
   Assume   ${\color{black}(\mathbf{H}_{\mathbf{1}})}$, ${\color{black}(\mathbf{H}_{\mathbf{2}})}$ and  ${\color{black}(\mathbf{H}_{\mathbf{3}})}$   hold.  Let $v \in \HL$ and let $u=v+ \bv_0\in \mbWL(\upsigma^-, \upsigma^+)$.  Set 
 $\tilde V(v)= V(v+\bv_0)${\color{black}.} 
  We have the upper bounds, for some constants $\rK_2>0$ and $\rK_3>0$
 \begin{equation}
 \label{technip}
 \left\{  
 \begin{aligned}
 &\tilde V(v(x)) \leq \rK_2\left( \vert v (x)\vert^2+ \vert v  (x) \vert^m+ \chi_0(x)  \right)\\
 &\vert  \nabla_u{\tilde V}(v)  \vert \leq \rK_3\left( \vert v \vert+ \vert v \vert^m+ \chi_1 (x)\right){\color{black},}
  \end{aligned}
 \right. 
 \end{equation}
 where $\chi_0$  and $\chi_1$ are smooth functions   with  \emph{compact support} in $\Lambda_\rL $.  In  particular, for any $v \in \HL$, we have, for some constants $\rK_4(p)>0$, 
 \begin{equation}
  \label{Lptude}
  \left\{
  \begin{aligned}
&\tilde V (v)\in L^p(\Lambda_\rL), {\rm \ for \  }  1\leq p <+\infty  {\rm \ and \ }\nabla_v \tilde V(v)\in L^p(\Lambda_\rL), {\rm \ for \  } 2\leq p < +\infty, {\rm \ with \ } \\
&\Vert  \tilde V(v) \Vert_{L^p(\Lambda_\rL)} \leq \rK_4(p)\left(  \Vert v \Vert_{H^1(\Lambda_\rL)}^2+  \Vert v \Vert_{H^1(\Lambda_\rL)}^m+1\right)
  {\rm \ for \  }  1\leq p <+\infty, {\rm \ and \ } \\
  &\Vert \nabla_u \tilde V(v) \Vert_{L^p(\Lambda_\rL)} \leq \rK_4(p)\left(  \Vert v \Vert_{H^1(\Lambda_\rL)}+  \Vert v \Vert_{H^1(\Lambda_\rL)}^m+1\right)
  {\rm \ for \  } 2\leq p < +\infty.
 \end{aligned}
 \right.
 \end{equation}
 \end{corollary}

  \subsection{Properties of the functional  $\E$}
  
 \begin{lemma}  Assume   ${\color{black}(\mathbf{H}_{\mathbf{1}})}$  ${\color{black}(\mathbf{H}_{\mathbf{2}})}$ and  ${\color{black}(\mathbf{H}_{\mathbf{3}})}$. 
 The  function $\tilde V$ is well defined, continuous, twice differentiable from $H^1_{\rm per}(\Lambda_\rL)$  to $L^p(\Lambda_\rL)$, for any $1\leq p <+\infty$ and 
   \begin{equation}
 \rd \tilde  V(v)(w)=\nabla_v V(u)\cdot w,  {\rm \ for \ } w \in H^1_{\rm per}(\Lambda_\rL), {\rm \ where \ } u=v+ \bv_0.
 \end{equation}
  \end{lemma}
 
 \begin{proof}
 It relies on the fact that, thanks to estimates \eqref{technip}, the function   $\tilde V$ is a  Nemitskii {\color{black}operator.} We skip the details.
\end{proof}
 
   Next, we consider the functional $\mathbb V$ defined on $H^1_{\rm per}(\Lambda_\rL)$, by 
 $$\mathbb V (v)= \int_{\Lambda_\rL} V(v+\bv_0)=\int_{\Lambda_\rL} V(u), {\rm \ for \ any  \ } v \in \HL, {\rm \ with \ } u=v+ \bv_0. $$
  
 \begin{proposition}
 \label{deriveV}
   The functional $\mathbb V$ is $C^1$ on $\HL$ and 
   $$ \rd \mathbb V (v){\color{black}\cdot} w=\int_{\Lambda_\rL} \nabla_v V(u){\color{black}\cdot} w,  {\rm \ for \ any \ } w \in \HL, 
 {\rm \ where \ } u=v+\bv_0.
 $$
  \end{proposition}
  
    \subsection{Proof of proposition \ref{differentiability}}
    We observe that the quadratic term 
    $\displaystyle{\rbJ(u)=\frac{1}{2}\int_{\Lambda_\rL}\vert \nabla u (x)\vert^2 \rd x}$ is also $C^1$, with
    $$
    \rd J(u)w=\int_{\Lambda_\rL} \nabla u\cdot \nabla w, {\rm \ for \ } w \in \HL.
    $$ 
    Combining this with the result of Proposition \ref{deriveV} gives the conclusion.


  \section{Some properties  of solutions}  
  \label{propsol}
 
 \subsection{Remarks on the periodicity condition}
 {\bf Periodisation and translations}.
 For a map $f:\Lambda_\rL \to \R^k$, we define its periodisation $f_{\rm per}$ according to \eqref{periodisation}, that is 
 $$f_{\rm per}(x+n \rL\rbe_1)=f(x), {\rm \  for \  any \ } x \in \Lambda_\rL, {\rm and \ any \ } n \in \Z.$$
 If $g$ is a map from $\R^2\to \R^k$, and $a \in \R^2$,  we denote by $\tau_a(g)$ its translation by $a$; that is, $\tau_a(g)$ is defined by
 $$
 \tau_a(g)(x)=g(x-a), {\rm \ for \  any \ }  x \in \R^2,  
 $$
 and finally, for $f:\Lambda_\rL \to \R^k$ and $a \in \R^2$, we denote by $\tilde{\tau}_a(f): \Lambda_\rL\to \R^k$ the restriction of $\tau_a(f_{\rm per})$ to $\Lambda_\rL$.  We have
 the {\color{black}following} elementary observation, whose proof is left to the reader:
 
 \begin{lemma}
 \label{elem}
 Let $f\in L^1(\Lambda_\rL)$ and $a \in \R^2$. We have the identities
 \begin{equation}
 \label{elem0}
 \int_{\Lambda_\rL} f(x) \rd x=\int_{\Lambda_\rL+\{a\} }f_{\rm per}(x) \rd x =\int_{\Lambda_\rL} {\tilde {\tau}_a} f(x) \rd x.
 \end{equation}
 Moreover $\Vert \tilde{\tau}_a(f)\Vert_{L^p(\Lambda_\rL)}=\Vert f\Vert_{L^p(\Lambda_\rL)}$.  If $f \in \HL$, then $f_{\rm per} \in H^1_{\rm loc} (\R^2)$, 
 $\nabla \tilde{\tau}_a(f)=\tilde{\tau}_a(\nabla f)$  and 
 $\Vert \tilde{\tau}_a(f)\Vert_{H^1(\Lambda_\rL)}=\Vert f\Vert_{H^1(\Lambda_\rL)}$.
 \end{lemma}
 Notice that for  the fact that $f_{\rm per} \in H^1_{\rm loc} (\R^2)$, we use the trace theorem.
 
 
 \bigskip
 \noindent
 {\bf  Translating the differential $\rd \E$}.  We have:
 \begin{lemma}
 \label{derper}
   Let    $u \in \mbWL(\upsigma^-, \upsigma^+)=\{\bv_0\}+\HL$. We have, for any $a \in \R^2$
   \begin{equation}
   \label{devper0}
   \rd \E(\tilde{\tau}_a(u)) \left(w\right)=\rd \E(u) \left(\tilde{\tau}_{-a} (w) \right), {\rm \ for \ any \ } w \in \HL.
      \end{equation}
   \end{lemma}
   \begin{proof} We have, in view of proposition \ref{differentiability}, for any  $w \in \HL$
    \begin{equation}
    \label{weakiper}
    \begin{aligned}
  \rd \E(\tilde{\tau}_a(u)) (\tilde{\tau}_a(w))&=    
  \int_{\Lambda_\rL } \left( \nabla \tilde{\tau}_a(u)\cdot \nabla \tilde{\tau}_a(w)+  \nabla_u V (\tilde{\tau}_a(u))\cdot \tilde{\tau}_a(w) \right)\rd x\\
   &=\int_{\Lambda_\rL } \tilde \tau_a \left( \nabla u{\color{black}\cdot} \nabla w+ \nabla_u V(u)\cdot w\right)\rd x 
   =\int_{\Lambda_\rL }  \nabla u{\color{black}\cdot} \nabla w+ \nabla_u V(u)\cdot w \rd x\\
   & =\rd \E (u)(w).
  \end{aligned}
  \end{equation}
    We notice that $\tilde{\tau}_a \circ \tilde{\tau}_{-a}={\rm Id}$.  Hence, applying \eqref{weakiper} with the map $\tilde{\tau}_{-a}(w)$ we obtain the desired identity {\color{black}\eqref{devper0}}.
    \end{proof}

\begin{corollary}  
\label{translatos}
Suppose that $u\in \{\bv_0\}+ \HL$ is a critical point of $\E$.   We have, for any $a\in \R^2$, 
$$\Vert \rd E(u_a) \Vert_{{\HL}^*}=\Vert \rd E(u) \Vert_{{\HL}^*}, {\rm \ where \ } u_a=\tilde {\tau}_{a}(u).$$
In particular, if $u$ is a critical point of $\E$, then 
$u_a=\tilde {\tau}_{a}(u)$ is  also a  critical point of $\E$, i.e. $\rd\E(u_a)=0$.
\end{corollary}
\subsubsection{Proof of Proposition \ref{regularita}}
{\color{black}Testing the critical-point identity with the periodisation $\sum_{n\in\Z}\varphi(\cdot+n\rL\rbe_1)$ of any $\varphi\in C_c^\infty(\R^2,\R^k)$ gives $-\Delta\rbu+\nabla_uV(\rbu)=0$ weakly on $\R^2$. Standard local elliptic regularity then gives $\rbu\in C^2(\R^2)$.} Taking $a=(\frac \rL 2, 0)$, we see that $u_a$ is  $C^2$ near $x_1=0$. Going back to $u$ this yields the Neumann condition 
\begin{equation}
\label{neumann}
\frac{\partial u}{\partial x_1}(-\frac{\rL}{2}, x_2)= \frac{\partial u}{\partial x_1 }(\frac \rL 2, x_2), {\rm \ for \ any \ }  x_2 \in \R.
\end{equation}
Since 
$$\R^2= \underset {n \in {\color{black}\Z}}\cup \Lambda_\rL^n, {\rm \ where \ }  \Lambda_\rL^n=\left [\frac{2n-1}{2}\rL, \frac{2n+1}{2}\rL\right]\times \R, $$
we have for any $\varphi \in \mathcal D(\R^2)$, using integration by parts
\begin{equation*}
\begin{aligned} 
\int_{\R^2} {\color{black}\bigl(}\nabla \rbu\cdot \nabla \varphi+  \nabla_u V (\rbu){\color{black}\cdot} \varphi{\color{black}\bigr)\rd x} &
=\underset {n \in \Z}\sum \int_{\Lambda_\rL^n} {\color{black}\bigl(}\nabla \rbu\cdot \nabla \varphi+  \nabla_u V (\rbu){\color{black}\cdot} \varphi{\color{black}\bigr)\rd x} \\
&=\underset {n \in \Z}\sum \int_{\Lambda_\rL^n} \left(-\Delta \rbu+\nabla_u V (\rbu)  \right){\color{black}\cdot}\varphi{\color{black}\,\rd x}\\
 &-\underset {n \in \Z} \sum \int_{\{\frac{2n-1}{2}\rL\} \times \R} \varphi (\frac{2n-1}{2}\rL, x_2){\color{black}\cdot} \frac{\partial \rbu}{\partial x_1}(\frac{(2n-1)}{2}\rL, x_2){\color{black}\,\rd x_2}\\
 &+ {\color{black}\underset {n \in \Z}\sum} \int_{\{\frac{2n+1}{2}\rL\} \times \R} \varphi (\frac{2n+1}{2}\rL, x_2){\color{black}\cdot} \frac{\partial \rbu}{\partial x_1}(\frac{(2n+1)}{2}\rL, x_2){\color{black}\,\rd x_2}\\
&=0, 
\end{aligned}
\end{equation*}
the boundary terms canceling thanks to the periodicity condition.
Hence $\displaystyle{-\Delta \rbu+\nabla_u V (\rbu)}=0$ on $\R^2$ in the sense of distributions{\color{black}, which} yields \eqref{Rdeux}.

 \subsection{Basic uniform estimates for solutions to \eqref{themainthing}}
\subsubsection {Local regularity} 
 Recall that $u\in \{\bv_0\}+ \HL$ is a critical point of $\E$ if and only if 
    \begin{equation}
\label{weaki}
\int_{\Lambda_\rL }{\color{black}\bigl(}\nabla u\cdot \nabla w+  \nabla_u V (u){\color{black}\cdot} w{\color{black}\bigr)\rd x}=0, {\rm \ for \ any \ } w \in \HL.
\end{equation}
 If we take $w=\varphi\in \mathcal D \left(\mathring{ \Lambda}_\rL\right)$, we see that   
 $$
 \int_{\Lambda_\rL } {\color{black}\bigl(}\nabla u\cdot \nabla \varphi+  \nabla_u V (u){\color{black}\cdot} \varphi{\color{black}\bigr)\rd x}=\langle -\Delta u+\nabla_u V (u), \varphi\rangle _{{\color{black}\mathcal D', \mathcal D}}, 
 $$
  and hence, combining with \eqref{weaki}, we are led to
 \begin{equation}
 \label{distribos}
 \Delta u=\nabla_uV(u)=\underset{j=1}{\overset{k}\sum}\frac{\partial V}{\partial u_j}(u)\,  \rbf_j, 
 {\rm \ in \ }\mathcal D'\left(\mathring {\Lambda}_\rL\right),  {\rm \ where \ } \rbf_i=(0, \ldots, 1, \ldots, 0)\in \R^k.
 \end{equation}

{\color{black}Standard regularity theory gives the following result.}

\begin{lemma}     Assume   ${\color{black}(\mathbf{H}_{\mathbf{1}})}$  ${\color{black}(\mathbf{H}_{\mathbf{2}})}$. If  $u\in {\color{black}\mathbb W_\rL}$  is  a solution to \eqref{themainthing} satisfying  $\E(u, \Lambda_\rL)<+\infty$, then,
$u \in C^1_{\rm loc} (\Lambda_\rL)$.
\end{lemma}
\begin{proof}  {\color{black}If} $u \in \mathbb W_\rL$, then the right hand side of \eqref{distribos} is in $L_{\rm loc}^p(\R^2)$, for any $p\geq 2$. {\color{black}The} result follows from standard elliptic theory{\color{black}.}
\end{proof}

\subsubsection{Using the invariance by translation}
\label{invariance}
We provide in this part, uniform estimates  on  a solution $u$,  its gradient as well as its second derivative, which  depend only on the nature of the potential $V$, but not on the  solution itself.  In the proofs,  we invoke the fact that it suffices to obtain the  respective estimates  only for the   origin ${\rm \bf O}=(0,0)$.   Indeed, given any $a \in \Lambda_\rL$, we have, by definition 
\begin{equation}
\label{lactose}
\tau_{-a}(u)(\rbO)=u(a), {\rm \ for \ any \ } a \in \Lambda_\rL.
\end{equation}
  By  Corollary \ref{translatos}
   the translated map $u_a\equiv \tilde {\tau}_{-a}(u)$ is also a solution to \eqref{distribos}, with $\E(u_a)=\E(u)$.  Hence, obtaining estimates at $\rbO$, {\color{black}for} arbitrary solutions, yields corresponding estimates for   $u(a)$.


\subsubsection{An upper bound on solutions  to \eqref{themainthing} depending only on $V$}

The following uniform estimate is given in \cite{bethu} for domains. We specify it for periodic solutions:
\begin{lemma} 
\label{maximodo}
 Assume   ${\color{black}(\mathbf{H}_{\mathbf{1}})}$  ${\color{black}(\mathbf{H}_{\mathbf{2}})}$ and  ${\color{black}(\mathbf{H}_{\mathbf{3}})}$. There exists a constant ${\rbM}_0$ depending only on the potential $V$, such that, if  $u\in \mathbb W_\rL$  is  a solution to \eqref{themainthing} satisfying  $\E(u, \Lambda_\rL)<+\infty$, then, we have
\begin{equation}
\label{maximodo0}
\vert u(x) \vert \leq \rbM_0, {\rm \ for \ } x \in \Lambda_\rL.
\end{equation}
\end{lemma}

\begin{proof}  
 We argue as in \cite{BBH0, bethu}. We compute, using equation \eqref{elip0}
       \begin{equation}
       \label{multiplie}
       \begin{aligned}
       \frac 12\Delta \vert u \vert^2&=u\cdot \Delta u+ \vert \nabla u \vert^2   
       = u \cdot \nabla_u  V( u)+ \vert \nabla u \vert^2  \\
     &  \geq u \cdot \nabla_u  V( u),   {\rm \ on \ } \Lambda_\rL. 
       \end{aligned}
       \end{equation}
     On the other hand, it follows from assumption ${\color{black}(\mathbf{H}_{\mathbf{2}})}$,   that there exists some constant $C_\infty \geq 0$ such that 
     \begin{equation}
     \label{convexita}
     y{\color{black}\cdot}\nabla  V( y )\geq \upalpha_\infty \vert y \vert^2- C_\infty,  {\rm \ for \ any \ } y \in \R^k,
     \end{equation}
       Hence, combining \eqref{multiplie} and \eqref{convexita} we obtain the  inequality
     \begin{equation*}
    -\frac 12\Delta \vert u \vert^2  +\upalpha_\infty  \left( \vert u\vert^2 -\frac{ C_\infty}{\upalpha_\infty}\right) \leq 0 {\rm \ on \ } \Lambda_\rL.
     \end{equation*}
     We introduce next the function $\displaystyle{W=\vert u \vert^2-\left(\frac{C_\infty}{\upalpha_\infty}+2{\rbR}_0\right)}$, where $\rbR_0$ is defined in \eqref{grandr0}. We are led to  the differential  inequality
       \begin{equation}
       \label{prince}
    -\frac 12\Delta W   +\upalpha_\infty  W \leq 0  {\rm \ on \ } \Lambda_\rL. 
     \end{equation}
   {\color{black}Invoking Lemma \ref{asympto2D0} and the finite gradient energy, we may choose sequences $a_n\to-\infty$ and $b_n\to+\infty$ such that
   \[
   \Vert u(\cdot,a_n)-\upsigma^-\Vert_{H^1(-\rL/2,\rL/2)}
   +\Vert u(\cdot,b_n)-\upsigma^+\Vert_{H^1(-\rL/2,\rL/2)}\longrightarrow0.
   \]
   Indeed, the squared slice norms are integrable on the corresponding half-lines. The one-dimensional Sobolev embedding gives uniform convergence on these slices. Since $\nabla V(\upsigma^\pm)=0$, \eqref{convexita} gives $|\upsigma^\pm|^2\leq C_\infty/\upalpha_\infty$. Thus, as $\rbR_0>0$, for all sufficiently large $n$,}
   \begingroup\color{black}
   \begin{equation}
   \label{aumilieu}
   W(x_1,a_n)<0,\qquad W(x_1,b_n)<0,
   \qquad x_1\in[-\rL/2,\rL/2].
   \end{equation}
   \endgroup
    Next assume  by contradiction that there exists some point $\rba\in \Lambda_\rL$ such that 
    \begin{equation}
    \label{patapouf}
    W(\rba) > 0.
    \end{equation}
    {\color{black}Choose $n$ large enough so that $\rba\in\Gamma_\rL(a_n,b_n)$ and \eqref{aumilieu} holds. By compactness and periodicity in $x_1$, $W$ attains a positive maximum at an interior point of this cylinder, which we again denote by $\rba$.} Invoking  the invariance principle of subsection \ref{invariance}, we may assume without loss of generality that $\rba=\rbO$.   Since $\rbO$ is a maximum point, we deduce that$ -\Delta W(\rbO)\geq 0$ and hence, 
     by assumption \eqref{patapouf} , we are led to $\displaystyle{-\frac12\Delta W(\rbO)+\upalpha_\infty  W(\rbO) > 0}$, a contradiction with \eqref{prince}. The proof is complete.
   \end{proof}
\subsubsection{Uniform gradient estimates}

    \begin{proposition}
    \label{classic}
     Assume   ${\color{black}(\mathbf{H}_{\mathbf{1}})}$, ${\color{black}(\mathbf{H}_{\mathbf{2}})}$ and  ${\color{black}(\mathbf{H}_{\mathbf{3}})}$ and $\rL\geq 1$.  Let $u\in \mathbb W_\rL$  be a solution of \eqref{themainthing}. Then $u$ is differentiable  on $\Lambda_\rL$  and there exists a constant 
    $\rbC_{\rm grd}>0$, depending only on $V$, 
     such that
     \begin{equation}
   \label{classicic0}
 \vert    \nabla u  \vert (x) \leq \rbC_{\rm grd}, {\rm \ for \  any \ } x \in \Lambda_\rL.
  \end{equation}
    \end{proposition}
   
 \begin{proof} In view of the invariance principle of subsection \ref{invariance}, it suffices to bound ${\color{black}\nabla}u(\rbO)$.  We invoke Lemma A.1 of \cite{BBH0}, which  asserts that, if  $v$ is a solution on  some domain $\mathcal U$ of 
 $\R^n$ of  $-\Delta v=f$,  then we have the inequality
 \begin{equation}
 \label{bbhoo}
 \vert \nabla v \vert^2(x) \leq \rbC_{\rm 0}\left (\Vert f \Vert_{L^{\infty}(\mathcal U)} \Vert v \Vert_{L^{\infty}(\mathcal U)}  +
 \frac{1}{{\rm dist } (x, \partial \mathcal U)^2} \Vert v \Vert_{L^{\infty}(\mathcal U)}^2\right),   {\rm \ for \ all \ } x \in \mathcal  U.
 \end{equation}
 We apply inequality \eqref{bbhoo} to the solution ${\color{black}u_{\rm per}}$, with  source term  $f=-\nabla_u V ({\color{black}u_{\rm per}}) $  on the domain 
 $\mathcal U=\B(1)$, for $x=\rbO$: This  yields, invoking Lemma \ref{maximodo}
  \begin{equation}
  \vert    \nabla u  \vert^2 (\rbO) \leq \rbC_{\rm 0} \left ( \Vert \nabla_uV(u)\Vert_\infty  \Vert u \Vert_{L^{\infty}}+ \Vert u \Vert_{L^{\infty}}^2  \right)
  \leq \rbC_{\rm 0} \left (\rbN_0 \rbM_0+\rbM_0^2 \right)\equiv \rbC_{\rm grd}^2, 
  \end{equation}
  where $\displaystyle{\rbN_0=\sup\{ \vert \nabla V(\rby)\vert, \vert \rby\vert \leq \rbM_0.\}}$
  \end{proof}

   \subsubsection{ An upper bound for the second derivative}
    We might iterate the previous estimate on solutions at  order $2$. 

   \begin{lemma} 
    \label{regu}   
   Assume that the potential $V$ satisfies assumptions ${\color{black}(\mathbf{H}_{\mathbf{1}})}$, ${\color{black}(\mathbf{H}_{\mathbf{2}})}$ and  ${\color{black}(\mathbf{H}_{\mathbf{3}})}$, and let $\rL\geq 1$. Then, if $ u\in \{\bv_0\}+ {\color{black}\HL}$ is a solution to \eqref{weaki}, then  $u\in C^2(\Lambda_\rL)$ and $\rbu=u_{\rm per}\in C^2(\R^2)$.  We have the estimate, for some constant $\rbC_{\rm nd}>0$
\begin{equation}
\label{smoothy}
\vert u(x)\vert+\vert  \nabla u (x)\vert+ \vert \nabla^2 u (x)\vert  \leq \rbC_{\rm nd}, \forall  x \in \Lambda_\rL 
\end{equation}
\end{lemma}
   
   \begin{proof}      Differentiating the equation with respect to $x$, we obtain
    \begin{equation}
    \label{secondder}
    -\Delta \left( \frac{\partial u}{\partial x_i}\right) ={\color{black}-}\underset{1\leq j, \ell\leq k} \sum\frac{\partial^2 V}{\partial u_j\partial u_\ell}(u) 
\frac{\partial u_\ell}{\partial x_i}\rbf_j, 
       {\rm \ on \ }\mathring {\Lambda}_\rL,  {\rm \ in \ the \ sense \ of \ distributions}.
   \end{equation}
We conclude as in Proposition \ref{classic},   using again estimate \eqref{bbhoo} for equation  \eqref{secondder} and  invoking the uniform bounds  \eqref{maximodo0} and \eqref{classicic0}.
 \end{proof}  

\subsubsection{Remarks on homoclinic solutions}
\begin{definition}
We say that a solution $u\in \WL$ is \emph{homoclinic} if $u \in \mathbb W_\rL(\upsigma^\star, \upsigma^\star)$ for some
$\upsigma^\star \in \Sigma$. 
\end{definition}
The main result in this part is:

\begin{proposition} 
\label{mono}
Assume   ${\color{black}(\mathbf{H}_{\mathbf{1}})}$, ${\color{black}(\mathbf{H}_{\mathbf{2}})}$, and ${\color{black}(\mathbf{H}_{\mathbf{3}})}$. {\color{black}Let $\rL\geq1$ and let $u\in\mathbb W_\rL(\upsigma^\star,\upsigma^\star)$, for some $\upsigma^\star\in\Sigma$, be a non constant solution of \eqref{themainthing}.} Then, we have 
\begin{equation}
\E(u, \Lambda_\rL) \geq  \upgamma_{\rm hc} \log \rL+ \rbC_{\rm hc}, 
\end{equation}
where $\upgamma_{\rm hc}>0$ and $\rbC_{\rm hc}$  are  constants depending only on $V$.
\end{proposition}
  
The proof relies on several ingredients presented as separate results.

\begin{lemma} 
\label{givesx0}
Assume that  ${\color{black}(\mathbf{H}_{\mathbf{1}})}$, ${\color{black}(\mathbf{H}_{\mathbf{2}})}$, and ${\color{black}(\mathbf{H}_{\mathbf{3}})}$  hold{\color{black}, and let $\rL\geq1$}. Let $u \in 
 \mathbb W_\rL(\upsigma^\star, \upsigma^\star)$ be a \emph{non constant} solution to \eqref{themainthing}. Then, there exists a 
 {\color{black}point} $\rbx_0\in \Lambda_\rL$ such that  
 \begin{equation}
 \label{argumentation}
 \vert u(\rbx_0)-\upsigma^\star\vert  = \frac{\upmu_0}{2}.
 \end{equation}
 Set $\displaystyle{r_0=\inf \{\frac14,  \frac{\upmu_0}{4\rbC_{\rm grd}}\}}$. Then we have 
 \begin{equation}
 \label{integreV}
 \int_{\B\left(\rbx_0, r_0\right)}V(u_{\rm per} (\rbx) )\rd \rbx  \geq \upgamma_{\rm pot} , {\rm \ where \ }
  \upgamma_{\rm pot} =\frac{\pi{\color{black}\lambda_0}}{64} \upmu_0^2 r_0^2>0.
 \end{equation}
\end{lemma}
\begin{proof} We first claim that there exists a point $\rbx_1\in \Lambda_\rL$ such that 
\begin{equation}
 \label{argumentation0}
 \vert u(\rbx_1)-\upsigma^\star \vert \geq  \upmu_0.
 \end{equation}
{\it Proof of the claim \eqref{argumentation0}}. We argue by contradiction. 
If \eqref{argumentation0} is not true, then  we have
\begin{equation}
\label{argumentation1}
\vert u(\rbx)-\upsigma^\star\vert  \leq \upmu_0, {\rm \ for \ any \  } \rbx \in \Lambda_\rL
\end{equation}
We take $w=u-\upsigma^\star$ as test function in \eqref{weaki} so that we are led to
\begin{equation}
\label{formeintegrale}
\int_{\Lambda_\rL} {\color{black}\bigl(}\vert \nabla u \vert^2 + \nabla_u V(u(\rbx))\cdot (u(\rbx)-\upsigma^\star){\color{black}\bigr)}\rd \rbx =0.
\end{equation}
In view of Proposition \ref{potto} and \eqref{argumentation1}, we have 
$$
\frac12 {\color{black}\lambda_0} \vert u(\rbx)-\upsigma^\star \vert^2 \leq \nabla_u V(u(\rbx))\cdot (u(\rbx)-\upsigma^\star), {\rm \ for \  any \ }
 x \in \Lambda_\rL.
$$
Integrating and going back to \eqref{formeintegrale}, we derive
$$
\int_{\Lambda_\rL} {\color{black}\bigl(}\vert \nabla u \vert^2 +\frac12 {\color{black}\lambda_0} \vert u(\rbx)-\upsigma^\star \vert^2{\color{black}\bigr)} \rd \rbx \leq 0, 
$$
hence $u$ is constant, $u=\upsigma^\star$, a contradiction, so that the claim is established.

\smallskip
We invoke next Lemma \ref{asympto2D0} to assert that there exists a point $\rbx_2 \in \Lambda_\rL$ such that 
$$ \vert u(\rbx_2)-\upsigma^\star \vert \leq \frac{\upmu_0}{4}.$$
We consider the segment $[\rbx_1, \rbx_2]\subset \Lambda_\rL${\color{black}.} By the continuity of the function $\rbx \mapsto \vert u(\rbx)-\upsigma^\star \vert$  we deduce by the intermediate value theorem, the existence of an element $\rbx_0$ satisfying \eqref{argumentation}.

\smallskip
Finally, we turn  to \eqref{integreV}.   In view of Proposition  \ref{classic}, we have 
\[
{\color{black}\vert u_{\rm per} (\rbx_0)-u_{\rm per} (\rbx) \vert \leq \rbC_{\rm grd} \vert \rbx_0-\rbx\vert,}
\]
so that 
$$\frac{\upmu_0}{4} \leq \vert {\color{black}u_{\rm per}}(\rbx)-\upsigma^\star\vert  \leq \upmu_0, {\rm \ if  \ }  
\vert \rbx-\rbx_0\vert \leq r_1=\frac{ \upmu_0}{4\rbC_{\rm grd}}.$$
Invoking Proposition  \ref{potto} once more, we deduce
 $$ \frac{{\color{black}\lambda_0}}{64} \upmu_0^2\leq V({\color{black}u_{\rm per}(\rbx)}),{\rm \  if \  \rbx }\in \B(x_0, r_0). $$
Setting $r_0=\inf\{\frac 14, r_1\}$, and integrating on $\B(\rbx_0, r_0)$ we obtain the desired result.
\end{proof}

 The classical Pohozaev identity reads   in dimension $2$ (see for instance \cite{bethu}, Lemma 3.3 for the proof) {\color{black}as follows:}

\begin{lemma}  
    \label{poho} Assume   ${\color{black}(\mathbf{H}_{\mathbf{1}})}$ {\color{black}and} ${\color{black}(\mathbf{H}_{\mathbf{2}})}$. Let $u$  be  a solution of \eqref{elip0} on $\R^2$. We have, for any  radius $0<r$      
    \begin{equation}
   \label{poho1}
   \left\{
   \begin{aligned}
 \frac{1}{r} \int_{\B^2(x,  r)} V(u) &=\frac 14\int_{\partial \B^2(x,  r) } 
  \left(
  \left\vert \frac{\partial  u}{\partial \tau}\right \vert^2 -\left\vert \frac{\partial  u}{\partial r} \right\vert^2 
  +2V(u)
  \right) {\rm d} \tau, \\
  &=\frac{1}{2} \int_{\partial \B^2(x,  r) } e(u)(\ell) \rd \ell-
  \frac{1}{2} \int_{\partial \B^2(x,  r) }
  \left\vert \frac{\partial  u}{\partial r} \right\vert^2 (\ell) \rd \ell
  \end{aligned}
  \right.
    \end{equation}
    where $\tau$ denotes arclength on $\partial \B(x, r)$, so that $\displaystyle{ {\rm d} \tau=r\rd \theta}$ corresponds to the  $\mathcal H^1$-measure on $\S^1(r)$ and 
    $\displaystyle{\frac{\partial  }{\partial \tau}=\frac{1}{r} \frac{\partial  }{\partial \theta}}$.
    \end{lemma}

\begin{corollary}  
\label{pohocanta} 
Assume   ${\color{black}(\mathbf{H}_{\mathbf{1}})}$ and  ${\color{black}(\mathbf{H}_{\mathbf{2}})}$. Let $u$  be  a solution of \eqref{elip0} on $\R^2$. We have, for any  radius $0<r$   
\begin{equation}
\label{superU}
\frac{\rd }{\rd r} \E\left(\B(\rbx, r)\right) \geq \frac{2}{r} \int_{\B^2(\rbx,  r)} V(u), {\rm \ for \ any \ } \rbx \in \R^2.
\end{equation}
\end{corollary}
\begin{proof}
Invoking  Fubini's theorem,  and  \eqref{poho1}, we have, for any $x \in \R^2$, 
\begin{equation}
\begin{aligned}
\frac{\rd }{\rd r} \E\left(\B(\rbx, r)\right)&=\int_{\partial \B(\rbx, r)} e(u)(\ell) \rd \ell=\frac{{\color{black}2}}{r} \int_{\B^2(\rbx,  r)} V(u)+ 
 {\color{black}\int}_{\partial \B^2(x,  r) }
  \left\vert \frac{\partial  u}{\partial r} \right\vert^2(\ell) \rd \ell \\
  &\geq \frac{2}{r} \int_{\B^2(\rbx,  r)}V(u).
\end{aligned}
\end{equation}
\end{proof}

\begin{proof} [Proof of  Proposition \ref{mono}]  Let $\rbx_0$ be given by Lemma \ref{givesx0}. In view of the invariance principle of subsection \ref{invariance}, we may assume that $\rbx_0=\rbO${\color{black}. We use the periodic extension of $u$, still denoted by $u$.} In view of  \eqref{superU}, 
we have
\begin{equation}
\label{routedu}
\frac{\rd }{\rd r} \E\left(\B({\color{black}\rbO}, r)\right) \geq \frac{2}{r} \int_{\B^2({\color{black}\rbO},  r)} V(u) \geq \frac{2}{r}\int_{\B^2(\rbx_0,  r_0)}V(u) \geq {\color{black}\frac{2}{r}}\upgamma_{\rm pot}, {\rm \  for  \ }   r_0\leq r \leq  \frac{\rL}{2}.
\end{equation}
Integrating \eqref{routedu} from $r_0$ to 
$\displaystyle{\frac{\rL}{2}}$, we obtain
$$\E\left (u,\B(\rbO, \frac \rL2)\right)  \geq \E\left (u,\B(\rbO, r_0)\right)+2\upgamma_{\rm pot} \log\left(  \frac{\rL}{2r_0} \right) $$

{\color{black}Since $\E(u,\B(\rbO,r_0))\geq0$ and $\B(\rbO,\rL/2)\subset\Lambda_\rL$, the result follows with $\upgamma_{\rm hc}=2\upgamma_{\rm pot}$ and $\rbC_{\rm hc}=-2\upgamma_{\rm pot}\log(2r_0)$.}
\end{proof}

\subsection{Behavior on expanding strips: Uniform bounds}
When $u$ is a solution to \eqref{themainthing}, {\color{black}the result of} Lemma \ref{asympto2D0} can be improved as follows:

\begin{lemma} 
\label{asympto2D0bis}
 Assume that the potential $V$ satisfies assumptions ${\color{black}(\mathbf{H}_{\mathbf{1}})}$, ${\color{black}(\mathbf{H}_{\mathbf{2}})}$, ${\color{black}(\mathbf{H}_{\mathbf{3}})}$, and let $\rL>1$.
  Let $u \in \mathbb W_\rL$ be a solution to \eqref{themainthing}.  Then there exist elements $\upsigma^-$ and $\upsigma^+$ in $\Sigma$ such that 
\begin{equation}
\label{asympto2D00}
 \sup\left \{ \vert u(x_1, x_2)-\upsigma^\pm \vert, \, x_1 \in [-\frac \rL 2, \frac \rL 2]\right \} \to 0 {\rm \ as \ } x_2\to\pm \infty. 
\end{equation}
\end{lemma}

The proof relies on the following: 
\begin{lemma}
\label{corned}
Assume that the potential $V$ satisfies assumptions ${\color{black}(\mathbf{H}_{\mathbf{1}})}$, ${\color{black}(\mathbf{H}_{\mathbf{2}})}$, ${\color{black}(\mathbf{H}_{\mathbf{3}})}$, and let $\rL>1$.  Let $a \in \R\cup\{-\infty\}$, $b \in \R \cup \{+\infty\}$  such that $a \leq b$.  {\color{black}There exist constants} $\upeta_2>0$ and $\rbK_5>0$ {\color{black}such that, if} $u \in \mathbb W_\rL$ {\color{black}is} a solution to \eqref{themainthing} satisfying
\begin{equation}
\label{borneo}
\rbD_\rL(u, s) \leq \upeta_2, {\rm \ for \ } s \in [a, b],
\end{equation}
then there exists  a unique $\upsigma^\star\in \Sigma$ such that 
\begin{equation}
\label{bono}
\vert u(x)-\upsigma^\star\vert^2 \leq \rbK_5 \underset{s \in [a, b]}\sup \rbD_\rL(u, s),{\rm \ for \ any \ } x \in \Gamma_\rL(a, b).
\end{equation}
\end{lemma}
\begingroup\color{black}
\begin{proof}
We use the periodic extension of $u$ to $\R^2$, still denoted by $u$.
It satisfies $\Delta u=\nabla V(u)$ in the weak sense: this follows directly by testing \eqref{weaki} with the periodisation
$\sum_{n\in\Z}\varphi(x+n\rL(1,0))$ of any $\varphi\in C_c^\infty(\R^2,\R^k)$.
All constants below depend only on $V$.
Fix $p=(p_1,s)\in\Gamma_\rL(a,b)$ and let $Q=p+(-1/2,1/2)^2$ and $Q'=p+(-1/4,1/4)^2$.
Since $\rL>1$, periodicity gives $E_Q:=\E(u,Q)\leq\rbD_\rL(u,s)\leq\upeta_2$.
By Lemma \ref{uniclo}, there is $\upsigma_Q\in\Sigma$ such that, with $w=u-\upsigma_Q$,
\begin{equation}
\label{corned-local-hone}
\Vert w\Vert_{H^1(Q)}^2\leq (\rbC_2+2)E_Q.
\end{equation}
Assume $\upeta_2\leq1$. By \eqref{technical0} and the two-dimensional Sobolev embeddings, we obtain
\begin{equation}
\label{corned-local-source}
\begin{aligned}
\Vert\nabla V(\upsigma_Q+w)\Vert_{L^4(Q)}
&\leq C\bigl(\Vert w\Vert_{L^4(Q)}+\Vert w\Vert_{L^{4m}(Q)}^m\bigr)\\
&\leq C\bigl(E_Q^{1/2}+E_Q^{m/2}\bigr)\leq C E_Q^{1/2}.
\end{aligned}
\end{equation}
The interior $W^{2,4}$ estimate for $\Delta w=\nabla V(\upsigma_Q+w)$ therefore yields
\begin{equation}
\label{corned-local-elliptic}
\begin{aligned}
\Vert w\Vert_{C^1(Q')}&\leq C\Vert w\Vert_{W^{2,4}(Q')}\\
&\leq C\bigl(\Vert w\Vert_{L^4(Q)}+\Vert\nabla V(\upsigma_Q+w)\Vert_{L^4(Q)}\bigr)
\leq C E_Q^{1/2}.
\end{aligned}
\end{equation}
In particular, for a constant $\rbK_5>0$ depending only on $V$,
\begin{equation}
\label{corned-local-pointwise}
\vert u(p)-\upsigma_Q\vert^2\leq\rbK_5 E_Q\leq\rbK_5\rbD_\rL(u,s).
\end{equation}
Choose $\upeta_2\leq\min\{\upeta_1,1\}$ small enough that
$\rbK_5\upeta_2<\upgamma_0/16$, where $\upgamma_0$ is defined in \eqref{infdist}.
Thus $u(\Gamma_\rL(a,b))$ lies in the union of the pairwise disjoint balls of radius $\sqrt{\upgamma_0}/4$ centred at the wells.
Since $u$ is continuous and $\Gamma_\rL(a,b)$ is connected, the well $\upsigma_Q$ is the same $\upsigma^\star$ for every $p$.
The claimed estimate follows from \eqref{corned-local-pointwise}; uniqueness follows from the same separation of the wells.
\end{proof}
\endgroup

\begin{proof} [Proof of Lemma \ref{asympto2D0bis} completed] {\color{black}There exists $a_0>0$ such that} 
$$ \E(u, \Gamma_\rL(a_0-1, +\infty){\color{black})} \leq \upeta_2, {\rm \ and \ }  \E(u, \Gamma_\rL(-\infty, -a_0+1){\color{black})} \leq \upeta_2.$$
We deduce from Lemma \ref{corned} and Lemma \ref{asympto2D0} that, {\color{black}for $x\in\Gamma_\rL(a+1,+\infty)$ in the first line and $x\in\Gamma_\rL(-\infty,a-1)$ in the second,}
\begin{equation}
\label{lala}
\left\{
\begin{aligned}
\vert u(x)-\upsigma^+   \vert^2 &\leq \rbK_5 \E(u, \Gamma_\rL(a, +\infty)) \to 0, {\rm \ as \ } a\to + \infty, {\rm \ and \ } \\
\vert u(x)-\upsigma^-   \vert^2 &\leq \rbK_5 \E(u, \Gamma_\rL(-\infty, a )) \to 0, {\rm \ as \ } a\to - \infty.
\end{aligned}
\right.
\end{equation}
Hence the result.
\end{proof}



\subsection{Behavior on expanding strips: Energy estimates}

\begin{lemma}
\label{caccioplo}
Assume that the potential $V$ satisfies assumptions ${\color{black}(\mathbf{H}_{\mathbf{1}})}$, ${\color{black}(\mathbf{H}_{\mathbf{2}})}$, ${\color{black}(\mathbf{H}_{\mathbf{3}})}$, and let $\rL>1$.  Let $a \in \R\cup\{-\infty\}$, $b \in \R \cup \{+\infty\}$  such that $b-a\geq 2$, and $R>0$.  There exist constants $\upeta_3>0$ and $\rbK_6>0$, such that, for any $u \in \mathbb W_\rL$ that is a solution to \eqref{themainthing}, if \begin{equation}
\label{caccioplo-small}
\rbD_\rL(u, s) \leq \upeta_3, {\rm \ for \ any \ } s \in (a-R, b+R),
\end{equation}
then we have  
\begin{equation}
\label{framadate}
\E\left(u, \Gamma_\rL(a, b)\right)  \leq \frac{\rbK_6}{R^2} \E(u, \Gamma_\rL (a-R, b+R){\color{black})}.
\end{equation}
\end{lemma}
\begin{proof} As a preliminary remark, we observe  that, if we impose $\upeta_3 \leq \upeta_2$, as well as 
$\rbK_5\upeta_3 \leq \upmu_0^2$
then we may apply {\color{black}\eqref{bono}} to assert that there exists a unique $\upsigma^\star \in \Sigma$ such that 
\begin{equation}
\label{bonosse}
\vert u(x)-\upsigma^\star  \vert^{\color{black}2} \leq \rbK_5 {\color{black}\upeta_3}\leq \upmu_0^{\color{black}2} {\rm \ for \ any \ } x \in \Gamma_\rL(a-R, {\color{black}b}+R).
\end{equation}
The core of the proof is  then  similar to the classical  proof of Caccioppoli's inequality. Consider a  smooth plateau function  $\chi_R$ defined on $\R$ with  the following properties
\begin{equation*}
\left\{
\begin{aligned}
&\chi_R(x)=1 {\rm \  for \  } x \in [a, b],  \chi_R(x)=0, {\rm \ for \ } x \in \R\setminus [a-R, b+R], \\
& 0\leq \chi_R(x) \leq 1, {\rm \ for \ any \ } x\in \R, \\
& \vert  \chi_R'(x) \vert \leq \frac{2}{R}, {\rm \ for \ any \ } x\in \R.
\end{aligned}
\right.
\end{equation*}
Consider next the test function $w_R$ defined by 
\begin{equation}
\label{testus}
 w_R(x_1, x_2)=\left( u(x_1, x_2)-\upsigma^\star\right) \chi_R^2(x_2), {\rm \ for \ } (x_1, x_2) \in \Lambda_\rL{\rm \ so \ that \ }  w_R\in \HL.
 \end{equation}
  Using $w_R$ as test function in \eqref{weaki}, we obtain
 \begin{equation}
\label{weaki2}
\int_{\Lambda_\rL }{\color{black}\bigl(}\nabla u\cdot \nabla w_R+\nabla_u V(u){\color{black}\cdot}w_R{\color{black}\bigr)\,\rd x}=0, 
\end{equation}
We expand the first term on the left-hand side  as
\begin{equation}
\label{ristreto}
\begin{aligned}
 \nabla u{\color{black}\cdot} \nabla w_R&=\nabla u{\color{black}\cdot} \nabla ((u-\upsigma^\star)\chi_R^2) \\
 &={\color{black}\chi_R\nabla u\cdot\nabla\bigl((u-\upsigma^\star)\chi_R\bigr)
 +\chi_R'(x_2)\chi_R(u-\upsigma^\star)\cdot u_{x_2}} \\
 &={\color{black}\bigl\vert\nabla\bigl((u-\upsigma^\star)\chi_R\bigr)\bigr\vert^2
 -\chi_R'(x_2)(u-\upsigma^\star)\cdot\partial_{x_2}\bigl((u-\upsigma^\star)\chi_R\bigr)} \\
 &\quad {\color{black}+\chi_R'(x_2)(u-\upsigma^\star)\cdot\partial_{x_2}\bigl((u-\upsigma^\star)\chi_R\bigr)
 -\vert\chi_R'(x_2)\vert^2\vert u-\upsigma^\star\vert^2} \\
 &=\vert \nabla ((u-\upsigma^\star) \chi_R)\vert^2-\vert \chi_R'(x_2)\vert^2 {\color{black}\vert u-\upsigma^\star\vert^2}.
\end{aligned}
\end{equation}
In view of \eqref{bonosse} we may invoke  Proposition \ref{potto}, to assert that 
\begin{equation}
\label{quidbono}
\frac12\lambda_0\left  \vert u-\upsigma^\star\right\vert ^2 \leq  \nabla_u V (u){\color{black}\cdot} (u-\upsigma^\star)  {\rm \ on \ } \Gamma_\rL(a-R, b+R)
\end{equation} 
Combining  \eqref{weaki2},\eqref{ristreto} and \eqref{quidbono}, we are led to
\begin{equation}
\label{noloso}
\begin{aligned}
\int_{\Lambda_\rL}(\vert \nabla ({\color{black}(u-\upsigma^\star)\chi_R})\vert^2+\frac 12 \lambda_0\vert u-\upsigma^\star \vert^2 \chi_R(x_2)^2) \rd x
&\leq  \int_{\Lambda_\rL}\vert {\color{black}\chi_R'(x_2)}\vert^2 \vert u-\upsigma^\star\vert^2{\color{black}\,\rd x} \\
&\leq \frac{4}{R^2} \int_{\Gamma_\rL(a-R, b+R)}\vert u-\upsigma^\star \vert^2  \rd x.  \\
\end{aligned}
\end{equation}
Using Proposition \ref{potto} once more, we have 
\begin{equation}
\label{zermolo0}
\begin{aligned}
 \int_{\Gamma_\rL(a-R, b+R)}\vert u-\upsigma^\star \vert^2  \rd x  &\leq \frac{16}{\lambda_0} \int_{\Gamma_\rL(a-R, b+R)}V(u)  \rd x \\
&\leq \frac{16}{\lambda_0 } \E(u, \Gamma_\rL(a-R, b+R)), 
\end{aligned}
\end{equation}
whereas 
\begin{equation}
\label{zermolo1}
\begin{aligned}
\int_{\Lambda_\rL}(\vert \nabla ({\color{black}(u-\upsigma^\star)\chi_R})\vert^2+\frac 12 \lambda_0\vert u-\upsigma^\star \vert^2 {\color{black}\chi_R(x_2)^2}) \rd x
&\geq \int_{\Gamma_\rL(a, b)}\left(\vert \nabla u\vert^2+\frac 12 \lambda_0\vert u-\upsigma^\star \vert^2 \right)\rd x \\
&\geq \int_{\Gamma_\rL(a, b)}\left(\vert \nabla u\vert^2+\frac {\lambda_0}{2\lambda_{\rm max}^+} V(u) \right)\rd x \\
& \geq \frac {\lambda_0}{2\lambda_{\rm max}^+}\E(u, \Gamma_\rL(a, b)).
\end{aligned}
\end{equation}
Combining \eqref{zermolo0} and \eqref{zermolo1}, we obtain \eqref{framadate} with 
$\displaystyle{\rbK_6=\frac{{\color{black}128}\lambda_{\rm max}^+}{\lambda_0^2}}$. 
\end{proof}

\subsection{Decay estimates for solutions}
As a consequence of Lemma \ref{caccioplo}, we obtain the following decay property.

\begin{corollary}
\label{decay}
Assume that the potential $V$ satisfies assumptions ${\color{black}(\mathbf{H}_{\mathbf{1}})}$, ${\color{black}(\mathbf{H}_{\mathbf{2}})}$, ${\color{black}(\mathbf{H}_{\mathbf{3}})}$, and let $\rL>1$.   Let $u \in \mathbb W_\rL$ be  a solution to \eqref{themainthing}. There {\color{black}exist} some $a_0^+  \in \R$, $a_0^-\in \R$  and ${\color{black}\upsigma^+,\upsigma^-}\in \Sigma $ such that, for any ${\color{black}R>0}$, we have 
\begin{equation}
\label{decadence}
\left\{
\begin{aligned}
&\vert u(x)-{\color{black}\upsigma^\pm} \vert^2 \leq \frac{\rbK_7}{R^2}, {\rm \ for \ }  x \in \Gamma_\rL (a_0^++R, +\infty) ( {\color{black}\text{resp.}\ } x \in \Gamma_\rL(-\infty , a_0^--R) ) \\
&\E(u, \Gamma_\rL (a_0^++R, +\infty){\color{black})}\leq \frac{\rbK_6}{R^2}\upeta_3, ( {\color{black}\text{resp.}\ } 
\E(u, \Gamma_\rL ({\color{black}-\infty, a_0^--R}){\color{black})}\leq \frac{\rbK_6}{R^2}\upeta_3).
\end{aligned}
\right.
\end{equation}
\end{corollary}

\begin{proof}  Since $\E(u, \Lambda_\rL) <+\infty$, there exists some $\tilde {a_0}^+\in \R$ such that  $\E(u, \Gamma_\rL (\tilde {a}_0^+-1, +\infty){\color{black})}<\upeta_3$. It follows in particular that we have  
$$ \rbD_\rL(u, s) <\upeta_3, {\rm \ for \ any \ } s\in (\tilde {a_0}^+, +\infty[.$$
Given any $R>0$, we are in a position to apply Lemma \ref{caccioplo} with $a=\tilde {a}_0^++R$, $b=+\infty$, to assert that
\begin{equation}
\label{oubli}
\E(u, \Gamma_\rL (\tilde{a}_0^++R, +\infty){\color{black})}\leq \frac{\rbK_6}{R^2}\upeta_3{\color{black},}
\end{equation}
hence
$$\rbD_\rL(u, s)<\frac{\rbK_6}{R^2}\upeta_3, {\rm \ for \ } s>\tilde {a}_0^++R+1 {\color{black}.}$$
We set $a_0^+ =\tilde {a}_0^+ +1$ and $\rbK_7={\rbK_6\rbK_5}\upeta_3$. The first inequality in \eqref{decadence} for  $a_0^+$ then follows from {\color{black}\eqref{bono}}, whereas the second {\color{black}follows} directly from \eqref{oubli}. 
\end{proof}

\begin{remark}
\label{decadence01}
{\rm As a consequence of  Corollary \ref{decay},  setting $ \tilde R_0=\sup \{ \vert a_0^+\vert , \vert a_0^-\vert\}$, we have 
for $R\geq 2 \tilde R_0${\color{black}, $R>0$}, the inequality $\displaystyle{R-\vert a_0^\pm\vert \geq \frac{R}{2}}$, so that 
\begin{equation}
\label{decadence02}
\left\{
\begin{aligned}
&\vert {\color{black}u}(x)-{\color{black}\upsigma^\pm} \vert^2 \leq \frac{4\rbK_7}{R^2}, {\rm \ for \ }  
{\color{black}x\in\Gamma_\rL(R,+\infty)\ ({\rm resp.}\ x\in\Gamma_\rL(-\infty,-R))},  \\
&\E\left({\color{black}u},\Lambda_\rL \setminus  \Gamma_\rL (-R,  R)\right)\leq  \frac{8\rbK_6}{R^2}\upeta_3, 
\end{aligned}
\right.
\end{equation}

}
\end{remark}

\subsection{Proof of Proposition \ref{lwb0}: Lower bounds on the energy of  solutions}
The proof of Proposition \ref{lwb0} is a direct consequence of Lemma \ref{caccioplo}. Indeed, assume that 
$\E(\rbu, \Lambda_\rL) \leq\upeta_3$. Then, assumption {\color{black}\eqref{caccioplo-small}} is immediately satisfied for any $s \in \R$, and we may use \eqref{framadate} for any ${\color{black}a,b\in\R\text{ with }b-a\geq2}$ and $R>0$.   This yields
\begin{equation}
\label{framadate0}
\E\left(\rbu, \Gamma_\rL(a, b)\right)  \leq \frac{\rbK_6}{R^2} \E({\color{black}\rbu}, \Gamma_\rL (a-R, b+R){\color{black})}
\leq \frac{\rbK_6}{R^2} \E(\rbu,\Lambda_\rL).
\end{equation}
Since $R>0$ is arbitrary, we let $R\to +\infty$. This yields
$$
\E\left(\rbu, \Gamma_\rL(a, b)\right) =0, {\rm \ for \ any \ } a<b, {\rm \and \ hence \ } \E\left(\rbu, \Lambda_\rL\right)=0.
$$
The proof is hence complete, choosing $\upeta_{\rm min}=\upeta_3$.
\qed

\section{Some properties  of approximate solutions}
  We set $H=H^1_{\rm per}(\Lambda_\rL)$ and denote $H^\star$ its topological dual.  We recall that, given an open set $\Omega$, we \emph{may identify} the topological  dual of $H_0^1(\Omega)$  with a subset of $\mathcal D'(\Omega)$, namely the vector space of distributions $f \in \mathcal D'(\Omega)$  such that 
$$\Vert f \Vert_{H^{-1}(\Omega)}\equiv \sup \left \{  \vert  \langle f, \varphi \rangle_{\mathcal D', \mathcal D}{\color{black}\vert}, \, 
\Vert \varphi \Vert_{H^1(\Omega)}\leq  1, \varphi \in C_{\rm c}^\infty (\Omega) \right\}
<+ \infty.
$$
Any element in $\left(H^1_0(\Omega)\right)^\star$ yields  a distribution in $H^{-1}(\Omega)$, and conversely,  using the fact that  $C_c^\infty(\Omega)$ is dense in $H^1_0(\Omega)$ we may identify a distribution $f\in H^{-1}(\Omega)$ with an element in $\left(H^1_0(\Omega)\right)^\star$. Moreover, we have 
\begin{equation}
\label {identification}
\Vert f \Vert_{\left(H^1_0(\Omega)\right)^\star}=\Vert f \Vert_{H^{-1}(\Omega)}.
\end{equation}

\subsection{Regularity for approximate solutions}
Approximate solutions  have  some properties close to the properties of solutions, modulo perturbation terms. {\color{black}In this subsection, $u$ denotes its periodic extension to $\R^2$ whenever it is evaluated outside $\Lambda_\rL$.}
  \begin{lemma}
  \label{sansdec0}
  Assume that the potential $V$ satisfies assumptions ${\color{black}(\mathbf{H}_{\mathbf{1}})}$, ${\color{black}(\mathbf{H}_{\mathbf{2}})}$ and  ${\color{black}(\mathbf{H}_{\mathbf{3}})}$, and let $\rL>1$.
  Let $u\in \WL$. {\color{black}There exist maps $u_{\rm reg}\in W^{2,r}_{\rm loc}(\B(\rbO,1))$ and $u_{\rm irr}\in H^1_0(\B(\rbO,1))$ such that} 
  \begin{equation}
  \label{sansdec}
  u=u_{\rm reg}+\uir
  \end{equation}
  with $\uir \in H_0^1(\B(\rbO, 1))$, and, {\color{black}for} some constant $C(r)>0$ depending only on $r>2$ and $V$ and some universal constant $C>0$, we have 
  \begin{equation}
  \label{estdec}
  \left\{
  \begin{aligned}
  \Vert \ureg  \Vert_{W^{2, r}(\B(\rbO, \frac34){\color{black})} }&\leq C(r) \left(\E(u)+ \E(u)^{2m}+1 \right) {\rm\ and \ }\\
  \Vert \uir  \Vert_{H^1(\B(\rbO, 1))}&\leq C \Vert \rd E(u)\Vert_{H^\star}.
 \end{aligned}
  \right.
  \end{equation}
  \end{lemma}
  
 For the proof, we rely on the following:
\begin{lemma}  
\label{higherint}
 Assume that the potential $V$ satisfies assumptions ${\color{black}(\mathbf{H}_{\mathbf{1}})}$, ${\color{black}(\mathbf{H}_{\mathbf{2}})}$ and  ${\color{black}(\mathbf{H}_{\mathbf{3}})}$, and let $\rL>1$. 
Let $u \in \{\bv_0\}+\HL$. Then we have,
\begin{equation}
-\Delta u=-\nabla_u V(u)+ f, {\rm on \ }  \B(\rbO, 1)
\end{equation}
where $f  \in H^{-1}(\B(\rbO, 1))$ with 
\begingroup\color{black}
\begin{equation}
\label{Hmoins1}
\Vert f\Vert_{H^{-1}(\B(\rbO,1))}
\leq\sqrt2\,\Vert\rd\E(u)\Vert_{H^\star}.
\end{equation}
\endgroup
\end{lemma}

\begingroup\color{black}
\begin{proof}
For $\varphi\in C_c^\infty(\B(\rbO,1),\R^k)$, set
$\Phi(x)=\sum_{j\in\Z}\varphi(x+j\rL(1,0))$.
Periodicity and Proposition \ref{differentiability} give
$\langle-\Delta u+\nabla V(u),\varphi\rangle=\rd\E(u)(\Phi)$.
Since $\rL>1$, at most two terms overlap at any point, and
$\Vert\Phi\Vert_H\leq\sqrt2\,\Vert\varphi\Vert_{H^1(\B(\rbO,1))}$.
Thus $f=-\Delta u+\nabla V(u)$ satisfies \eqref{Hmoins1}.
\end{proof}
\endgroup

\begin{proof}[Proof of Lemma \ref{sansdec0}] Let $\uir \in {\color{black}H_0^1(\B(\rbO,1))}$ be the unique solution to the Dirichlet problem
\begin{equation}
\label{Laxmilgram}
-\Delta \uir =f,  {\rm \ in \ } \mathcal D'(\B(\rbO, 1)), \, {\rm \ and \ } \uir=0 {\rm \ on \ } \partial \B(\rbO, 1),
\end{equation}
obtained for instance, by the Lax-Milgram theorem. We have
\begin{equation}
\label{LM1}
\Vert \uir\Vert_{H_0^1(\B(\rbO, 1))} \leq  C \Vert f \Vert_{H^{-1}(\B(\rbO, 1){\color{black})}}\leq C\Vert\rd\E(u)\Vert_{H^\star},
\end{equation}
where we used \eqref{Hmoins1} for the last inequality.  Writing $u=\uir+\ureg$,  it follows that $\ureg$  satisfies the equation 
$$-\Delta\ureg=-\nabla_u V(u) {\rm \ in \ } \mathcal D'(\B(\rbO, 1)), {\rm \ and \ } \ureg=u {\rm \ on \ } \partial \B(\rbO, 1). $$ 
Hence,  by standard elliptic estimates, it follows that, for any $r>2$, we have

\begin{equation*}
\begin{aligned}
\Vert \ureg \Vert_{W^{2, r}(\B(\rbO, \frac34))}  
&\leq {\color{black}C(r)}\left(\Vert u \Vert_{{\color{black}H^1}(\B(\rbO, 1))}+\Vert \nabla_u V(u) \Vert_{L^r(\B(\rbO, 1))}\right) \\
&\leq C(r) \left(\E(u)+ \E(u)^{2m}+1 \right), \\
{\rm and \ } \Vert \uir  \Vert_{H^1(\B(\rbO, 1){\color{black})}}&\leq C \Vert \rd E(u)\Vert_{H^\star}, 
\end{aligned}
\end{equation*}
yielding hence \eqref{estdec}.
\end{proof}
\subsection{Caccioppoli inequalities for approximate solutions}

The Caccioppoli inequality for approximate solutions is technically more involved. It can be stated as follows:

\begin{proposition}
\label{caccioplon}
Assume that the potential $V$ satisfies assumptions ${\color{black}(\mathbf{H}_{\mathbf{1}})}$, ${\color{black}(\mathbf{H}_{\mathbf{2}})}$, ${\color{black}(\mathbf{H}_{\mathbf{3}})}$, and let $\rL>1$.  Let $a \in \R\cup\{-\infty\}$, $b \in \R \cup \{+\infty\}$  such that $b-a\geq 1$, and $N\in \N^\star$.  There exist constants $\upeta_4>0$ and $\rbK_6>0$, such that, for any $u \in \mathbb W_\rL$ {\color{black}with $\delta:=\Vert\rd\E(u)\Vert_{H^\star}\leq1$}, if
\begin{equation}
\label{bonoloc}
\rbD_\rL(u, s)\leq  \upeta_4, {\rm \ for \  any  \ } s \in (a-N^2, b+N^2), 
\end{equation}
then we have, for a constant $\rbC_4>0$
\begin{equation}
\label{ruegamma}
\E\left(u, \Gamma_\rL(a, b)\right)  \leq \frac{\rbK_6}{N} \E(u, \Gamma_\rL (a-N^2, b+N^2))+ \rbC_4 \Vert \rd E(u)\Vert_{H^\star}.
\end{equation}
{\color{black}The constants depend only on $V$.}
\end{proposition}
\begin{proof}  We  divide the proof into steps.
{\color{black}Taking $\rbK_6\geq1$, the case $N=1$ is immediate. Henceforth assume $N\geq2$.}

\smallskip
\noindent
{\it Step 1}. Assume that   $\upeta_4 \leq \upeta_1$. Then, if \eqref{bonoloc} is satisfied,  there  exists a unique $\upsigma^\star \in \Sigma$ such that 
\begin{equation}
 \label{majorana20}
 \int_{\Gamma_\rL (a-N^2, b+N^2)}\vert u-\upsigma^*\vert^2\leq 20\rbC_2\E\left(u, \Gamma_\rL (a-N^2, b+N^2)\right).
 \end{equation}
 {\color{black}Indeed, every unit square $Q$ in $\Gamma_\rL(a-N^2,b+N^2)$ has $\E(u,Q)\leq\upeta_4$. Corollary \ref{corollo} gives one common well $\upsigma^\star$ and
\[
\Vert u-\upsigma^\star\Vert_{H^1(Q)}^2\leq C\E(u,Q)\leq C\upeta_4.
\]
A covering by unit squares with overlap multiplicity at most four proves \eqref{majorana20}. Such coverings also exist for $1<\rL<2$ and for infinite endpoints. The same argument applies to every substrip of $\Gamma_\rL(a-N^2,b+N^2)$ of height at least one, with constants independent of its height and of $\rL$.}  We next consider the restriction of 
 $u$ to $\Gamma_\rL (a-N^2, b+N^2)$, and write
 \begin{equation}
 \label{Vr}
 u=\upsigma^\star + v, {\rm \ on \ } \Gamma_\rL (a-N^2, b+N^2).
 \end{equation}
It follows from \eqref{majorana20} that $v  \in H^1_{\rm per} (\Gamma_\rL (a-N^2, b+N^2))$, where 
\begin{equation*}
H^1_{\rm per} (\Gamma_\rL (a-N^2, b+N^2))=
\left\{ 
\begin{aligned}
&w \in H^1(\Gamma_\rL (a-N^2, b+N^2)), \\
&  w(-\frac{\rL}{2}, x_2)=w(+\frac{\rL}{2}, x_2), x_2\in (a-N^2, b+N^2)
\end{aligned}
\right\}, 
\end{equation*}
equipped with the standard $H^1$ norm. It follows from \eqref{majorana20} that
\begin{equation}
\Vert v \Vert_{H^1}^2\leq C \E\left(u, \Gamma_\rL (a-N^2, b+N^2)\right).
\end{equation}

\smallskip
\noindent
{\it Step 2: An averaging argument}. {\color{black}First suppose that $a,b$ are finite.} We claim that there exists $\tilde a \in [a-N^2, a] $ and $\tilde b \in [b, b+N^2]$ such that
\begin{equation}
\label{clameur}
\E\left(u,  A_N\right) \leq 
\frac{1}{N} \E(u, \Gamma_\rL(a-N^2, b+N^2){\color{black})},
\end{equation}
where $A_N=\Gamma_\rL (\tilde a-N, \tilde a)\cup \Gamma_\rL (\tilde b, \tilde b+N{\color{black})}$.
\begin{proof}[Proof of  the claim \eqref{clameur}]
\renewcommand{\qedsymbol}{}
\begin{equation*}
\left\{
\begin{aligned}
&\underset{i=1}{\overset {N} \sum}\left( \E\left(u, \Gamma ( a-Ni, a-N(i-1){\color{black})}\right)\right)=\E\left(u, \Gamma (a-N^2, a)\right), \\
&\underset{i=0}{\overset {N-1} \sum}\left( \E\left(u, \Gamma ( b+Ni, b+N(i+1){\color{black})}\right)\right)=\E\left(u, \Gamma (b, b+N^2)\right).
\end{aligned}
\right.
\end{equation*}
By an averaging argument,   there exists $n^- \in \{1, \ldots, N\}$, and $n^+ \in \{0, \ldots, N-1\}$ such that 
\begin{equation*}
\left\{
\begin{aligned}
&\E\left(u, \Gamma ( a-Nn^-, a-N({\color{black}n^-}-1){\color{black})}\right) \leq \frac{1}{N}\E\left (u, \Gamma_\rL (a-N^2, a)\right) \\
&\E\left(u, \Gamma ( b+N{\color{black}n^+}, b+N(n^++1){\color{black})}\right) \leq \frac{1}{N}\E\left (u, \Gamma_\rL (b, b+N^2{\color{black})}\right) \\
\end{aligned}
\right.
\end{equation*}
Setting $\tilde a=a-N(n^- -1)$ and  $\tilde b=b+ n^+ N$, the result follows.
\end{proof}
{\color{black}If $a=-\infty$, set $\tilde a=-\infty$ and omit the left component of $A_N$ and the corresponding averaging argument. If $b=+\infty$, proceed similarly on the right. Thus \eqref{clameur} remains valid, including when $A_N$ is empty.}

\smallskip
\noindent
{\it Step 3: The test function $w_N$.} We  follow  the proof of Lemma \ref{caccioplo}, except that now, the function $u\in {\color{black}\mathbb W_\rL}$ is not supposed to satisfy \eqref{themainthing}, and hence uniform estimates are not available.  
We use   again the function $w_N$ defined in \eqref{testus} with $\upsigma^\star$ given above, $R=N$ and $a$ and $b$ replaced by $\tilde a$ and $\tilde b$, that is $w_N=(u-\upsigma^\star) {\color{black}\chi_N^2}$, where
\begin{equation*}
\left\{
\begin{aligned}
&\chi_N(x)=1 {\rm \  for \  } x \in [\tilde a, \tilde b],  \chi_N(x)=0, {\rm \ for \ } x \in \R\setminus [{\color{black}\tilde a}-N, {\color{black}\tilde b}+N], \\
& 0\leq \chi_N(x) \leq 1, {\rm \ for \ any \ } x\in \R, \\
&  \vert \chi_N'(x) \vert \leq \frac{2}{N}, {\rm \ for \ any \ } x\in \R.
\end{aligned}
\right.
\end{equation*}
{\color{black}Set $v_N=(u-\upsigma^\star)\chi_N$. At an infinite endpoint no cutoff is imposed. The functions $v_N$ and $w_N=\chi_Nv_N$, extended by zero where necessary, belong to $H^1_{\rm per}(\Lambda_\rL)$.}
Using $w_N$ as a test function in \eqref{weakasse}, we obtain
\begin{equation}
\label{weakin}
\int_{\Lambda_\rL }{\color{black}\bigl(}\nabla u\cdot \nabla w_N+  \nabla_u V (u){\color{black}\cdot} w_N{\color{black}\bigr)\rd x}=\langle \rd \E(u), w_N\rangle_{H^\star, H}.
\end{equation}
 Using formula \eqref{ristreto}, \eqref{weakin} yields
\begingroup\color{black}
\begin{equation}
\label{nolosobis}
\begin{aligned}
\int_{\Lambda_\rL}\left(|\nabla v_N|^2+
\nabla_uV(u)\cdot(u-\upsigma^\star)\chi_N^2\right)\rd x
&\leq\int_{\Lambda_\rL}|\chi_N'|^2|u-\upsigma^\star|^2\rd x\\
&\quad+|\langle\rd\E(u),w_N\rangle_{H^\star,H}|.
\end{aligned}
\end{equation}
\endgroup
Notice also that 
\begingroup\color{black}
\begin{equation}
\label{VN}
\Vert w_N\Vert_H\leq C\Vert v_N\Vert_H
\leq C\bigl(\Vert v\Vert_{H^1(\Gamma_\rL(\tilde a,\tilde b))}+\Vert v\Vert_{H^1(A_N)}\bigr).
\end{equation}
\endgroup

\smallskip
\noindent
{\it Step 4: Revisiting the term  $\nabla_u V(u)\cdot u$.}   We claim that we have,  for any $y \in \R^k$
\begin{equation}
\label{quidbonobis}
\frac12\lambda_0\left  \vert y-\upsigma^\star\right\vert ^2 \leq  \nabla_u V (y){\color{black}\cdot} (y-\upsigma^\star) 
+\rbC_{\rm nl} \vert y-\upsigma^\star\vert^{m+1}, {\rm \ where \ }
\rbC_{\rm nl}=\rbK_0+\lambda_0 \upmu_0^{-m+1}.
\end{equation} 
\begin{proof}[Proof of the claim]
\renewcommand{\qedsymbol}{}
By Proposition \ref{potto}, we have 
\begin{equation}
\label{quidbono-approx}
\frac12\lambda_0\left  \vert y-\upsigma^\star\right\vert ^2 \leq  \nabla_u V (y){\color{black}\cdot} (y-\upsigma^\star), {\rm \ if   \ }  \vert  y-\upsigma^\star\vert \leq \upmu_0, 
\end{equation} 
so that \eqref{quidbonobis} is satisfied {\color{black}in} that case.  Next, consider the case $\vert  y-\upsigma^\star\vert \geq \upmu_0$. In that case,  we may apply inequality \eqref{luistrio} to assert that 
\begin{equation}
\label{luisquatuor}
\left\{
\begin{aligned}
\vert \nabla_u V (y){\color{black}\cdot} (y-\upsigma^\star)\vert  &\leq \rbK_0  \vert y-\upsigma^\star \vert^{m+1}, {\rm \ for \ }  \vert  y-\upsigma^\star\vert \geq \upmu_0,  \\
 \frac{1}2\lambda_0\left  \vert y-\upsigma^\star\right\vert ^2 &\leq \frac{\lambda_0}{\upmu_0^{m-1}} \vert y-\upsigma^\star\vert ^{m+1}, 
 {\rm \ for \ }  \vert  y-\upsigma^\star\vert \geq \upmu_0. 
\end{aligned}
\right.
\end{equation}
Combining the two inequalities, we obtain the result \eqref{quidbonobis}.
\end{proof}

\smallskip
\noindent
 {\it Step 5: Expanding  \eqref{nolosobis}}.
  We use the inequality \eqref{quidbonobis} with $y=u(x)$. Integrating on $\Lambda_\rL$, this yields
\begin{equation*}
\frac{\lambda_0}{2} \int_{\Lambda_\rL} \vert u-\sigma^\star\vert^2{\color{black}\chi_N^2}\leq \int_{\Lambda_\rL} \nabla_u V (u){\color{black}\cdot} w_N + \rbC_{\rm nl}\int_{\Lambda_\rL} \vert u-\upsigma^\star\vert^{m+1} {\color{black}\chi_N(x_2)^2}.
\end{equation*}
Going back to \eqref{nolosobis}, we obtain
\begingroup\color{black}
\begin{equation}
\label{nolosoter}
\begin{aligned}
\int_{\Lambda_\rL}\left(|\nabla v_N|^2+
\frac{\lambda_0}{2}|v|^2\chi_N^2\right)\rd x
&\leq\int_{\Lambda_\rL}\left(|\chi_N'|^2|v|^2+
\rbC_{\rm nl}|v|^{m+1}\chi_N^2\right)\rd x\\
&\quad+|\langle\rd\E(u),w_N\rangle_{H^\star,H}|.
\end{aligned}
\end{equation}
\endgroup
 Since $\chi_N$ {\color{black}vanishes} outside $\Gamma_\rL(\tilde a, \tilde b)\cup A_N$, where $A_N$ is {\color{black}introduced} in {\color{black}\eqref{clameur}} and {\color{black}$\chi_N^{\prime}$ is supported in $A_N$}, we are led to 
\begingroup\color{black}
\begin{equation}
\label{tropgentil}
\begin{aligned}
\int_{\Gamma_\rL(\tilde a,\tilde b)}\left(|\nabla v|^2+\frac{\lambda_0}{2}|v|^2
-\rbC_{\rm nl}|v|^{m+1}\right)\rd x
&\leq\int_{A_N}\left(\frac4{N^2}|v|^2+
\rbC_{\rm nl}|v|^{m+1}\right)\rd x\\
&\quad+C\Vert\rd\E(u)\Vert_{H^\star}\Vert v_N\Vert_H.
\end{aligned}
\end{equation}
\endgroup

\smallskip
\noindent
{\it Step 6: Analysis of inequality \eqref{tropgentil}}.
Invoking the Sobolev embeddings
{\color{black}on unit squares and a covering of each unit-height strip with bounded overlap, we obtain a constant $C$ independent of $\rL\geq1$ and $s$ such that}
\begin{equation}
\label{embed}
{\color{black}\Vert v\Vert_{L^{m+1}(\Gamma_\rL(s,s+1))}
\leq C\Vert v\Vert_{H^1(\Gamma_\rL(s,s+1))},
\qquad \Gamma_\rL(s,s+1)\subset\Gamma_\rL(a-N^2,b+N^2).}
\end{equation}
{\color{black}First decrease $\upeta_3$, if necessary, so that the absorption condition below holds; the earlier small-energy statements remain valid. Then choose $\upeta_4$ sufficiently small that the estimate in Step 1 gives $\Vert v\Vert_{H^1(\Gamma_\rL(s,s+1))}\leq\upeta_3$,}
so that 
\begin{equation}
\begin{aligned}
\int_{\Gamma_\rL(s, s+1)}\vert v \vert^{m+1} \rd x &\leq  C^{m+1}  \Vert v \Vert_{H^1 (\Gamma_\rL (s, s+1))}^{m+1}\\
& \leq C^{m+1} \upeta_3^{m-1}\Vert v \Vert_{H^1 (\Gamma_\rL (s, s+1))}^{2}{\color{black}.}
\end{aligned}
\end{equation}
A suitable covering argument then yields
\begin{equation*}
\left\{
\begin{aligned}
\int_{\Gamma_\rL(\tilde a,  \tilde b)}\vert v \vert^{m+1} \rd x &\leq 2C^{m+1} \upeta_3^{m-1}\Vert v \Vert_{H^1 (\Gamma_\rL (\tilde a, \tilde{\color{black}b}))}^{2} \\
\int_{A_N}\vert v \vert^{m+1} \rd x &\leq 2C^{m+1} \upeta_3^{m-1}\Vert v \Vert_{H^1 (A_N)}^{2}.
\end{aligned}
\right.
\end{equation*}
Going back to \eqref{tropgentil}, and setting $\rbc_0=\inf \{1, \lambda_0\slash 2\}$, we are led to
\begin{equation}
\label{H1AN}
\begin{aligned}
\Vert v \Vert_{H^1 (\Gamma_\rL (\tilde a, \tilde b ))}^{2}\left(\rbc_0-2C^{m+1} \rbC_{\rm nl}\upeta_3^{m-1}\right)
\leq &\Vert v \Vert_{H^1 ({\color{black}A_N})}^{2}\left(1+2\rbC_{\rm nl}C^{m+1} \upeta_3^{m-1}\right) \\
&+C \Vert \rd E(u)\Vert_{H^\star} \Vert v_N \Vert_{H^1}.
\end{aligned}
\end{equation}
We impose that  $\displaystyle{\upeta_3^{m-1} <\frac{\rbc_0}{{\color{black}4}C^{m+1} \rbC_{\rm nl}}}$, so that 
$\displaystyle{\rbc_0-2C^{m+1} \rbC_{\rm nl}\upeta_3^{m-1}\geq \frac{\rbc_0}{2}}$, 
 and  set 
$\rbc_1=1+2\rbC_{\rm nl}C^{m+1} \upeta_3^{m-1}$. Then  \eqref{H1AN} yields
\begin{equation}
\label{fatuite}
\begin{aligned}
\frac{{\color{black}\rbc_0}}{2}\Vert v \Vert_{H^1 (\Gamma_\rL (\tilde a, \tilde b ))}^{2}&\leq \rbc_1 \Vert v \Vert_{H^1 ({\color{black}A_N})}^{2} 
+C \Vert \rd E(u)\Vert_{H^\star} \Vert v_N \Vert_{H^1}.\\
\end{aligned}
\end{equation}
{\color{black}By \eqref{VN} and Young's inequality, and then \eqref{fatuite},
\[
\begin{aligned}
C\delta\Vert v_N\Vert_H
&\leq\frac{\rbc_0}{4}\Vert v\Vert_{H^1(\Gamma_\rL(\tilde a,\tilde b))}^2
+\Vert v\Vert_{H^1(A_N)}^2+C\delta^2,\\
\frac{\rbc_0}{4}\Vert v\Vert_{H^1(\Gamma_\rL(\tilde a,\tilde b))}^2
&\leq(\rbc_1+1)\Vert v\Vert_{H^1(A_N)}^2+C\delta^2.
\end{aligned}
\]
}

On the other hand, we 
have, by \eqref{eimplesH0} \begin{equation}
\label{eimplesH02}
\Vert v \Vert_{H^1(A_N{\color{black})}}^2 \leq \rbC_{\rm cap}\E(u, A_N) \leq \frac{\rbC_{\rm cap}}{N} \E(u, \Gamma_\rL(a-N^2, b+N^2){\color{black})},
\end{equation} 
{\color{black}where} we have used \eqref{clameur} for the last inequality.  Combining \eqref{eimplesH02} and \eqref{fatuite}, we obtain
\begingroup\color{black}
\begin{equation}
\label{vosotros}
\frac{\rbc_0}{4}\Vert v\Vert_{H^1(\Gamma_\rL(\tilde a,\tilde b))}^2
\leq\frac{(\rbc_1+1)\rbC_{\rm cap}}{N}
\E\left(u,\Gamma_\rL(a-N^2,b+N^2)\right)+C\delta^2.
\end{equation}
\endgroup
{\color{black}The same local covering estimate with exponent $m$ gives
\[
\int_{\Gamma_\rL(\tilde a,\tilde b)}|v|^m\leq C\upeta_4^{(m-2)/2}\Vert v\Vert_{H^1(\Gamma_\rL(\tilde a,\tilde b))}^2.
\]
Thus Lemma \ref{technical}, after decreasing $\upeta_4$ if necessary, gives \eqref{lautre} also when $1\leq\tilde b-\tilde a<2$, and we have}
\begin{equation}
\label{lautre23}
\E(u, \Gamma_\rL (\tilde a, \tilde b){\color{black})}\leq 4\rbK_0 \Vert v \Vert_{H^1(\Gamma_\rL(\tilde a, \tilde b))}^2,
\end{equation}
so that combining \eqref{vosotros} and \eqref{lautre}, we are led to
$$
\E(u, \Gamma_\rL (\tilde a, \tilde b){\color{black})}\leq {\color{black}\frac{16\rbK_0(\rbc_1+1)\rbC_{\rm cap}}{\rbc_0 N}} \E(u, \Gamma_\rL(a-N^2, b+N^2){\color{black})}+{\color{black}C\delta^2};
$$
Since $[a, b] \subset  [\tilde a, \tilde b]$, {\color{black}and $\delta^2\leq\delta$ by assumption,} the conclusion follows.
{\color{black}We enlarge $\rbK_6$ if necessary so that the same constant is valid here and in Lemma \ref{caccioplo}.}
\end{proof}

\section{The fate of Palais-Smale sequences}
\label{Passy}
In this  section we consider sequences $(u_n)_{n \in \N}$ in $\WL(\upsigma^-, \upsigma^+)$  whose energies are bounded and for which the differentials tend to $0$; that is, there exists some constant $M>0$ such that 
\begin{equation}
\label{psoide}
\left\{
\begin{aligned}
&E(u_n) \leq M, {\rm \ for \ any \ } n \in \N\\
&\rd E(u_n) \to 0 {\rm \ in \ }  {\HL}^\star, {\rm \ as \ } n \to +\infty.
\end{aligned}
\right.
\end{equation}

\subsection{ Concentration points and clustering}
\label{generalcluster}
 For $n \in \N$, we introduce the set $\mathcal L_n \subset \Z$ of concentration points of $u_n$, namely 
\begin{equation}
\label{mathcalAn}
\mathcal L_n=\left\{\ell \in \Z, \rbD(u_n, \ell)=\E(u_n, \Gamma_\rL (\ell-1, \ell+1){\color{black})}> \frac 14\upeta_4\right\} \subset \Z, 
\end{equation}
where $\upeta_4$ is the constant appearing in Proposition \ref{caccioplon}. Since the total energy is finite and   points in $\Lambda_\rL$ belong to at most two distinct sets $\Gamma_\rL(\ell-1, {\color{black}\ell+1})$ with $\ell \in \Z$, except for a set of zero measure, we deduce that 
\begin{equation}
\mathfrak q_n=\sharp (\mathcal L_n) \leq \frac{8M}{\upeta_4}.
\end{equation}
\begingroup\color{black}
We may discard finitely many terms so that $\Vert\rd\E(u_n)\Vert_{H^\star}\leq1$.
For the distinct endpoints in Proposition \ref{chainofsolutions}, $\mathcal L_n$ is nonempty: otherwise the two neighboring integer strips give $\rbD_\rL(u_n,s)\leq\upeta_4/2\leq\upeta_1$ for every real $s$, and Corollary \ref{infy} would imply equal endpoints.
\endgroup
Passing possibly to a subsequence (still denoted $(u_n)_{n \in \N}$), we may  assume that  $\mathfrak q=\mathfrak q_n$ is independent of $n$ and write
\begin{equation}
\label{mathcalCn}
\left\{
\begin{aligned}
&\mathcal L_n=\{\ell_{n, 1}, \ell_{n, 2}, \ldots, \ell_{n, \mathfrak q}\}, {\rm \ with \ } \ell_{n, 1}<\ell_{n, 2} <\ldots, \ell_{n, \mathfrak q} {\rm \ and \  set \ } \\
&\rbd_{n, j}=\ell_{n, j+1}-\ell_{n, j}, {\rm \ for \ } j=1, \dots, \mathfrak q-1, {\rm with } \ \rbd_{n, \mathfrak q}=+\infty.
\end{aligned}
\right.
\end{equation}
Passing possibly to a further subsequence, we may suppose that $\rbd_{n,j}$ has a limit in $\R^+\cup\{+\infty\}$, that is 
\begin{equation}
\rbd_{n, j} \to \rbd_j^\star\in \R^+\cup\{+\infty\},  {\rm \ as \ }   n \to +\infty, {\rm \ for \ } j=1, \ldots, \mathfrak q.
\end{equation}
Next, we  form  clusters of consecutive points which stay within finite distance from each other.  More precisely, we decompose the set of indices as follows, with $0<\mathfrak p_\star\leq \mathfrak q$, 
\begin{equation}
\label{clustering}
\left\{
\begin{aligned}
&\{1, \ldots, \mathfrak q\}=\underset{k=1}{\overset {\mathfrak p_\star} \cup} \mathcal J_k, {\rm \ where \ }  
\mathcal J_k=\{j_k, j_{k}+1, \ldots, j_{k}+\mathfrak m_k\}, \mathfrak m_k\in \N, {\rm \ with \ }\\
& j_{k}+\mathfrak m_k <j_{k+1}, \ \rbd_{j_k+\mathfrak m_k}^\star=\infty,  {\rm \ and \ } \\
&\rbd_{j_k+m}^\star<+\infty, {\rm \ for \ } m=0, \ldots,\mathfrak m_k-1.
\end{aligned}
\right.
\end{equation}
Finally, we set, for $n \in \N$
\begin{equation}
\rba_{n, k}=\ell_{n, j_k}, {\rm \ for \ } k=1, \ldots,  \mathfrak p_\star.
\end{equation}
Passing possibly to a further subsequence, we may assume that $\mathfrak p_\star$ does not depend on $n$.  Notice that it follows from \eqref{clustering} that there exists some $n_\star \in \N$, such that  we have 
  \begingroup\color{black}
\begin{equation}
\label{dismax}
\begin{aligned}
0\leq\ell_{n,j_k+m}-\ell_{n,j_k}&\leq R_k^\star-1,
\qquad R_k^\star=2+\sum_{r=0}^{\mathfrak m_k-1}\rbd_{j_k+r}^\star,\\
&m=0,\ldots,\mathfrak m_k,\quad n\geq n_\star.
\end{aligned}
\end{equation}
\endgroup

In particular, we have
\begin{equation}
\label{cestpetit}
\rbD(u_n, s) \leq \frac12 \upeta_4, {\rm \ for \ } 
s \in \R\setminus \underset {k={\color{black}1}}{\overset{\mathfrak p_\star}\cup} [\rba_{n, k}-R_k^\star, \rba_{n, k}+R_k^\star], {\rm \ for \ } n \geq {n_\star}.
\end{equation}

  \begin{lemma} 
\label{clusterisation}
We have   the following properties:
\begin{equation}
\label{clusterisation0}
\left\{
\begin{aligned}
& \rba_{n, k} \in \Z, {\rm \ for \ } k=1, \ldots, \mathfrak p_\star, {\rm \ and \ } \rba_{n, k}<\rba_{n, k+1}, {\rm \ for \ } 
k=1, \ldots, \mathfrak p_\star-1, \\
& \mathbbmss d_{n, k}\equiv\rba_{n, k+1}-\rba_{n, k} \to +\infty {\rm \ as \ } n \to +\infty, {\rm \ for \ } k=1, \ldots, \mathfrak p_\star-1, \\
&\underset {n \to + \infty} \liminf   \, \rbD\left(u_n, \rba_{n, k}\right)\geq \frac 14 \upeta_4  {\rm \ and \ there\ exists \ } R_\star>0, {\rm \ and \ } n_0 \in \N, {\rm \ such \ that \ } \\
&{\color{black}\rbD_\rL(u_n,s)\leq\frac12\upeta_4}, {\rm\ for\ }
s\in\R\setminus\underset{k=1}{\overset{\mathfrak p_\star}\cup}B(\rba_{n,k},R_\star),
{\rm\ and\ }n\geq n_0.
\end{aligned}
\right.
\end{equation}
\end{lemma}
   \begin{proof} The only statement which has not been proven yet is the last inequality.   We choose 
 \begin{equation}
  \label{rstar}
  R_\star=\sup\{R_k^\star, k=1,\ldots,  \mathfrak p_\star\}.
  \end{equation}
   {\color{black}The last inequality follows from \eqref{cestpetit}: outside the cluster intervals, the two neighboring integers are not concentration points, and their strips cover $\Gamma_\rL(s-1,s+1)$.}
     \end{proof}

\subsection{Convergence on compact sets}
\label{compacita}
\begin{proposition}
\label{local-compactness}
Assume that the potential $V$ satisfies assumptions ${\color{black}(\mathbf{H}_{\mathbf{1}})}$, ${\color{black}(\mathbf{H}_{\mathbf{2}})}$ and  ${\color{black}(\mathbf{H}_{\mathbf{3}})}$, and let $\rL>1$. Let  
$(v_n)_{n \in \N}$ in $\WL$ be a sequence satisfying \eqref{psoide}. Given  any $-\infty<a<b<+\infty$, there exists a subsequence of $(v_n)_{n \in \N}$, still denoted $(v_n)_{n \in \N}$ and $v \in \WL$ such that 
\begin{equation*}
\left\{
\begin{aligned}
&\rd \E(v)=0, {\rm i. e.} \  v {\rm \ is \ a \ weak \ solution \ of \eqref{themainthing}},  \\
& v_n \to  v {\rm \ strongly \  in \ }  H^1 (\Gamma_{\rL}(a, b)), {\rm \ as \ } n \to +\infty, {\rm \ and \ } \\
&\E(v_n, \Gamma_\rL (a, b)) \to \E(v, \Gamma_\rL (a, b)), {\rm \ as \ } n \to +\infty.
\end{aligned}
\right.
\end{equation*}
\end{proposition}

\begin{proof}  We divide the proof into steps.

\smallskip
\noindent
{\it Step 1: Local weak convergence in $\HL$.}
{\color{black}The energy bound and coercivity of $V$ give uniform $H^1$ bounds on every bounded strip. }Extracting possibly a subsequence,  still denoted $(v_n)_{n \in \N}$, we may assume that $(v_n)_{n \in \N}$ converges weakly in {\color{black}$H^1_{\rm loc}(\Lambda_\rL)$} to some map $v\in \ H^1 (\Gamma_{\rL}(c, d))$ as $n \to +\infty$, for any $-\infty <c<d<+\infty$.  Hence, by compact embedding $H^1 (\Gamma_{\rL}(c, d))\hookrightarrow  {\color{black}L^p(\Gamma_{\rL}(c,d))}$,  we  deduce that 
$v_n \to v$ strongly in ${\color{black}L^p(\Gamma_{\rL}(c,d))}$. Passing possibly to a further subsequence, we may assume that $v_n $ converges to $v$ pointwise for almost every $x \in \Lambda_\rL$.   {\color{black}Using the strong $L^p$ convergence for every finite $p$ and the polynomial growth bounds in Lemma \ref{technical}} we deduce that 
\begin{equation}
\label{Vn}
V(v_n) \to  V(v) {\rm \ in \ } \ L^1_{\rm loc} (\Lambda_\rL) {\rm \ and \ }   \nabla_uV(v_n) \to  \nabla_uV(v) {\rm \ in \ } \ L^1_{\rm loc} (\Lambda_\rL)
{\rm \ as \ } n  \to +\infty.
\end{equation}
{\color{black}Lower semicontinuity and Fatou's lemma give $\E(v)\leq M$, and periodicity passes to the limit; hence $v\in\WL$.} {\color{black}Passing to the limit in $\rd\E(v_n)(\varphi)$ for every smooth $\rL$-periodic test function $\varphi$ compactly supported in $x_2$ gives $\rd\E(v)(\varphi)=0$. By density, $\rd\E(v)=0$, so $v$ is a weak solution to \eqref{themainthing}.}

\medskip
\noindent
{\it Step 2:  Strong convergence on compact sets.} Consider the map $\rbv_n=v_{n,{\rm per}} $  defined on $\R^2$ using periodisation.  Let $a \in \R^2$  be given.   Applying Lemma \ref{sansdec0} to the translated map {\color{black}$\tilde\tau_{-a}(v_n)$} {\color{black}and translating back}, we deduce the decomposition
\begin{equation}
\label{decompon}
\left\{
\begin{aligned}
&\rbv_n=\virna+ \vregna {\rm \ on \ } \B(a, 1)  {\rm  \ such \ that, \ for  \ }  r>2 \\
&\Vert \vregna \Vert_{W^{2, r}(\B(a, \frac34))} 
\leq C(r) \left(M + M^{2m}+1 \right), \\
&\Vert \virna  \Vert_{H^1(\B(a, 1){\color{black})}}\leq C \Vert \rd E(v_n)\Vert_{H^\star} \to 0 {\rm \ as \ } n \to +\infty.
\end{aligned}
\right.
\end{equation}
In view of the compact embedding $W^{2, r}\left(\B(a, \frac34)\right)\hookrightarrow H^1\left(\B(a, \frac34)\right)$ and the third line in \eqref{decompon}, we deduce, passing possibly to a further subsequence,  that 
$\rbv_n$ {\color{black}converges} strongly to $v$ on the ball $\B(a, \frac34)$. Next consider a finite covering of $\Gamma_\rL(a,b)$ by a family of balls 
$\{\B(a_i, \frac34)\}_{i \in J}$, where $J$ is a finite set of indices.  Extracting possibly a further subsequence, we may assume that {\color{black}$(v_n)_{n\in\N}$}
{\color{black}converges} strongly to $v$ on each ball $\B(a_i, \frac34)$, hence strongly on $\Gamma_\rL(a, b)$.  {\color{black}Finally, the asserted energy convergence follows from the strong $H^1$ convergence and \eqref{Vn}.}
\end{proof}
\subsection{Extracting pseudo-profiles}
Here we go back to the situation described in Subsection \ref{generalcluster}, in particular Lemma \ref{clusterisation}.  For $k=1, \ldots, \mathfrak p_\star$,  consider the translated map $u_{n, k}$ defined on $\Lambda_\rL$ by 
\begin{equation}
\label{unk}
u_{n, k}(x_1, x_2)=u_n(x_1, x_2 +\rba_{n, k}), {\rm \ } x \in \Lambda_\rL, {\rm \ for \  } (x_1, x_2)\in \Lambda_\rL,   
\end{equation}
so that $u_{n, k}$ verifies \eqref{psoide}. As a consequence of Lemma \ref{clusterisation}, we have also 
\begin{equation}
\label{densiten}
\left\{
\begin{aligned}
&\rbD_\rL (u_{n, k}, 0)\geq\frac{1}{4} \upeta_4 ,{\rm \  for \ any \ }  n \in \N, {\rm \ and \ }\\
&\rbD_\rL(u_{n, k}, s) \leq {\color{black}\frac12\upeta_4}, {\rm \ for \ }  R_\star< \vert s \vert \leq R, 
 {\rm  \ for \ fixed \ } R   {\rm \ and  \ } n \ {\color{black}\text{sufficiently large.}}
\end{aligned}
\right.
\end{equation}

\begin{proposition} 
\label{thelastofus}
Assume that the potential $V$ satisfies assumptions ${\color{black}(\mathbf{H}_{\mathbf{1}})}$, ${\color{black}(\mathbf{H}_{\mathbf{2}})}$ and  ${\color{black}(\mathbf{H}_{\mathbf{3}})}$, and let $\rL>1$. Let  
$(u_n)_{n \in \N}$ in $\WL(\upsigma^-, \upsigma^+)$ be a sequence satisfying \eqref{psoide}.  Let $k \in \{1, \ldots, \mathfrak p_\star\}$, and let $(u_{n, k})_{n \in \N}$ be defined by \eqref{unk}. Passing possibly to a further subsequence, still denoted $(u_{n, k})_{n \in \N}$, there exists \emph{a non constant} finite energy solution $\rbu_k$ to \eqref{themainthing} such that $u _{n, k}$ converges weakly to  $\rbu_k$ in  $H^1_{\rm loc}(\Lambda_\rL)$ as  $n \to +\infty$ , and such that, for any $R>0$, $u _{n, k}$ converges strongly to {\color{black}$\rbu_k$} in $H^1_{\rm per}(\Gamma_{\rL}(-R, R))$ as $n \to +\infty$. 
\end{proposition}

\begin{proof} This is a direct consequence of {\color{black}Proposition \ref{local-compactness}} applied to the sequence $(u_{n, k})_{n \in \N}$.  The fact that 
$\rbu_k$ is not constant is a consequence of the first inequality of \eqref{densiten} and the strong convergence of $u_{n, k}$ towards $\rbu_k$ in $H^1\left(\Gamma_\rL(-1, 1)\right)$ as $n \to +\infty$, so that
$$
\rbD_\rL(\rbu_k, 0)\geq \frac{1}{4} \upeta_4, 
$$
and hence $\rbu_k$ cannot be a constant map.
\end{proof}

\subsection{Energy estimate off the cluster}

\begin{proposition}  
\label{dehors0}
Let $N \in \N^\star$. {\color{black}For $k=0,\ldots,\mathfrak p_\star$ (and $k=1,\ldots,\mathfrak p_\star-1$ in the separation inequality), there exists} $n_1(N)\in \N$ such that, if $n \geq n_1(N)$, then we have
\begin{equation}
\label{dehors}
\left\{
\begin{aligned}
& \mathbbmss d_{n, k}=\rba_{n, k+1}-\rba_{n, k} \geq  4(R_\star+N^2) {\rm \ and \ } \\
&\E\left(u_n, \Sigma (n, k, N)\right)\leq  \frac{\rbK_6 M}{N}+  \rbC_4 \Vert \rd E(u_n)\Vert_{H^\star}.
\end{aligned}
\right.
\end{equation}
where we have set  
\begin{equation}
\label{Sigmank}
\left\{
\begin{aligned}
\Sigma (n, k, N)&= \Gamma_\rL \left(\rba_{n, k}+(R_\star+ N^2),\  \rba_{n, k+1}-(N^2+R_\star)\right), {\rm \ for \ } k=1, \ldots, \mathfrak p_\star-1, \\
\Sigma (n, 0, N)&= \Gamma_\rL \left(-\infty, \  \rba_{n,1}-(N^2+R_\star)\right) {\rm \ and \ }\\
\Sigma (n, \mathfrak p_\star, N)&= \Gamma_\rL \left(\rba_{n, \mathfrak p_\star}+(R_\star+ N^2),\  {\color{black}+\infty}\right).
\end{aligned}
\right.
\end{equation}
\end{proposition}
\begin{proof} The first inequality of  \eqref{dehors} is a direct consequence of the second statement in \eqref{clusterisation0}. For the second we invoke Proposition \ref{caccioplon}, with $a=\rba_{n, k}+(R_\star+ N^2)$ and 
$b={\color{black}\rba_{n,k+1}}-(R_\star+ N^2)$.  {\color{black}In view of \eqref{clusterisation0} and $\Vert\rd\E(u_n)\Vert_{H^\star}\leq1$,} assumption \eqref{bonoloc} of Proposition \ref{caccioplon} is satisfied.  Hence inequality \eqref{ruegamma} holds, which yields  the second inequality in \eqref{dehors}.
{\color{black}For $k=0$ and $k=\mathfrak p_\star$, the same application uses $a=-\infty$ and $b=+\infty$, respectively.}
\end{proof}
\subsection{Proof of Proposition \ref{chainofsolutions}}
  So far, we have already defined $\mathfrak p_\star$,  the sequences of points $(\rba_{n, k})_{n \in \N}$  in $\Lambda_\rL$,  as well as the solution $\rbu_k$ of \eqref{themainthing}, so that  the first  required properties stated in Proposition \ref{chainofsolutions} are all established, in particular thanks to Lemma \ref{clusterisation} and Proposition \ref{thelastofus}.  {\color{black}It remains to establish} the convergence \eqref{chainofthought} and identity \eqref{theidentity}. 
  
  \begingroup\color{black}
Before proving the convergence, we establish \eqref{losange}. On each of the existing regions $\Sigma(n,k,1)$, $k=0,\ldots,\mathfrak p_\star$, Lemma \ref{locDL2} and Corollary \ref{EimpliesH} give a single well, with the squared $H^1$ norm of the difference from that well bounded by $\rbC_{\rm cap}M$. Every bounded region has height at least two for large $n$. Passing to a subsequence, these finitely many wells are independent of $n$; denote them by $\upsigma_k$. The exterior ones are $\upsigma_0=\upsigma^-$ and $\upsigma_{\mathfrak p_\star}=\upsigma^+$, since a nonzero constant is not square integrable on a half-cylinder.
For fixed $R_\star+1<S<T$, Proposition \ref{thelastofus} transfers the bound on $\Sigma(n,k,1)$ to $\int_{\Gamma_\rL(S,T)}|\rbu_k-\upsigma_k|^2\leq\rbC_{\rm cap}M$. Letting $T\to+\infty$ and using Lemma \ref{asympto2D0} identifies the right endpoint of $\rbu_k$ as $\upsigma_k$. The left endpoint is similarly $\upsigma_{k-1}$. This proves \eqref{losange}; in particular $\upsigma_k=\upsigma_k^+$ for $k\geq1$. Uniqueness in Lemma \ref{locDL2} shows that the well on every smaller region $\Sigma(n,k,N)$, $N\geq1$, is also $\upsigma_k$. A diagonal subsequence common to the finitely many profiles gives local strong convergence on every bounded shifted strip.
\endgroup
\subsubsection {Proof of the convergence \eqref{chainofthought}}
   We set 
  \begingroup\color{black}
\begin{equation}
\label{gn}
\begin{aligned}
g_n&=u_n-\rbU_n,\qquad
\rbu_{k,n}=\rbu_k(\cdot-\rba_{n,k}\rbe_2),\\
\rbU_n&=\sum_{k=1}^{\mathfrak p_\star}\rbu_{k,n}
-\sum_{k=1}^{\mathfrak p_\star-1}\upsigma_k^+.
\end{aligned}
\end{equation}
\endgroup

  {\color{black}The endpoint matching and Lemma \ref{affinage} imply $g_n\in\HL$.} We have to show that the $H^1$ norm  of $g_n$ can be made arbitrarily small, that is, given $\eps>0$, there exists $n_\eps \in \N$ such that if $n\geq {\rm \bf n}_\eps$, then 
  \begin{equation}
  \label{Egn}
  \Vert g_n \Vert_{H^1(\Lambda_\rL)}^2\leq \eps, {\rm \ for \  } n \geq \rbn_\eps.
\end{equation}
We decompose $\Lambda_\rL$, for $N$ to be determined in the course of the proof, as 
\begingroup\color{black}
\begin{equation}
\label{partition}
\left\{\begin{aligned}
\Lambda_\rL&=\left(\bigcup_{k=1}^{\mathfrak p_\star}\mathcal C_{k,n}\right)\cup\Sigma(n,N),\\
\mathcal C_{k,n}&=\Gamma_\rL(\rba_{n,k}-(R_\star+N^2),\rba_{n,k}+R_\star+N^2),\\
\Sigma(n,N)&=\bigcup_{k=0}^{\mathfrak p_\star}\Sigma(n,k,N).
\end{aligned}\right.
\end{equation}
\endgroup

We will call the sets $\mathcal C_{k, n}$ the cluster sets centered at ${\color{black}\rba_{n, k}}$.
We divide the proof of \eqref{Egn} in several steps.

\smallskip
\noindent
{\it Step 1: Vanishing of energy between the clusters}.  We apply Proposition \ref{dehors0} with  some integer $N$ to be determined in the course of the proof, but verifying, at this stage the lower bound
\begin{equation}
\label{Nlarge}
\frac{\rbK_6 M}{N} \leq \frac{\eps}{64\rbC_{\rm cap}(\mathfrak p_\star+1)},
\end{equation}
where {\color{black}$\rbC_{\rm cap}>0$} is the constant appearing in Corollary \ref{EimpliesH}.  Since, by assumption 
$\displaystyle{\Vert \rd\E(u_n ) \Vert_{H^\star}}$ tends to $0$ as {\color{black}$n\to+\infty$}, 
there exists some $\rbn_{\eps, 1}>0$ such that 
\begin{equation}
\label{debacle}
\rbC_4 \Vert d\E(u_n ) \Vert_{H^\star}\leq \frac{\eps}{64\rbC_{\rm cap}(\mathfrak p_\star+1)}, {\rm \ if \ } n \geq \rbn_{\eps, 1}.
\end{equation}
Proposition \ref{dehors0} then yields, if $n \geq  \sup \{ n_{\eps, 1}, {\color{black}n_1(N)}\}$
\begin{equation}
\label{deroute}
\begin{aligned}
\E\left(u_n, \Sigma (n, k, N)\right)&\leq  \frac{\rbK_6 M}{N}+  \rbC_4 \Vert \rd E({\color{black}u_n})\Vert_{H^\star} \\
& \leq \frac{\eps}{32\rbC_{\rm cap}(\mathfrak p_\star+1)}.
\end{aligned}
\end{equation}
We next invoke Corollary  \ref{EimpliesH}, to  deduce that
\begin{equation}
\int_{\Sigma (n, k, N)} \vert \nabla u_n\vert^2+ \vert u_n{\color{black}-\upsigma_k}\vert^2
\leq \frac{\eps}{32(\mathfrak p_\star+1)}, 
\end{equation}
and hence,  considering the  $\Sigma(n, N)=\underset {k=0}{\overset{\mathfrak p_\star}\cup} \Sigma (n, k, N)$ defined in \eqref{partition}, we obtain, if $N$ satisfies \eqref{Nlarge} and $n \geq \rbn_{\eps, 1}$ given in \eqref{debacle}
\begingroup\color{black}
\begin{equation}
\label{judo}
\sum_{k=0}^{\mathfrak p_\star}\int_{\Sigma(n,k,N)}
\left(|\nabla u_n|^2+|u_n-\upsigma_k|^2\right)\leq\frac{\eps}{32}.
\end{equation}
\endgroup

\smallskip
\noindent
{\it Step 2: Bounds on  the solutions $\rbu_k$ off the clusters}.   In view of Corollary \ref{decay} and Remark \eqref{decadence01}, there exists  $\tilde R>0$, such that, for any $k=1, \ldots, \mathfrak p_\star$,  we have, for $N^2\geq 2\tilde R$, 
\begingroup\color{black}
\begin{equation}
\label{decadence2}
\left\{\begin{aligned}
&|\rbu_k(x)-\upsigma_k^\pm|^2\leq\frac{4\rbK_7}{N^4},
\quad\hbox{for }\pm x_2\geq N^2,\\
&\E\left(\rbu_k,\Lambda_\rL\setminus\Gamma_\rL(-N^2,N^2)\right)
\leq\frac{8\rbK_6\upeta_3}{N^4}.
\end{aligned}\right.
\end{equation}
\endgroup

If we choose $N$ so large that \eqref{Nlarge} is satisfied and such that moreover, we have
\begin{equation}N^2\geq2\tilde R,\end{equation}
\begingroup\color{black}
then Corollary \ref{EimpliesH}, applied separately to the two tails and their respective wells, gives
\[
\begin{aligned}
&\int_{\Gamma_\rL(N^2,+\infty)}\left(|\nabla\rbu_k|^2+|\rbu_k-\upsigma_k^+|^2\right)\\
&\quad+\int_{\Gamma_\rL(-\infty,-N^2)}\left(|\nabla\rbu_k|^2+|\rbu_k-\upsigma_k^-|^2\right)
\leq\frac{8\rbC_{\rm cap}\rbK_6\upeta_3}{N^4}.
\end{aligned}
\]
Here $N$ is also chosen large enough for the low-density hypothesis on these tails. On $\Sigma(n,k,N)$, the corrected sum satisfies
\[
\rbU_n-\upsigma_k
=\sum_{j\leq k}(\rbu_{j,n}-\upsigma_j^+)
+\sum_{j>k}(\rbu_{j,n}-\upsigma_j^-).
\]
Each summand lies on the appropriate tail. The regions $\Sigma(n,k,N)$ are disjoint, so the squared norm of a sum of at most $\mathfrak p_\star$ terms gives
\[
\sum_{k=0}^{\mathfrak p_\star}\int_{\Sigma(n,k,N)}
\left(|\nabla\rbU_n|^2+|\rbU_n-\upsigma_k|^2\right)
\leq\frac{8\mathfrak p_\star^2\rbC_{\rm cap}\rbK_6\upeta_3}{N^4}.
\]
\endgroup
We then impose that $N$ is sufficiently large so that
{\color{black}$8(\mathfrak p_\star+1)^3\rbC_{\rm cap}\rbK_6\upeta_3/N^4\leq\eps/16$}, and hence
\begingroup\color{black}
\begin{equation}
\label{majorunef}
\sum_{k=0}^{\mathfrak p_\star}\int_{\Sigma(n,k,N)}
\left(|\nabla\rbU_n|^2+|\rbU_n-\upsigma_k|^2\right)\leq\frac{\eps}{16}.
\end{equation}
\endgroup
\noindent
\smallskip
{\it Step 3: Estimates off  the clusters completed}. Combining \eqref{judo} with \eqref{majorunef} we are led to 
\begin{equation}
\label{zombie}
\int_{\Sigma (n, N)} \vert\nabla (u_n-  \rbU_{n})\vert^2 +
 \vert  u_n-  \rbU_{n}\vert^2 \leq \frac{\eps}{4}.
\end{equation}

\noindent
\smallskip
{\it Step 4: Estimates near the clusters}.   For $k={\color{black}1}, \ldots, \mathfrak p_\star$, with $N$ determined previously, we consider the  clusters 
$\mathcal C_{k, n}={\color{black}\Gamma_\rL} (\rba_{n, k}-(R_\star+N^2), \rba_{n, k}+R_\star+N^2)$.  It follows from Proposition \ref{thelastofus}, {\color{black}after translating $\mathcal C_{k,n}$ by $-\rba_{n,k}\rbe_2$ to the fixed strip}
$$
\int_{\mathcal C_{k,n}} \vert \nabla (u_n-\rbu_{k,n})\vert^2+ \vert u_n-\rbu_{k,n}\vert^2\to 0, {\rm \ as \ } n \to +\infty.
$$
Hence, we may find some $n_{\eps, 2}$ such that, if $n\geq n_{\eps, 2}$, then 
\begin{equation}
\label{waouhh}
\int_{\mathcal C_{k,n}} \vert \nabla (u_n-\rbu_{k,n})\vert^2+ \vert u_n-\rbu_{k,n}\vert^2
\leq \frac{\eps}{16(\mathfrak p_\star+1)},
\end{equation}
Finally, we consider {\color{black}$\tilde{\rbU}_{k,n}=\rbU_n-\rbu_{k,n}
=\sum_{j<k}(\rbu_{j,n}-\upsigma_j^+)+\sum_{j>k}(\rbu_{j,n}-\upsigma_j^-)$.
}Arguing as in Step 2, we have
$$
\int_{\mathcal C_{k,n}} \vert\nabla  \tilde \rbU_{k, n}\vert^2 +
 \vert  \tilde  \rbU_{k, n}\vert^2 {\color{black}\leq} \frac{\eps}{16(\mathfrak p_\star+1)}
$$
so that
\begin{equation}
\label{zumo}
\underset{k={\color{black}1}} {\overset{\mathfrak p_\star}\sum} \int_{\mathcal C_{k,n}} \vert\nabla  \tilde \rbU_{k, n}\vert^2 +
 \vert  \tilde  \rbU_{k, n}\vert^2 \leq \frac{\eps}{16}.
\end{equation}
Combining \eqref{zumo} with \eqref{waouhh}, we are led, for $n \geq n_{\eps, 2}$  to
\begin{equation}
\label{zumo2}
\underset{k={\color{black}1}} {\overset{\mathfrak p_\star}\sum} \int_{\mathcal C_{k,n}} \vert\nabla {\color{black}(u_n-\rbU_n)}\vert^2 +
 \vert  u_n-  \rbU_{n}\vert^2 \leq \frac{\eps}{4}.
\end{equation}

\smallskip
\noindent
{\it Step {\color{black}5} : Adding the estimates}. Combining {\color{black}\eqref{zumo2}} with \eqref{zombie} we obtain the desired estimate \eqref{Egn}, {\color{black}taking $n_\eps\geq\max\{n_1(N),n_{\eps,1},n_{\eps,2}\}$.}

\subsubsection{Proof of identity \eqref{theidentity}}
\label{lidentite}
 We introduce again an arbitrary $\eps>0$. We are going to prove that, for $n \geq n_{\eps, 3}$ to be determined, we have
 \begin{equation}
 \label{ecremage}
 \left \vert \E(u_n)-\underset{k=1}{\overset{\mathfrak p_\star}\sum} \E(\rbu_k)\right \vert \leq {\color{black}\eps}.
\end{equation}
We divide again the proof into steps.

\smallskip
\noindent
{\rm \it Step 1: Energy off  the cluster}. \begingroup\color{black}
We use Proposition \ref{dehors0} directly. Choose
\begin{equation}
\label{Nlarge2}
N\geq\rbN_{\eps,2}=\frac{64(\mathfrak p_\star+1)\rbK_6M}{\eps},
\end{equation}
and, for this fixed $N$, take $n\geq n_1(N)$ so large that
$\rbC_4\Vert\rd\E(u_n)\Vert_{H^\star}\leq\eps/[64(\mathfrak p_\star+1)]$.
Adding \eqref{dehors} over $k=0,\ldots,\mathfrak p_\star$ gives
\begin{equation}
\label{albanie}
\E(u_n,\Sigma(n,N))\leq\frac{\eps}{32}.
\end{equation}
\endgroup
\smallskip
\noindent
{\rm \it Step 2: Energy near the cluster}.     {\color{black}In} view of the continuity of the energy functional, for {\color{black}each fixed $N$ and $k=1,\ldots,\mathfrak p_\star$}, we have 
\begingroup\color{black}
\begin{equation}
\label{cantoche}
\E(u_n,\mathcal C_{k,n})\longrightarrow
\E\left(\rbu_k,\Gamma_\rL(-R_\star-N^2,R_\star+N^2)\right)
\quad\hbox{as }n\to+\infty.
\end{equation}
\endgroup

On the other hand 
$$
\E(\rbu_k, \Gamma_\rL(-R_\star{\color{black}-}N^2, R_\star+ N^2){\color{black})} \to \E(\rbu_k), {\rm \ as \ } N \to +\infty.
$$
Hence, there exists some $\rbN_{\eps, 3}\in \N$, depending only on $\eps$,  such that 
\begingroup\color{black}
\begin{equation}
\label{canopee}
0\leq\E(\rbu_k)-\E\left(\rbu_k,\Gamma_\rL(-R_\star-N^2,R_\star+N^2)\right)
\leq\frac{\eps}{8\mathfrak p_\star},\quad N\geq\rbN_{\eps,3}.
\end{equation}
\endgroup

{\color{black}Combining \eqref{cantoche} with \eqref{canopee}, for each fixed $N\geq\rbN_{\eps,3}$ there exists $\rbn_{\eps,2}(N)$ such that, for $n\geq\rbn_{\eps,2}(N)$,}
\begingroup\color{black}
\begin{equation}
\label{greece}
\left|\sum_{k=1}^{\mathfrak p_\star}\E(u_n,\mathcal C_{k,n})
-\sum_{k=1}^{\mathfrak p_\star}\E(\rbu_k)\right|\leq\frac{\eps}{4},
\quad N\geq\rbN_{\eps,3},\quad n\geq\rbn_{\eps,2}(N).
\end{equation}
\endgroup

\smallskip
\noindent
{\rm \it Step {\color{black}3}: Proof of \eqref{theidentity} completed}.  {\color{black}We first fix an integer $N\geq\sup\{\rbN_{\eps,2},\rbN_{\eps,3}\}$, and then choose $n_{\eps,3}$ large enough for both \eqref{albanie} and \eqref{greece} to hold.} Combining \eqref{greece} and \eqref{albanie}, we obtain the desired result for $n \geq \rbn_{\eps, 3}$. {\color{black}This} completes the proof of Proposition \ref{chainofsolutions}.
\section{One-dimensional profiles}
\label{one-d}

\subsection{First properties}

In this section, we review some useful facts concerning one-dimensional profiles, {\color{black}that is \emph{finite energy solutions to \eqref{heterocline}}.} Most of the   results  found in this section are not new and can be found in {\color{black}\cite{alabrogui, schatzman, bonafoux, rabinowitz}}{\color{black}. Consider the problem}
\begin{equation}
\label{heterocline20}
\left\{
\begin{aligned}
	\frac{{\rm d}^2 \rbq }{{\rd} s^2}(s)&=\nabla_w V(\rbq (s)) {\rm \ for  \ } s \in  \R, \\
	\rbq(\pm \infty)&=\upsigma^\pm, {\rm \ where \ } \upsigma^\pm\in \Sigma, 
\end{aligned}
\right.
\end{equation} 
so that ${\color{black}\rbq} \in \mathbb X(\upsigma^-, \upsigma^+)$. Notice that \emph{all results} in Sections \ref{functframe} and  \ref{propsol} transfer directly to one-dimensional profiles: It suffices to consider the function $\mathbbmss w$ defined on $\Lambda_2$ by
$\mathbbmss w(x_1, x_2)={\color{black}\rbq(x_2)}$, for $x_2\in \R$, $\displaystyle{x_1\in [-1,1]}$.  For instance, applying Proposition \ref{lwb0} to the map  $\mathbbmss  w$ yields, if ${\color{black}\rbq}$ is not constant
\begin{equation}
\label{weirdos}
\E({\color{black}\rbq}, \R) \geq \frac{\upeta_{\min}}{2}.
\end{equation}
Likewise, we obtain a regularity result by applying Lemma \ref{regu}.

\begin{lemma}
\label{regu1D}
Let $\rbq$ be a solution of \eqref{heterocline20}. We have, 
\begin{equation}
\label{smoothy-1d}
\vert \rbq(s)\vert+\vert  \rbq'  (s)\vert+ \vert \rbq'' (s)\vert  \leq \rbC_{\rm nd}, \forall  s \in \R, 
\end{equation}
where $\rbC_{\rm nd}$ is the constant provided in Lemma \ref{regu}.
\end{lemma}

\subsection{Some notation}
 We set 
\begin{equation}
\label{mathcalZ}
\mathcal Z(\upsigma^-, \upsigma^+)=\{w \in \mathbbmss X(\upsigma^-, \upsigma^+), w  {\rm \ is \ a \ minimizing \  solution \ to \ } \eqref{heterocline20}\}.
\end{equation}
As seen in many places, if $w$ is a solution to \eqref{heterocline20}, then so is the translated map $w_a=\tau_a(w)$ defined by
$\tau_a(w)(s)=w(s-a)$, for $a\in \R$, with $\E(w_a, \R)=\E(w, \R)$.   Recall that we have introduced in \eqref{equivclass}  the equivalence class
$$\mathcal C(w)=\{z \in \mathbbmss X(\upsigma^-, \upsigma^+), {\rm \ such \ that \ } \exists a\in \R, {\rm \ such \ that \ } 
z(x)=w(x-a), \forall x \in \R\}
$$
 If $F$  and $G$ are two {\color{black}subsets} of $\mathbbmss X(\upsigma^-, \upsigma^+)$, we set, for $q \in \mathbbmss X(\upsigma^-, \upsigma^+)$, and $s\in \{0, 1\}$
 \begin{equation*}
 \begin{aligned}
 {\rm  \bf dist}_s(F, G)&=\inf\{ \Vert q-v \Vert_{H^s(\R)}, q \in G,  v \in F\} {\rm \ and \ }\\
 {\rm  \bf dist}_s(q, F)&={\rm  \bf dist}_{{\color{black}s}}(\{q\}, F)=\inf\{ \Vert v -q\Vert_{H^s(\R)},  v \in F\}, 
 \end{aligned}
 \end{equation*}
 with the usual convention $H^0(\R)=L^2(\R)$. 
 Given $v \in \mathbbmss X(\upsigma^-, \upsigma^+)$,  and $w\in  \mathcal Z(\upsigma^-, \upsigma^+)$, the quantity 
 $\displaystyle{ \dist_s\left(v, \mathcal C(w)\right)}$ plays a distinguished role in the rest of the discussion.   For $v$ and $w$ two elements in
 $\mathbbmss X(\upsigma^-, \upsigma^+)$, we set
 \begin{equation}
 \label{orbitedist}
 \mathfrak D_s (v, w)=\dist_s (v, \mathcal C(w))=\dist_s (w, \mathcal C(v))=\dist_s (\mathcal C(v), \mathcal C(w)).
 \end{equation}

\begin{lemma} 
\label{separatiste}
Let  $v$ and $w$ be two elements in $\mathbbmss X(\upsigma^-, \upsigma^+)$. We have, for $s=0, 1$
\begin{equation} 
\mathfrak D_s (v, w)=0 {\rm \ if \ and \ only \ if \ } v \in C(w),  \left({\rm or, equivalently, \ } w \in \mathcal C(v)\right).
\end{equation}
\end{lemma} 
 \begin{proof} We have 
 $$\mathfrak D_s (v, w)=\underset{a \in \R} \inf F_s(a) {\rm  \ where \ }  F_s(a)=\Vert v(\cdot)-w(\cdot-a) \Vert_{H^s(\R)}, {\rm \ for \ }
 a \in \R.
  $$
\begingroup\color{black} The real-valued function $F_s$ is continuous on $\R$, and $F_s(a)\to+\infty$ as $a\to\pm\infty$ if $\upsigma^+\not=\upsigma^-$. If $\upsigma^+=\upsigma^-=\upsigma$, then
\[
F_s(a)^2\longrightarrow\Vert v-\upsigma\Vert_{H^s}^2+\Vert w-\upsigma\Vert_{H^s}^2
\quad\hbox{as }|a|\to\infty,
\]
since the scalar product of a fixed $H^s$ function with translates of another tends to zero. Thus, if $\mathfrak D_s(v,w)=0$ and $v,w$ are not both the same constant, a minimizing sequence of shifts is bounded. Passing to a subsequence and using continuity gives some $a_0\in\R$ such that $v=\tau_{a_0}(w)$. The constant case and the converse are immediate.
\endgroup \end{proof}
   Given $\uprho >0$ and $v \in \mathbbmss X(\upsigma^-, \upsigma^+)$, 
   we introduce, for $s=0, 1$,  the neighborhood $\mathcal  N_s(v, \uprho)$  defined by 
   $$
   \mathcal N_s(v, \uprho)=\{w \in \mathbbmss X(\upsigma^-, \upsigma^+), {\rm \ such \ that \ } \dist_s(w, \mathcal C(v))\leq \uprho\}.
   $$
  
   \begin{lemma} 
   \label{glisko}
 Let $s \in \{0, 1\}$ and   let  $v$ and $w$ be two elements in $\mathbbmss X(\upsigma^-, \upsigma^+)$ such that $\mathfrak D_s(v, w)\not=0$. If $0<\uprho <\frac12 \mathfrak  D_s(v, w) $, then, we have 
   \begingroup\color{black}
\begin{equation}
\label{doesnothold}
\left\{
\begin{aligned}
&\dist_s\bigl(\mathcal N_s(v,\uprho),\mathcal N_s(w,\uprho)\bigr)
+2\uprho\geq\mathfrak D_s(v,w),\\
&\mathcal N_s(v,\uprho)\cap\mathcal N_s(w,\uprho)=\emptyset
\qquad\hbox{if }0<\uprho<\tfrac12\mathfrak D_s(v,w).
\end{aligned}
\right.
\end{equation}
\endgroup
   If {\color{black}$\tilde v\in\mathcal N_s(w,\uprho)$ and $\tilde w\in\mathcal N_s(v,\uprho)$}, then, we have
  \begingroup\color{black}
\begin{equation}
\label{ultimatum}
\Vert\tilde w-\tilde v\Vert_{H^s(\R)}
\geq\dist_s\bigl(\mathcal N_s(v,\uprho),\mathcal N_s(w,\uprho)\bigr)
\geq\mathfrak D_s(v,w)-2\uprho.
\end{equation}
\endgroup
   \end{lemma}
  \begingroup\color{black}
\begin{proof}
Let $f\in\mathcal N_s(v,\uprho)$, $g\in\mathcal N_s(w,\uprho)$ and $\delta>0$.
Choose $z_v\in\mathcal C(v)$ and $z_w\in\mathcal C(w)$ such that
$\Vert f-z_v\Vert_{H^s(\R)}\leq\uprho+\delta$ and
$\Vert g-z_w\Vert_{H^s(\R)}\leq\uprho+\delta$.
The triangle inequality for the norm gives
\begin{equation}
\mathfrak D_s(v,w)\leq\Vert z_v-z_w\Vert_{H^s(\R)}
\leq\Vert f-g\Vert_{H^s(\R)}+2\uprho+2\delta.
\end{equation}
Letting $\delta\to0$ and taking the infimum over $f$ and $g$ proves the first inequality in \eqref{doesnothold}.
Since $\mathfrak D_s(v,w)-2\uprho>0$, the neighborhoods are disjoint.
Finally, \eqref{ultimatum} follows from the definition of the distance between sets.
\end{proof}
\endgroup

\subsection{Some consequences of assumption ${\color{black}(\mathbf{H}_{\mathbf{4}})}$}
\label{H4}
The statement of assumption  ${\color{black}(\mathbf{H}_{\mathbf{4}})}$  can be stated as 
\begin{equation}
\label{H4-1}
\mathcal Z(\upsigma_0^-, \upsigma_0^+)=\mathcal  C(w_1)\cup \mathcal C(w_2), {\rm \ with \ }  
\mathcal  C(w_1)\cap \mathcal C(w_2)=\emptyset, 
\end{equation}
where $w_1$ and $w_2$ are two distinct \emph{minimizing} solutions to \eqref{heterocline20}.   
Consider  the function $\Uptheta: \R^+\to \R^+$, defined for $\uprho\in \R^+$, by 
 \begingroup\color{black}
\begin{equation}
\label{Uptheta1}
\begin{aligned}
\Uptheta(\uprho)
&=\inf\left\{\E(v,\R):\,v\in\mathbbmss X(\upsigma_0^-,\upsigma_0^+),\
\dist_1\bigl(v,\mathcal Z(\upsigma_0^-,\upsigma_0^+)\bigr)\geq\uprho\right\}\\
&=\inf\left\{\E(v,\R):\,v\in\mathbbmss X(\upsigma_0^-,\upsigma_0^+),\
\dist_1(v,\mathcal C(w_i))\geq\uprho,\ i=1,2\right\}.
\end{aligned}
\end{equation}
\endgroup
The following extends the strict-gap statement in \cite{schatzman}, Corollary 3.2:

\begin{lemma}  
\label{Uptheta0}
Assume that  ${\color{black}(\mathbf{H}_{\mathbf{1}})}$, ${\color{black}(\mathbf{H}_{\mathbf{2}})}$,   ${\color{black}(\mathbf{H}_{\mathbf{3}})}$ and ${\color{black}(\mathbf{H}_{\mathbf{4}})}$ hold.   The function $\Uptheta$ is non-decreasing, with $\Uptheta(0)=\bm(\upsigma_0^-, \upsigma_0^+)$. Moreover, we have 
 \begin{equation}
 \label{Uptheta2}
 \Uptheta (\uprho )>\mathbbmss m (\upsigma_0^-, \upsigma_0^+), {\rm \ for \ any \ } \uprho >0.\\
  \end{equation}
 \end{lemma}
\begin{proof}
Monotonicity and the value at zero follow from the definition.
If \eqref{Uptheta2} failed, there would be a minimizing sequence $(v_n)$ with
$\dist_1(v_n,\mathcal Z(\upsigma_0^-,\upsigma_0^+))\geq\uprho>0$.
By Ekeland's principle on $v_0+H^1(\R)$, after an $o(1)$ perturbation in $H^1$
we may assume that the one-dimensional differentials tend to zero and the distance remains at least $\uprho/2$.
Regard this sequence as maps on $\Lambda_2$ independent of $x_1$.
Averaging test functions in $x_1$ shows that it is also a Palais--Smale sequence there.
Proposition \ref{chainofsolutions} gives a chain of nonconstant profiles,
each independent of $x_1$ as a local weak limit of translated sequence elements.
Dividing the energy identity by $2$, their one-dimensional energies sum to
$\bm(\upsigma_0^-,\upsigma_0^+)$.
Deleting profiles with equal endpoints leaves a chain joining $\upsigma_0^-$ to $\upsigma_0^+$,
whose energy is at least this minimum by \eqref{v3}; each deleted profile would have positive energy.
Thus none is deleted, and \eqref{strictosenso} forces a single profile.
It is a minimizer, hence belongs to $\mathcal Z(\upsigma_0^-,\upsigma_0^+)$.
The strong $H^1$ decomposition in Proposition \ref{chainofsolutions} then gives
$\dist_1(v_n,\mathcal Z(\upsigma_0^-,\upsigma_0^+))\to0$, a contradiction.
\end{proof}

It is sometimes interesting to renormalize $\E$ and $\Uptheta$, introducing, for $\uprho>0$
\begin{equation}
\label{renormalized}
\tilde \Uptheta(\uprho)=\Uptheta (\uprho)-\bm {\rm \ and \ }
\tilde \E(v, \R)=\E(v, \R)-\bm, {\rm \ for \ } v \in \mathbbmss X(\upsigma^-, \upsigma^+), 
\end{equation}
 so that the  function $\tilde \Uptheta$ is non-decreasing, with $\tilde \Uptheta(0)=0$.
Lemma \ref{Uptheta0} can then be cast as follows:

\begin{corollary} 
\label{voisinage0}
Assume that  ${\color{black}(\mathbf{H}_{\mathbf{1}})}$, ${\color{black}(\mathbf{H}_{\mathbf{2}})}$,   ${\color{black}(\mathbf{H}_{\mathbf{3}})}$ and ${\color{black}(\mathbf{H}_{\mathbf{4}})}$ hold.   For any  
$\uprho>0$, if  $v \in \mathbbmss  X(\upsigma_0^-, \upsigma_0^+)$ satisfies the energy bounds
\begin{equation}
\label{voisinage1}
 \tilde \E(v, \R) < \tilde \Uptheta (\uprho), 
 \end{equation}
then, for some $i\in \{1, 2\}$,  we have $v \in \mathcal N_1(w_i, \uprho)$, i.e.
$\dist_1( v, \mathcal C(w_i))\leq \uprho$.
\end{corollary}
\subsection{Some consequences of assumption ${\color{black}(\mathbf{H}_{\mathbf{5}})}$}
 \subsubsection{The spectral gap property}
 Given $v \in \mathbb X(\upsigma_0^-, \upsigma_0^+)$ a solution to \eqref{heterocline}, we consider in this paragraph the  second derivative $\rd^2\E(v)$ of $\E$ at $v$ and the  linearized operator  $\rbA(v)$ of \eqref{heterocline} near $v$  defined by 
  \begin{equation}
  \label{linearized}
  \rbA(v)(\rbh)=-\frac{\rd^2 \rbh}{\rd s^2}+\nabla^2 V(v)\rbh, {\rm \ for \ } \rbh \in {\color{black}H^2(\R)}.
 \end{equation}
{\color{black}The same formula defines a weak operator from $H^1(\R)$ to $H^{-1}(\R)$.} They are related by the identity for $\rbh, \rbk$ in    $H^2(\R)$ 
\begin{equation}
\label{secondo}
 \langle \rbA(v) \rbh, \rbk \rangle_{L^2(\R)}= \rd^2\E({\color{black}v})(\rbh, \rbk)=
 \int_{\R}  \rbh'(s)\cdot \rbk'(s)+\nabla^2V (v(s)) (\rbh(s), \rbk(s))\rd s, 
\end{equation}
so that $\rd^2\E(v)$ is a symmetric quadratic form, and $\rbA(v)$ is self-adjoint.  As already  seen, if $v$ is a minimizing solution, then $\rd^2\E$ is positive, that is 
\begin{equation}
\label{pasdef}
\rd^2\E(v)(\rbh, \rbh)={\color{black}\langle \rbA(v) \rbh, \rbh \rangle_{H^{-1},H^1}} \geq 0, {\rm \ for \  any\ }  \rbh\in H^1(\R).
\end{equation}

If \eqref{pasdef}  is {\color{black}satisfied}, then  we  have 
\begin{equation}
\label{cooloss}
\ker \rbA(v)\equiv \{\rbh \in H^1(\R), {\rm {\color{black}\ s.t.\ }} \rbA(v)(\rbh)=0\} =\{ \rbh \in H^1(\R), {\rm \ s.t. \ }  \rd^2\E(v) (\rbh, \rbh)=0\}.
\end{equation}
 Indeed, if \eqref{pasdef} is {\color{black}satisfied}, we may apply the Cauchy-Schwarz inequality for the quadratic form $\rd^2\E$,  so that,  if 
 $\rd^2\E(v)(\rbh, \rbh)=0$,  we have, for any
 $\rbk \in H^1(\R)$, 
$$
{\color{black}|\rd^2\E(v)(\rbh,\rbk)|\leq
\bigl(\rd^2\E(v)(\rbh,\rbh)\bigr)^{1/2}
\bigl(\rd^2\E(v)(\rbk,\rbk)\bigr)^{1/2}=0,}
$$
 and hence  $ \rbA(v) \rbh=0$. Moreover,  since solutions $v \in \mathbb X(\upsigma_0^-, \upsigma_0^+)$ to \eqref{heterocline} are uniformly bounded in view of 
 Lemma  \ref{maximodo}, there exists a constant $\rbC_{\rm nd}>0$ such that, for {\color{black}any $s\in\R$ and $z\in\R^k$}, we have
 \begin{equation}
 \label{sec}
\left   \vert \nabla^2 V(v(s))(z, z)\right \vert \leq \rbC_{\rm nd}{\color{black}|z|^2}.
   \end{equation}
 Notice that the  elements of the  Kernel $\ker \rbA(v)$ are $C^2$ functions, by regularity theory for differential {\color{black}equations}. We consider its orthogonal for the $L^2$ scalar product, i.e.
 $$
 \ker \rbA(v)^\perp=\{\rbh \in H^1(\R), {\rm \ such \ that \ } \langle \rbh, z\rangle_{L^2(\R)}=0, {\rm \ for \ any \ } z \in \ker \rbA\}.
 $$
   
 \begin{proposition}
 \label{spectralou}
  Let $w\in \mathbb X(\upsigma_0^-, \upsigma_0^+)$ be a minimizing solution to \eqref{heterocline}. There exists a number 
  $\alpha_\star(w)>0$ such that 
\begin{equation}
 \label{claimL200}
 \rd^2\E(w)(\rbh, \rbh)\geq  \alpha_\star(w) \Vert \rbh \Vert_{H^1(\R)}^2, {\rm \ if \ } \rbh \in \ker\rbA^\perp.
 \end{equation}
 \end{proposition}
\begin{proof}    
Since we have invariance by translation, we  cannot invoke the standard results for  compact self-adjoint {\color{black}operators}. However, we are able to adapt the proof, in particular thanks to \eqref{potto}, which yields 
$$  \nabla_\rby^2 V(\rby) \geq \frac{\lambda_0^-}{2}{\rm Id}_{\R^k}  {\rm \ for  \ }  \rby \in {\color{black}\B(\upsigma_0^-,\upmu_0) \cup\B(\upsigma_0^+, \upmu_0)}. $$
Hence, there exists some $S_0>0$ such that   $w(s)\in \B(\upsigma_0^-, \upmu_0)$, for $s \leq -S_0$, $w(s)\in \B(\upsigma_0^+, \upmu_0)$, for $s \geq S_0$, so that  
\begin{equation}
\label{humus}
\left\{
\begin{aligned}
& \nabla_\rby^2 V(w(s)) \geq \frac{\lambda_0^-}{2}{\rm Id}_{\R^k},{\rm \ for \ } s \in \R\setminus [{\color{black}-}S_0, S_0],
 {\rm \ and  \ therefore\ } \\
& \nabla_\rby^2 V(w(s))(\rbh, \rbh) \geq \frac{\lambda_0^-}{2}\vert \rbh\vert^2, {\rm \ for \ any \ } \rbh \in \R^k, {\rm \ and \ } s \in \R\setminus [{\color{black}-}S_0, S_0].
\end{aligned}
\right.
\end{equation}
We first prove a weaker version of \eqref{claimL200}

\smallskip
\noindent
{\it Step 1}. {\color{black}\hspace{0pt}}{\it There exists a number $ \alpha_0>0$ such that 
 \begin{equation}
 \label{claimL20}
 \rd^2\E(w)(\rbh, \rbh)\geq  \alpha_0 \Vert \rbh \Vert_{L^2(\R)}^2, {\rm \ if \ } \rbh \in \ker \rbA^\perp. 
 \end{equation}
 }
 \noindent
 {\it Proof of \eqref{claimL20}}.
  We argue by contradiction and  assume that 
 there is a sequence $(\rbh_n)_{n \in \N}$  of elements of $\ker \rbA^\perp$ such that
 \begin{equation}
 \label{oignon}
 \Vert \rbh_n \Vert_{L^2({\color{black}\R})}=1, 
 \rd^2\E(w) (\rbh_n, \rbh_n)\underset { n \to \infty} \to 0,  {\rm \ and \ } \rbh_n \in (\ker A)^\perp, {\rm \ for \ any \ } n \in \N.
 \end{equation}
\begingroup\color{black}
Since the quadratic energies and the $L^2$ norms in \eqref{oignon} are bounded, $(\rbh_n)_{n\in\N}$ is bounded in $H^1(\R)$ thanks to \eqref{sec}.
\endgroup  Passing possibly to a subsequence, still denoted 
   $(\rbh_n)_{n \in \N}$, we may therefore  assume that there exists some $\rbh \in H^1(\R)$   such that 
 $\langle {\color{black}\rbh}, {\color{black}w'}\rangle=0$ and 
 \begin{equation}
 \label{vide} 
 \left\{
 \begin{aligned}
 &\rbh_n    \underset{n \to + \infty} \rightharpoonup \rbh {\rm \ in \ }  H^1(\R) {\rm \ weakly, \ } 
 {\rm  \ with \ } 
   \rbh \in \ker A^\perp, {\rm \ and \ }  \\
  & \rbh_n \to \rbh  {\rm \ strongly \ in \ } L^2(K), {\rm \ for \ any \ compact \ set \ } K \subset \R.
\end{aligned}
\right.
\end{equation}
It follows  from the fact that $w$ is assumed to be minimizing that $\rd^2\E(w)$ is a non-negative quadratic form.  We {\color{black}have} by lower-semicontinuity and assumption \eqref{oignon} 
 \begin{equation}
 \label{nullosse}
 0\leq \rd^2\E(w)(\rbh, \rbh) \leq \underset {n \to +\infty}\liminf\rd^2\E({\color{black}w})(\rbh_n, \rbh_n)=0, {\rm \ hence\ }
\rd^2\E(w)(\rbh, \rbh)=0.
\end{equation}
Combining \eqref{nullosse} and \eqref{cooloss}, we deduce that $\rbh\in \ker \rbA$, so that $\rbh\in \ker \rbA\cap \ker \rbA^\perp$, and hence $\rbh=0$. Choosing the compact set $K=[-S_0, S_0]$, we deduce from   \eqref{vide}, i.e. the strong convergence of $\rbh_n$ to $0$ in ${\color{black}L^2([-S_0, S_0])}$ as $n \to +\infty$, together with assumption \eqref{oignon},  that 
\begin{equation}
\label{coloc}
\underset {n \to +\infty} \lim \Vert \rbh_n\Vert_{L^2([-S_0, S_0])}=0 {\rm \ and  \ hence \ }
\underset {n \to +\infty} \lim \Vert \rbh_n\Vert_{L^2(\R\setminus [-S_0, S_0])}=1, 
\end{equation}
  so that, in view of \eqref{sec}, we have   $\nabla^2 V ({\color{black}w})(\rbh_n, \rbh_n)\to 0$ in $L^1([-S_0, S_0])$, as $n \to +\infty$, which leads to
\begin{equation}
\label{finitude}
\underset {n \to +\infty} \liminf \int_{-S_0}^{S_0} \left [\vert \rbh_n'(s) \vert^2+\nabla^2 V ({\color{black}w})(\rbh_n, \rbh_n) \right ]\rd s
\geq 0.
\end{equation}
On the other hand,  we have, in view of \eqref{humus}
\begin{equation}
\label{osamoelle}
\begin{aligned}
\int_{{\color{black}\R\setminus [-S_0, S_0]}} \left[ \vert \rbh'_n\vert^2+\nabla^2 V ({\color{black}w})(\rbh_n, \rbh_n)\right]\rd s  
&\geq \int_{{\color{black}\R\setminus [-S_0, S_0]}} \nabla^2 V (w(s) )(\rbh_n, \rbh_n) \rd s \\
&\geq \frac{\lambda_0^-}{2}\int_{{\color{black}\R\setminus [-S_0, S_0]}} (\rbh_n, \rbh_n) \rd s \\
& \geq  \frac{\lambda_0^-}{2} \Vert \rbh_n \Vert_{L^2(\R\setminus [-S_0, S_0])}^2
\end{aligned}
\end{equation}
Combining \eqref{osamoelle},  with \eqref{finitude} and \eqref{coloc}, we are led to
$$\
\underset {n \to +\infty} \liminf \int_{{\color{black}\R}} \left [\vert \rbh_n'(s) \vert^2+\nabla^2 V ({\color{black}w})(\rbh_n, \rbh_n) \right ]\rd s \geq 
  \frac{\lambda_0^-}{2} , 
$$
This contradicts \eqref{oignon} and hence \eqref{claimL20} is established.  
\qed

\smallskip
\noindent
{\it Proof of \eqref{claimL200}}.
 We notice first that, in view of {\color{black}\eqref{sec}}, we have
\begin{equation}
\begin{aligned}
&\Vert \rbh' \Vert^2_{L^2(\R)} \leq \rd^2\E(w) (\rbh, \rbh)+ \rbC_{\rm nd} \Vert \rbh \Vert_{L^2(\R)}^2, {\rm \ so \ that  } \\
&\Vert \rbh \Vert^2_{H^1(\R)}=\Vert \rbh' {\color{black}\Vert}^2_{L^2(\R)}+\Vert \rbh\Vert^2_{L^2(\R)} \leq 
\rd^2\E(w) (\rbh, \rbh)+ (\rbC_{\rm nd}+1) \Vert \rbh \Vert_{L^2(\R)}{\color{black}^2} \\
&\leq  \left(1+\frac{\rbC_{\rm nd}+1}{\alpha_0}\right) \rd^2\E(w) (\rbh, \rbh), {\rm \ if \ } \rbh \in \ker \rbA^\perp.
\end{aligned}
\end{equation}
Setting $\displaystyle{\upalpha_\star=\left(1+\frac{\rbC_{\rm nd}+1}{\alpha_0}\right)^{-1}}$, we obtain the desired result \eqref{claimL200}.

 \end{proof}

 \subsubsection{Projection on the orbit of a solution}
 For the proof of Proposition \ref{gapounet}, we rely on a powerful  tool introduced in \cite{schatzman}.
 
 \begin{lemma}
 \label{schatzman}
 Let $s \in \{0, 1\}$ {\color{black}and $i\in\{1,2\}$}. There exists a constant $\varrho_{s}^\star>0$ such that if ${\color{black}q\in\mathbbmss X(\upsigma_0^-, \upsigma_0^+)}$ satisfies
 $ {\color{black}\dist_s (q,\mathcal C(w_i))} \leq \varrho_s^\star$, 
 then there exists a unique real number ${\mathfrak m}_{\rm min, s}(q)$ such that
 \begin{equation}
 \label{objectifatteint}
\dist_s (q,C(w_i)) \equiv \underset{m \in \R}\inf \Vert q-w_i (\cdot-m)\Vert_{H^s(\R)}
=\Vert q-w_i (\cdot-{\mathfrak m}_{\rm min, s}(q))\Vert_{H^s(\R)}, 
\end{equation}
\end{lemma}
 
 The  first-order equation for the minimizer ${\mathfrak m}_{\rm min,s}(q)$ yields the following:
 
 \begin{lemma}
 \label{projection}
 Let $q$ and  ${\mathfrak m}_{\rm min,s}{\color{black}(q)}$ be as in  Lemma \ref{schatzman}.  Then, we have
  \begin{equation*}
 \left \langle \rbh({\mathfrak m}_{\rm min,s}(q)),  w'_i(\cdot-{\mathfrak m}_{\rm min,{\color{black}s}}(q))  \right\rangle_{H^s(\R)}= 0,
 \end{equation*}
where, for $m \in \R$, we have set $\rbh(m)=q-w_i(\cdot-m)$.
 \end{lemma}
 \begin{proof}
   We introduce the function $\Upphi : \R \to \R^+$ defined by 
    $$\Upphi (m)=\Vert \rbh(m) \Vert_{H^s(\R)}^2=
   \Vert q-w_i {\color{black}(\cdot-m)}\Vert_{H^s(\R)}^2 , {\rm \ for \ } m \in \R.
   $$  
   At ${\mathfrak m}_{\rm min,s}(q)$, the first derivative of $\Upphi$ with respect to $m$ vanishes. 
    We  need therefore  to expand 
 $\displaystyle{\Vert q-w_i {\color{black}(\cdot-m)}\Vert_{H^s(\R)}^2}$; however, only the difference of the two functions is in $H^s(\R)$, not the functions themselves (they do not belong to $L^2(\R)$). We write therefore
 \begin{equation}
  q=\tilde q+v_0{\color{black},} {\rm \ and \ } w_i={\color{black}\tilde w_i}+v_0, {\rm \ where \ } \tilde q \in H^s(\R), 
 {\color{black}\tilde w_i}\in H^s(\R),\\
 \end{equation}
 where $v_0$ is {\color{black}chosen smooth} as in \eqref{bv0}, and introduce the auxiliary functions, besides $\rbh$, 
  \begin{equation*}
  \tilde \rbh(m)=\tilde q-\tilde w_i (\cdot-m) {\rm \ and \ } {\bf  \upzeta}(m)=v_0-v_0(\cdot-m),{\rm \  so \ that \ }  
 \rbh(m)=\tilde \rbh(m)+\upzeta(m).
 \end{equation*}
   We notice that the functions $\rbh(m)$, $\tilde \rbh(m)$ and $\upzeta(m)$ are in $L^2(\R)$, for any $m \in \R$, and that, moreover, in view of the regularity properties of $w_i$ and $v_0$, the maps 
   $m \mapsto {\color{black}\tilde\rbh(m)}$ and $m \mapsto \upzeta(m)$ are in $C^1(\R, H^s(\R){\color{black})}$. The computation of the derivatives  actually yields 
 \begin{equation*}
 \frac{\rd  \rbh}{\rd m}(m)= w_i'(\cdot-m),  \frac{\rd  \tilde \rbh}{\rd m}(m)= \tilde w_i'(\cdot-m), {\rm \ and \ } 
 \frac{\rd \upzeta}{\rd  m}(m)=v_0'(\cdot-m), 
 \end{equation*}
 so that the derivatives of the functions are in $H^s(\R)$.
 Expanding $\Upphi(m)$, we obtain
\begin{equation}
\Upphi (m)=\Vert \rbh(m)\Vert^2=\Vert q-w_i (\cdot-m)\Vert^2=\Vert \tilde \rbh(m)\Vert^2+ \Vert \upzeta(m)\Vert^2 + 2\langle  \tilde \rbh(m), \upzeta(m)\rangle,
\end{equation}
so that
\begin{equation}
\begin{aligned}
\frac12\frac {\rd \Upphi}{\rd m} (m)&= \langle \frac{\rd  \tilde \rbh}{\rd m}(m), \tilde \rbh(m) \rangle+\langle \frac{\rd \upzeta}{\rd  m}(m), \upzeta(m)\rangle+
  \langle  \frac{\rd  \tilde \rbh}{\rd m}(m), \upzeta(m)\rangle+  \langle  \tilde \rbh(m), \frac{\rd \upzeta}{\rd m}(m)\rangle \\
  &=\langle  \frac{\rd  \tilde \rbh}{\rd m}(m), \tilde \rbh(m)+\upzeta(m)\rangle+\langle \frac{\rd \upzeta}{\rd  m}(m),\tilde \rbh(m)+ \upzeta(m)\rangle \\
  &=\langle  \frac{\rd  \tilde \rbh}{\rd m}(m)+ \frac{\rd \upzeta}{\rd  m}(m),  \rbh(m)\rangle=\langle  \frac{\rd   \rbh}{\rd m}(m),  \rbh(m)\rangle \\
 & =\langle w_i'(\cdot-m), \rbh(m)\rangle.
\end{aligned}
\end{equation}
Since  $\displaystyle{\frac {\rd \Upphi}{\rd m} (\mathfrak m_{\rm min, s})=0}$, the conclusion follows.
\end{proof}

\subsubsection{Proof of Proposition \ref{gapounet} completed}
 The {\color{black}argument} is by contradiction. If  the statement were not true, then, there {\color{black}would exist} a sequence of {\color{black}solutions} $(\rbq_n)_{n \in \N}$  to  \eqref{heterocline}, {\color{black}belonging} to $\mathbbmss X(\upsigma_0^-, \upsigma_0^+)\setminus \mathcal Z(\upsigma_0^-, {\color{black}\upsigma_0^+})$ such that 
 \begin{equation}
 \label{pasnouspasnous}
 \E(\rbq_n, \R) \to \bm (\upsigma_0^-, \upsigma_0^+), {\rm \ as \ } n \to +\infty.
 \end{equation}
 It follows from Corollary \ref{voisinage0} that $\displaystyle{\underset{i=1, 2}\inf \{\mathfrak D_0 (\rbq_n, w_i)\}\to 0}$, as $n \to +\infty$. Passing possibly to a subsequence, we may assume, without loss of generality that
 \begin{equation}
 \label{limqn}
 \mathfrak D_0 (\rbq_n, w_1 )\to 0 {\rm \ as \ } n \to + \infty, 
 \end{equation}
 so that there exists  some $n_0 \in \N$ such that, for $n \geq n_0$, we have 
 $\displaystyle{\mathfrak D_0 (\rbq_n, w_1 )} \leq \varrho_0^\star$, where $\varrho_0^\star$ is the constant given in  Lemma \ref{schatzman}.  {\color{black}In the rest of the proof set $i=1$. Hence we may define $\mathfrak m_n:={\mathfrak m}_{\rm min,0}(\rbq_n)$, where the latter is defined in \eqref{objectifatteint}. Moreover,} in view of Lemma \ref{projection}, we have 
  \begin{equation}
  \label{kernelo}
  \left \langle \rbh_n,  w'_i (\cdot-\mathfrak m_n) \right\rangle_{L^2(\R)}= 0,
  \end{equation}
    {\color{black}where} 
   $\rbh_n=\rbq_n-w_i (\cdot-\mathfrak m_n)$.
   We divide the proof into Steps.
   
   \medskip
   \noindent
   {\it   Step 1.  The spectral gap}. 
   Set $w_{i, n}=w_i (\cdot-\mathfrak m_n)$, so that $w_{i, n}$ is a minimizing solution to \eqref{heterocline}.
   In view of assumption ${\color{black}(\mathbf{H}_{\mathbf{5}})}$, we have $\ker \rbA (w_i)=\R w'_i$, which implies, by translation 
   $\ker \rbA (w_{i, n})=\R w'_{i, n}$.   Equation \eqref{kernelo}   implies that 
   $$ \rbh_n \in \ker \rbA(w_{i, n})^\perp.$$
    We are in a position to apply Proposition \ref{spectralou} to $\rbh_n$  to assert that
 \begin{equation}
 \label{trotinette}
   \rd^2\E({\color{black}w_{i,n}})(\rbh_n, \rbh_n)\geq  \alpha_\star(w_i) \Vert \rbh_n \Vert_{H^1(\R)}^2, {\rm \ if \ } {\color{black}\rbh_n \in \ker\rbA(w_{i,n})^\perp}
   \end{equation}
   
   \smallskip
   \noindent
   {\it Step 2. Energy estimates}. Recall that $\rbq_n$  satisfies  \eqref{heterocline} and so does 
   $w_{i, n}=w_i(\cdot -\mathfrak m_n)$.  Subtracting the two equations, we are led  to  
   $$-\frac{\rd^2  \rbh_n}{\rd s{\color{black}^2}}+ \nabla_u V(\rbq_n)-\nabla_u V(w_{i, n})=0, {\rm \ on \ } \R.$$
   Multiplying this equation by $\rbh_n$  and integrating on $\R$, we are led to 
   \begin{equation}
   \label{analysis}
  \int_{\R} \vert \rbh_n'\vert^2+ \left( \nabla_u V(\rbq_n)-\nabla_uV(w_{i, n})\right){\color{black}\cdot}\rbh_n{\color{black}\,\rd s}=0.
   \end{equation}
   
   \noindent
   {\it Step 3. Analysis of \eqref{analysis}}. Given ${\color{black}\rby} \in \R^k$ and $\rba=(a_1, \ldots, a_k)\in \R^k$, we have by  local  Taylor expansion
   $$\frac{\partial V}{\partial y_i}(\rby+\rba)- \frac{\partial V}{\partial  y_i}(\rby)
   =\underset {j=1}{\overset {k} \sum} \frac{\partial^2V}{\partial y_i \partial y_j}(\rby) {\color{black}a_j}+ 
   \underset {\rba \to 0} \rbo (\vert \rba \vert)
    $$
    and hence
    \begin{equation}
    \begin{aligned}
    \nabla_u [V(\rby+\rba)-V(\rby)]\cdot\rba&=
    \underset {i=1}{\overset {k}\sum} {\color{black}\underset {j=1}{\overset {k}\sum}}
     \frac{\partial^2V}{\partial y_i \partial y_j}(\rby) a_i a_j+  
   \underset {\rba \to 0} \rbo (\vert \rba \vert^2) \\
&=\rba \nabla^2 {\color{black}V}(\rby) \rba+  
   \underset {\rba \to 0} {\rbo} (\vert \rba \vert^2).
\end{aligned}
 \end{equation}
{\color{black}By \eqref{limqn} and the definition of $\mathfrak m_n$, $\Vert\rbh_n\Vert_{L^2}\to0$. The uniform derivative bound \eqref{smoothy-1d} then gives $\Vert\rbh_n\Vert_{L^\infty}\to0$. Thus the Taylor remainder is uniform along the bounded profiles. }Going back to \eqref{analysis}, we are led to
 \begin{equation}
   \label{analysis10}
 \rd^2\E(w_{i, n}{\color{black})}(\rbh_n, \rbh_n)= \int_{\R} \vert \rbh_n'(s)\vert^2+ \rbh_n(s)\nabla^2 {\color{black}V}(w_{i, n}(s))\rbh_n(s) \rd s=
  {\color{black}\underset {n\to+\infty}{\rm o}\bigl(\Vert\rbh_n\Vert_{L^2(\R)}^2\bigr)}.
   \end{equation}
 
 \smallskip
 \noindent
 {\it Step 4. End of the proof}.    Comparing \eqref{trotinette} and \eqref{analysis10}, we obtain a contradiction, so that  the Proposition is proved.
 \qed
   
\section{The mountain-pass geometry} 
\label{geometry}

\subsection{Upper bound on $\rbc_\rL$: Proof of Proposition  \ref{Collevel}}
The aim is to construct a comparison  path $\mathfrak p_0 \in \mathcal P$ such that 
\begin{equation}
\label{pathp0}
\underset{s\in [0, 1]}{\rm Max} \E(\mathfrak p_0(s), \Lambda_\rL)\leq  \rbC_{\rm pas}+ \rL{\color{black}\bm}, 
\end{equation}
for some constant  $\rbC_{\rm pas} >0$, to be determined. We
start with a map $\mathfrak q_0$  in $C^1(\Lambda_2, \R^k)$, where $\Lambda_2=[-1, 1]\times \R$, {\color{black}chosen constant in $x_1$ near each boundary,} such that 
\begin{equation}
\label{pifou}
\left\{
\begin{aligned}
\mathfrak q_0({\color{black}-1}, x_2)&=w_1(x_2) , {\rm \ for \ every \ } x_2\in \R \\
\mathfrak q_0({\color{black}1}, x_2)&=w_2(x_2) , {\rm \ for \ every \ } x_2\in \R, {\rm \ and \ } \\
     \E(\mathfrak q_0, \Lambda_2)&<+\infty.
\end{aligned}
\right.
\end{equation}
The map $\mathfrak q_0$  represents hence also a  continuous path in ${\color{black}C^1([-1, 1],\mathbbmss X(\upsigma_0^-, \upsigma_0^+))}$ connecting  $w_1$ to $w_2$. 
We extend $\mathfrak q_0$   as a  map ${\color{black}\R^2}\to \R^k$, still denoted $\mathfrak q_0$, by setting 
\begin{equation}
\label{pifou2}
\left\{
\begin{aligned}
\mathfrak q_0(x_1, x_2)&=w_1(x_2) , {\rm \ for \ every \ } x_1<-1\in \R \\
\mathfrak q_0(x_1, x_2)&=w_2(x_2) , {\rm \ for \ every \ } x_1>1\in \R,  \\
\end{aligned}
\right.
\end{equation}
so that, combining \eqref{pifou} and \eqref{pifou2}, {\color{black}we deduce} that $E(\mathfrak q_0, \Lambda_\rL)<+\infty$, for every $\rL>2$.  More precisely, we have
\begin{equation}
\label{pifou3}
\E(\mathfrak q_0, \Lambda_\rL)=\E(\mathfrak q_0, \Lambda_2)+ (\rL-2)\bm, {\rm \ for \  every \ }  \rL\geq 2.
\end{equation}
Next, we introduce, for given $\rL>0$,  a moving half-cylinder $\Lambda_\rL^+(z)$ defined by
\begin{equation}
\label{pifou4}
\Lambda_\rL^+(z)=[z, z+\frac \rL 2]\times \R, {\rm \ for \ } \rL>0, {\rm \ and \  } z \in \R, 
\end{equation}
 consider the restriction of $\mathfrak q_0$ to $\Lambda_\rL^+(z)$, translate it back to $z=0$, that is, we consider 
 the map $\mathfrak p_\rL(z): \Lambda_\rL^+(0)\to \R^k$ defined by
 \begin{equation}
 \label{pifou5}
 \mathfrak p_\rL(z) (x)={\color{black}\mathfrak q_0(x_1+z, x_2)}, {\rm \ for \ } x=(x_1, x_2)\in \Lambda_\rL^+(0), 
 \end{equation}
where $\be_1=(1, 0)$. We notice that
\begin{equation}
\label{pif0-half}
\left\{
\begin{aligned}
\mathfrak p_\rL(z) (x)&=\bw_2(x)={\color{black}w_2(x_2)}, {\rm \ for \ } x=(x_1, x_2),  x_1\in [0, \frac \rL 2], {\rm \ if \ } z\geq 1,\\
\mathfrak p_\rL(z) (x)&=\bw_1(x)={\color{black}w_1(x_2)}, {\rm \ for \ } x=(x_1, x_2),  x_1\in [0, \frac \rL 2],  {\rm \ if \ }  
z\leq -1-\frac{\rL}{2}.
\end{aligned}
\right.
\end{equation}
Moreover, we have  the energy bound
\begin{equation}
\E(\mathfrak p_\rL(z), \Lambda_\rL^+(0)) \leq {\color{black}\E(\mathfrak q_0, \Lambda_2)+ \frac{\rL}{2}\bm}.
\end{equation}
Finally, we extend the definition of $\mathfrak p_\rL(z)$ to the full cylinder $\Lambda_\rL$ by a reflection with respect to the axis $x_1=0$, that is we set 
\begin{equation}
\mathfrak p_\rL(z)(x_1, x_2)=\mathfrak p_\rL(z)(-x_1, x_2), {\rm \ for \ } x_1 \in [-\frac \rL 2, 0] , x_2 \in \R.
\end{equation}
We observe that, for every $z \in \R$, by construction $\mathfrak p_\rL(z)(-\frac \rL 2, x_2)=\mathfrak p_\rL(z)(+\frac \rL 2, x_2)$, and that 
 \begin{equation}
\label{pif0}
\left\{
\begin{aligned}
\mathfrak p_\rL(z) (x)&=\bw_2(x)={\color{black}w_2(x_2)}, {\rm \ for \ } x=(x_1, x_2),  x_1\in [0, \frac \rL 2], {\rm \ if \ } z\geq 1,\\
\mathfrak p_\rL(z) (x)&=\bw_1(x)={\color{black}w_1(x_2)}, {\rm \ for \ } x=(x_1, x_2),  x_1\in [0, \frac \rL 2],  {\rm \ if \ }  
z\leq -1-\frac{\rL}{2}.
\end{aligned}
\right.
\end{equation}
Moreover, we have 
\begin{equation}
\E(\mathfrak p_\rL(z), \Lambda_\rL^+(0)) \leq {\color{black}\E(\mathfrak q_0, \Lambda_2)+ \frac{\rL}{2}\bm}.
\end{equation}
$\E(\mathfrak p_\rL(z), \Lambda_\rL)= 2\E(\mathfrak p_\rL(z), \Lambda_\rL^+(0))$ and, by \eqref{pif0},  we have 
$\mathfrak p_\rL(z)=\bw_1$, for ${\color{black}z\leq-1-\frac{\rL}{2}}$ and $\mathfrak p_\rL(z)=\bw_2$, for $z\geq1$. Summing up, we obtain
\begin{equation}
\label{pathou}
\left\{
\begin{aligned}
&\mathfrak p_\rL(z) \in \mbWL ({\color{black}\upsigma_0^-, \upsigma_0^+})=H^1_{\rm per} ({\color{black}\Lambda_\rL})+ \bv_0, {\rm \ for \ any \ } z \in \R \\
&\E(\mathfrak p_\rL(z), \Lambda_\rL)\leq {\color{black}2\E(\mathfrak q_0, \Lambda_2)+\rL\bm}, {\rm \ for \  any \ } z \in \R, \\
&  \mathfrak p_\rL\in C^0 (\R,  \mbWL ({\color{black}\upsigma_0^-, \upsigma_0^+})), \\
&\mathfrak p_\rL(z)=\bw_1{\rm \ for \ } z< -1-\frac{\rL}2, {\rm \ and \ } \mathfrak p_\rL(z)=\bw_2,{\rm \  for \ } z\geq1.
\end{aligned}
\right.
\end{equation}
Finally, we conclude reparametrising the path as follows:
\begin{equation*}
\mathfrak p_0 (s)=\mathfrak p_\rL\left(-1-\frac{\rL}{2}+(2 +\frac {\rL}{2}) s\right), {\rm \ for \ } s\in [0, 1],
\end{equation*}
and deduce from \eqref{pathou} that $\mathfrak p_0$ has the desired property, setting 
$\displaystyle{\rbC_{\rm pass}=2\E\left(\mathfrak q_0, {\color{black}\Lambda_2}\right)}$.
\qed
\subsection{Lower bound on $\rbc_\rL$: Proof of Proposition \ref{lowerbound}}
The proof of Proposition \ref{lowerbound} relies in large part on corresponding results for the one-dimensional case.

\subsubsection{Slicing and reduction to one-dimensional problem} 
Let $\rL\geq 1$,  and set $\displaystyle{\rbI_\rL=[-\frac {\rL}{2}, \frac {\rL}{2}]}$. Consider an arbitrary   map 
$u\in \WL(\upsigma_0^-, \upsigma_0^+)$. For each 
$a \in \rbI_\rL$, we consider
 the restriction of $u$ to the line $\{a\}\times \R$. It  is well-defined, thanks to the trace theorem. We introduce therefore the map
 ${\rm   \bf Tr}_a(u): \R\to \R^k$ defined by
 \begin{equation}
 \label{thrace}
 {\rm   \bf Tr}_a(u)(s)=u(a, s), {\rm \ for \  every \ } a \in \rbI_\rL=[-\frac {\rL}{2}, \frac {\rL}{2}],  s \in \R.
 \end{equation}
 Moreover, the trace theorem implies:
 
 \begin{lemma}
 \label{thracont}
  The map  $a\mapsto \Tr_a$ is continuous from $\rbI_\rL\to {\mathbbmss v}_0+H^{\frac12}(\R)$. Moreover, for $i=1, 2$ the maps 
 $\Upsilon_i: \rbI_\rL\to \R^+$ defined by
 $\Upsilon_i(s)= \dist_0 (\Tr_s(u), \mathcal C(w_i))$ are  continuous functions.
 \end{lemma}
In view of Fubini's theorem, for almost every $a \in\rbI_\rL$, ${\rm \bf Tr}_a{\color{black}(u)}$ belongs to 
$\mathbb X({\color{black}\upsigma_0^-, \upsigma_0^+})$, hence the results of Section \ref{one-d} apply. For $i=1, 2$, we introduce the subsets 
\begin{equation}
\label{definemathA}
\left\{
\begin{aligned}
\mathbbmss A_i\left ( u,\uprho\right)&=\left\{a\in  \rbI_\rL ,{\color{black}\dist_0(\Tr_a(u),\mathcal C(w_i))\leq\uprho}\right\}, {\rm \ and }  \\
 \mathbbmss A_{\rm un}\left ( u,\uprho\right) &=\underset{i=1}{\overset {2}\cup} \mathbbmss A_i\left ( u,\uprho\right),   
    {\rm \ for \ } \uprho>0.
\end{aligned}
\right.
\end{equation}
The following is an important ingredient in the proof of Proposition \ref{lowerbound}.
 
\begin{proposition}
\label{niouininon}
 Assume that  ${\color{black}(\mathbf{H}_{\mathbf{1}})}$, ${\color{black}(\mathbf{H}_{\mathbf{2}})}$,   ${\color{black}(\mathbf{H}_{\mathbf{3}})}$ and ${\color{black}(\mathbf{H}_{\mathbf{4}})}$ hold. Let $\rL\geq 1$.  There exist constants  $0<\rbC_{\rm sink}$ and $0<\uprho_0$, 
 depending only on the potential $V$ such that the following property holds: Let $u \in \WL(\upsigma_0^-, \upsigma_0^+)$ satisfy
\begin{equation}
\E(u) \leq \rL \bm({\color{black}\upsigma_0^-, \upsigma_0^+})+ \rbC_{\rm sink}{\color{black}.} 
\end{equation}
Then   there exists  a unique   $i_0=i_0(u)\in \{1, 2\}$, such that we have 
\begin{equation}
\label{grognon}
\vert \mathbbmss A_{i_0}(u, \rho_0) \vert \geq \rL- \frac{1}{4},
 {\rm \ and \ } 
 \vert \mathbbmss A_{j_0} \vert =0 {\rm \ where \ } \{j_0\}=\{1, 2\} \setminus \{i_0\}.
\end{equation}
Moreover, the map $u\mapsto {\color{black}i_0(u)}$ is continuous on the level set  
\begin{equation}
\label{Esink}
\E^{\rbc_{\rm sink}}=\{ u \in \WL(\upsigma_0^-, \upsigma_0^+), \, \E(u, \Lambda_\rL) < \rbc_{\rm sink}\equiv \rL \bm+ \rbC_{\rm sink} \}
\end{equation}
\end{proposition}
 In other words, the proposition states that if the  energy is sufficiently small, then, for most numbers $a \in \rbI_\rL$, the trace
 $\Tr_a (u)$ is close to  the same   minimizers $w_1$ or  $w_2$, and none is close to the other one.   In this respect, it can be seen as an extension  to the two-dimensional case of Corollary \ref{voisinage0}.
  The proof  is divided into several Lemmas.
  
\begin{lemma}  
\label{journal}
Let $a$ and $b$ be in $\rbI_\rL$, with $a<b$, and assume that $\Tr_a (u) \in \mathcal N_0(z_1, \uprho)$ and 
$\Tr_b(u) \in \mathcal N_0(z_2, \uprho)$, where $z_1, z_2$ are two given  maps in $\mathbb X(\upsigma_0^-, \upsigma_0^+)$. \par Set 
$\displaystyle{\uprho_1(z_1, z_2) \equiv \frac 14 \mathfrak D_0 (z_1, z_2)}$. Then
we have, for  $\displaystyle{0\leq \uprho \leq \uprho_1 (z_1, z_2)}$
\begin{equation}
\label{perfecto}
\int_{[a, b]\times \R}\left \vert \frac{\partial u}{\partial x_1} (\ell, s) \right\vert^2 \rd \ell \rd s\geq
\frac{1}{4\vert b-a\vert} \mathfrak D_0 (z_1,z_2)^2.
\end{equation}
\end{lemma}  
  \begin{proof} We write, for $s\in \R$, and invoking the Cauchy-Schwarz inequality, 
  \begin{equation*}
  \begin{aligned}
\vert  \Tr_b(u)(s)-\Tr_a(u)(s)\vert^2&=\vert u(b, s)-u(a, s)\vert^2=\left(\int_a^b  \frac{\partial u}{\partial x_1}(\ell, s) \rd \ell\right)^2\\
 &\leq \vert b-a \vert  \int_a^b  \left \vert \frac{\partial u}{\partial x_1}(\ell, s)\right\vert^2 \rd \ell.
  \end{aligned}
  \end{equation*}
 Integrating with respect to $s$, we obtain 
 \begin{equation}
 \label{gratindau}
\int_{[a, b]\times \R}\left \vert \frac{\partial u}{\partial x_1} (\ell, s) \right \vert^2 \rd \ell \rd s \geq 
\vert b-a \vert^{-1} \Vert  \Tr_b(u)-\Tr_a(u)\Vert_{L^2(\R)}^2.
\end{equation}
Applying inequality \eqref{ultimatum} to $\Tr_a(u)$ and $\Tr_b(u)$, with $s=0$,  we are led to 
\begin{equation}
\label{derdeder}
\Vert  \Tr_b(u)-\Tr_a(u)\Vert_{L^2(\R)}^2 \geq \left(\mathfrak D_0(z_1, z_2)-2\uprho\right)^2.
\end{equation}
Combining \eqref{derdeder} with \eqref{gratindau},  together with the choice $0\leq \uprho<\uprho_1$, we derive the result.
\end{proof}
\begin{lemma} Assume that  ${\color{black}(\mathbf{H}_{\mathbf{1}})}$, ${\color{black}(\mathbf{H}_{\mathbf{2}})}$, {\color{black}${\color{black}(\mathbf{H}_{\mathbf{3}})}$}, and ${\color{black}(\mathbf{H}_{\mathbf{4}})}$ hold. Set 
\begin{equation}
\label{defrho1}
\uprho_0=\frac{1}{4}\mathfrak D_0(w_1, w_2).
\end{equation}
Let $\uprho \leq \uprho_0$, $u\in \WL(\upsigma_0^-, \upsigma_0^+)$, and assume that there exists, for $i=1, 2$, $c_i\in \rbI_\rL$ such that $\Tr_{c_i}(u)\in \mathcal N_0(w_i, \uprho)$. Then we have 
\begin{equation}
\label{perfectoide}
\E(u,  \Lambda_\rL) \geq \rL  \bm+\frac{1}{\sqrt{2}}\mathfrak D_0 (w_1, w_2) \sqrt{ \tilde\Uptheta (\uprho)}.
\end{equation}
\end{lemma}
\begin{proof}   We  may assume, without loss of generality, that $c_1<c_2$. We divide the proof in steps.

\medskip
\noindent
{\it Step 1. There 
exists $a_1, a_2 \in \rbI_\rL$ such that $(a_1, a_2){\color{black}\subset} \rbI_\rL \setminus \mathbbmss A_{\rm un}$, that is }
\begin{equation}
\label{tracage0}
\Tr_s (u) \not  \in  \mathcal N_0(w_1, \uprho)\cup \mathcal N_0(w_2, \uprho), {\rm \ for \ any \ } s \in (a_1, a_2).
\end{equation}
We consider the continuous  functions $\Upsilon_i$ introduced in Lemma \ref{thracont}, defined by $\Upsilon_i(s)= \dist_0 (\Tr_s(u), \mathcal C(w_i))$, for $s \in [-\rL\slash2, \rL \slash 2]$.  We have, by definition of $\mathcal N_0(w_1, \uprho)$
$$0\leq \Upsilon_1(c_1)\leq \uprho,  $$
whereas
\begin{equation}
\label{puppet}
\left\{
\begin{aligned}
\Upsilon_1(c_2)&= \dist_0 (\Tr_{c_2}(u), \mathcal C(w_1))
\geq \dist( \mathcal N_0(w_1, \uprho), \mathcal N_0(w_2, \uprho)) \\
&\geq \mathfrak D_0(w_1, w_2)-2\uprho \geq \frac{1}{2} \mathfrak D_0(w_1, w_2)\geq 2 \uprho.
\end{aligned}
\right.
\end{equation}
Hence we may set 
\begin{equation}
\left\{
\begin{aligned}
&a_1=\sup \{s \in [c_1, c_2], \Upsilon_1(s) \leq \uprho\}<c_2, {\rm \ and  \ similarly,  \  we  \ define\ } \\
&a_2=\inf \{s \in [a_1, c_2], \Upsilon_2(s) \leq \uprho\}>a_1.
\end{aligned}
\right.
\end{equation}
It follows from this definition that, for almost every  $s \in (a_1, a_2)$, we have $\Upsilon_1(s)>\uprho$ and $\Upsilon_2(s)>\uprho$, hence
$$
\Tr_s (u) \not \in  \mathcal N_0(w_1, \uprho)\cup \mathcal N_0(w_2, \uprho),
$$
 so that \eqref{tracage0} is established. 

\medskip
\noindent
{\it Step 2: We have}
\begin{equation}
\label{puppet2}
\int_{[a_1, a_2] \times \R} \frac12\left \vert \frac{\partial u}{\partial x_2} (\ell, s)\right\vert^2 + V(u)\rd \ell \rd s 
\geq \vert a_2-a_1\vert \cdot \Uptheta (\uprho).
\end{equation}
\begin{proof} In view of {\color{black}\eqref{tracage0}} and the definition \eqref{Uptheta1} of the function $\Uptheta$, we have, for almost every  $s \in (a_1, a_2)$
$$\E(\Tr_s(u), \R) \geq \Uptheta(\uprho).$$
Integrating over the interval $[a_1, a_2]$, we derive \eqref{puppet2}.
\end{proof}

\medskip
\noindent
{\it Step 3: We have}:
\begin{equation}
\label{perfectoide00}
\E(u, [a_1, a_2] \times \R) \geq \vert a_2-a_1 \vert \bm+\frac{1}{\sqrt{2}}\mathfrak D_0 (w_1, w_2)\sqrt{ \tilde\Uptheta (\uprho)}.
\end{equation}
\begin{proof}
Invoking \eqref{perfecto} for $a=a_1$ and $b=a_2$, and adding this to \eqref{puppet2},  we are led to
\begin{equation}
\label{little}
\begin{aligned}
\E(u, [a_1, a_2] \times \R) &\geq  \vert a_2-a_1\vert\cdot \Uptheta (\uprho)+\frac{1}{8\vert a_2-a_1\vert}\mathfrak D_0 (w_1, w_2)^2.\\
&\geq \vert a_2-a_1 \vert \bm+ \vert a_2-a_1\vert\cdot \tilde \Uptheta (\uprho)+\frac{1}{8 \vert a_2-a_1\vert}\mathfrak D_0 (w_1, w_2)^2. 
\end{aligned}
\end{equation}
We observe that the function $\displaystyle{ \lambda \mapsto \lambda\tilde  \Uptheta (\uprho)+\frac{1}{8\lambda} \mathfrak D_0 (w_1, w_2)^2}$ achieves its infimum on $\R^+$ at the point 
$\displaystyle{\lambda_0=\frac{1}{2\sqrt{2}} \frac{\mathfrak D_0 (w_1, w_2)}{\sqrt{\tilde \Uptheta (\uprho)}} }$.  This yields
$$
\lambda \tilde\Uptheta (\uprho)+\frac{1}{8\lambda} \mathfrak D_0 (w_1, w_2)^2 \geq \frac{1}{\sqrt{2}}
\mathfrak D_0 (w_1, w_2) \sqrt{\tilde\Uptheta (\uprho)}, {\rm \ for \  every \ } \lambda>0.  
$$
We apply the previous inequality for $\lambda=\vert a_2-a_1\vert$, and then combine it with \eqref{little}. This yields \eqref{perfectoide00}.
\end{proof}

\medskip
\noindent
{\it Step 4: We have}
\begin{equation}
\label{fastochoide}
\E(u, \Lambda_\rL \setminus [a_1, a_2] \times \R) \geq \left( \rL-\vert a_2-a_1\vert\right) \bm.
\end{equation}
\begin{proof} Since $\bm$ represents the minimal energy level, we have, for any $s \in \rbI_\rL$
$$
\E(\Tr_s (u), \R) \geq \bm
$$
Integrating  over $ \rbI_\rL\setminus [a_1, a_2]$, we obtain \eqref{fastochoide}.
\end{proof}

\medskip
\noindent
{\it Step 5. Proof of  \eqref{perfectoide} completed}. Combining \eqref{fastochoide} with \eqref{perfectoide00}, we derive \eqref{perfectoide}.
\end{proof}

\begin{proof} [Proof of  Proposition \ref{niouininon}]    We divide the proof into Steps.

\smallskip
\noindent
{\it Step 1}. {\it  Consider, for  $\upbeta>0$,  the set 
$\displaystyle{
\mathbbmss B(u, \upbeta)=\left \{ s\in \rbI_\rL, \, \tilde \E (\Tr_s(u), \R){\color{black}<} \upbeta\right\}}$.  If we impose 
 \begin{equation}
 \label{imposing}
0< {\rbC}_{\rm sink}\leq \frac 14 \tilde  \Uptheta (\uprho_0), 
 \end{equation}
then  $\displaystyle{\left\vert \mathbbmss B\left (u, \tilde\Uptheta (\uprho_0)\right) \right \vert\geq \rL -\frac{1}{4}}$.
\begin{equation}
\end{equation}
}
\begin{proof}
We have, by Fubini's theorem, with $\tilde \E(v, \R)=\E(v, \R)-\bm$ (see definition \eqref{renormalized})
\begin{equation}
\label{integrite}
\int_{-\frac{\rL}{2}}^\frac{\rL}{2} \tilde \E(\Tr_s(u), \R) \rd s  \leq \E(u, \Lambda_\rL)-\rL\bm \leq \rbC_{\rm sink}.
\end{equation}
We recall  that if $f$ is a non-negative function on an interval $I$, then we have, for any $\rbM>0$
\begin{equation}
\left\{ 
\begin{aligned}
\left\vert (\{s \in I , f(s)\geq  \rbM \}) \right \vert&\leq  \frac{1}{\rbM} \int_{I} f(s) \rd s,  {\rm  \ hence \ } \\
{\color{black}\left|\{s\in I:f(s)<\rbM\}\right|}&\geq  \vert I \vert-\frac{1}{\rbM}\int_I f(s) \rd s.
\end{aligned}
\right.
\end{equation}
Applying this principle to \eqref{integrite}, we deduce that 
$\displaystyle{{\color{black}\left|\rbI_\rL\setminus\mathbbmss B(u,\rbM)\right|} \leq  \frac{\rbC_{\rm sink}}{\rbM}}$.
Choosing $\rbM=\tilde \Uptheta (\uprho_0)$, and imposing \eqref{imposing}, we complete Step 1.
\end{proof}

\noindent
{\it Step 2. We have }
$\displaystyle{\mathbbmss B\left (u, \tilde\Uptheta (\uprho_0)\right)\subset  \mathbbmss A_{\rm un}(u, \uprho_0) }$, and hence
$\left\vert \mathbbmss A_{\rm un}\left (u, \uprho_0\right) \right \vert\geq \rL -\frac{1}{4}.$
\begin{proof} If $s \in \mathbbmss B\left (u, \tilde\Uptheta (\uprho_0)\right)$, then we have 
$\displaystyle{\tilde \E(\Tr_s (u), \R) < \tilde\Uptheta (\uprho_0)}$. It then follows from Corollary \ref{voisinage0} that 
$$
\Tr_s (u) \in \mathcal N_1(w_1, \uprho_0) \cup  \mathcal N_1(w_2, \uprho_0) \subset \mathcal N_0(w_1, \uprho_0) \cup  \mathcal N_0(w_2, \uprho_0).
$$
Hence, by definition \eqref{definemathA}, $\displaystyle{s \in \mathbbmss A_{\rm un} (u, \uprho_0)}$, so that the proof is complete.
\end{proof}

\noindent
{\it Step 3.   If we impose
\begin{equation}
\label{encoreune}
0< {\rbC}_{\rm sink}<\frac{1}{\sqrt{2}}\mathfrak D_0 (w_1, w_2) \sqrt{ \tilde\Uptheta ({\color{black}\uprho_0})}, 
\end{equation}
 then, one of the sets $\mathbbmss A_i (u, \uprho_0)$, $i=1, 2$ has zero measure.
 }
\begin{proof} We argue by contradiction. Assume that both sets have positive measure. Then we may find, for $i=1, 2$, $c_i\in \rbI_\rL$ such that $\Tr_{c_i}(u)\in \mathcal N_0(w_i, {\color{black}\uprho_0})$.  We deduce from {\color{black}inequality \eqref{perfectoide}} that
\begin{equation}
\label{perfectoidal}
{\color{black}\E(u,\Lambda_\rL)-\rL\bm} \geq \frac{1}{\sqrt{2}}\mathfrak D_0 (w_1, w_2) \sqrt{ \tilde\Uptheta (\uprho_0)}>{\rbC}_{\rm sink}.
\end{equation}
We obtain hence a contradiction.
\end{proof}

\noindent
{\it Step 4. Proof of  \eqref{grognon}}. 
In view of {\color{black}Steps 1 and 3}, we impose 
\begin{equation}
\label{choiceCsink}
\rbC_{\rm sink} < \inf\left\{\frac 14 \tilde  \Uptheta (\uprho_0),  \frac{1}{\sqrt{2}}\mathfrak D_0 (w_1, w_2) \sqrt{ \tilde\Uptheta (\uprho_0)}\right\}. 
\end{equation}  
We deduce from 
{\color{black}Steps 2 and 3} that $\displaystyle{\vert \mathbbmss A_{i_0} (u, \uprho_0)\vert \geq \rL-\frac14}$, so that \eqref{grognon} is established. 
\qed

\medskip
\noindent 
{\it Step {\color{black}5}.   Let $u$ and $v$ be two {\color{black}elements} in $\WL(\upsigma_0^-,\upsigma_0^+)$    satisfying \eqref{Esink} and  assume that ${\color{black}i_0(u)}=1$ and ${\color{black}i_0(v)}=2$. Then we have
\begin{equation}
\label{merci}
\Vert u-v \Vert_{L^2(\Lambda_\rL)}^2 \geq \left(\rL-\frac{1}{2}\right)\left (\mathfrak  D_0(w_1, w_2)-2\uprho_0\right)^2>0.
\end{equation}
}
\begin{proof} By assumption, we have $\displaystyle{\vert\mathbbmss A_{1} (u, \uprho_0)\vert \geq \rL-\frac14}$ and 
 $\displaystyle{\vert\mathbbmss A_{2} (v, \uprho_0)\vert \geq \rL-\frac14}$. Hence, we deduce  that 
 $$\vert \mathbbmss A_{1} (u, \uprho_0)\cap\mathbbmss A_{2} (v, \uprho_0)\vert \geq \rL-\frac{1}{2}{\color{black}.}$$
 On the other hand, we deduce from Lemma \ref{glisko} that, if 
 {\color{black}$a\in\mathbbmss A_1(u,\uprho_0)\cap\mathbbmss A_2(v,\uprho_0)$, then} 
 $$\Vert \Tr_a (u-v)\Vert_{L^2(\R)}^2 \geq \left (\mathfrak  D_0(w_1, w_2)-2\uprho_0\right)^2.$$
{\color{black}Integrating this inequality over $\mathbbmss A_1(u,\uprho_0)\cap\mathbbmss A_2(v,\uprho_0)$}, we obtain \eqref{merci}.
\end{proof}

\medskip
\noindent 
{\it Step {\color{black}6}.  Continuity of ${\color{black}i_0}$ and completion of the proof}.  It follows from \eqref{merci} that ${\color{black}i_0}$ is constant in the neighborhood of any $u$ satisfying \eqref{Esink}{\color{black},} hence ${\color{black}i_0}$ {\color{black}is} continuous. 

\end{proof}

\subsection{ Proof of Proposition \ref{lowerbound} completed}
The constant $\rbC_{\rm sink}$ has been determined in \eqref{choiceCsink}. To establish the {\color{black}result \eqref{sinky}}, we argue by contradiction and assume that it does not hold, which amounts to the existence of a path 
$\displaystyle{{\color{black}p} \in C^0([0, 1],\WL(\upsigma_0^-, \upsigma_0^+ ))}$, with $p(0)=\bw_1$, $p(1)=\bw_2$, and {\color{black}such that} 
\begin{equation}
\label{contradiction00}
 \E (p(s), \Lambda_\rL) <  \rL \bm+ \rbC_{\rm sink}, {\rm \ for \ any \ } s \in {\color{black}[0,1]}.
\end{equation}
That is, $p(s) \in \E^{c_{\rm sink}}$, where $\E^{c_{\rm sink}}$ is defined in \eqref{Esink}. We are therefore in a position to apply the last statement of Proposition \ref{niouininon}, so that the map $s\mapsto {\color{black}i_0}(p(s))$ is well-defined and continuous
on $[0, 1]$. However, we have
$$
{\color{black}i_0}(p(0))={\color{black}i_0}(\bw_1)=1 {\rm \ and \ } {\color{black}i_0}(p(1))={\color{black}i_0}(\bw_2)=2, 
$$
 a contradiction, so that the proof is complete.
 \qed

\section{Proof of  the main results}
\subsection{Summary of the results obtained at this stage }
As already mentioned in Proposition \ref{genecritic},  the number $\rc_\rL$ provided by formula \eqref{mountainpass} is a generalized critical value of $\E$, i.e. satisfies \eqref{palais}, with $c=\rc_\rL$, that is 
there exists a sequence $(u_n)_{n \in \N}$ such that 
\begin{equation}
\label{palais2}
\E(u_n, \Lambda_\rL) \to \rbc_\rL, {\rm \ as \ } n \to + \infty, {\rm \ and \ } \rd\E(u_n) \to 0 {\rm \ as \ } n \to +\infty.
\end{equation}
Applying Proposition \ref{chainofsolutions}, we deduce the existence of $\mathfrak p_\star$  non constant solutions $\rbu_k$, $k=1, \ldots, \mathfrak p_\star$ to \eqref{themainthing} , with $\mathfrak p_\star \in \N^\star$,  such that 
 \begin{equation}
  \label{mainidentity2}
  \bm(\upsigma_0^-, \upsigma_0^+)\rL + \rbC_{\rm  sink} \leq  \underset{k=1}{\overset {\mathfrak p_\star} \sum } \E(\rbu_k)=
  \rbc_\rL \leq \bm(\upsigma_0^-, \upsigma_0^+) \rL+ \rbC_{\rm pass}.
\end{equation}
Moreover \eqref{losange} holds true.   We recast \eqref{mainidentity2} as 
\begin{equation}
\label{castorama}
\left\{
\begin{aligned}
&\underset{k=1}{\overset {\mathfrak p_\star} \sum } \E(\rbu_k)=\bm(\upsigma_0^-, \upsigma_0^+) \rL+ \upgamma_\rL,  \  
 \, \rbu_k \in \mathbb W({\color{black}\upsigma_k^-, \upsigma_k^+}),  \upsigma_k^\pm \in \Sigma, {\rm \ and \ }\\
 &\upsigma_k^+=\upsigma_{k+1}^-, {\rm  \ for \ }  k=1, \ldots, \mathfrak p_\star-1, \, 
\upsigma_1^-=\upsigma_0^-, {\rm \ and \  } \upsigma_{\mathfrak p_\star}^+ =\upsigma_0^+,
\end{aligned}
\right.
\end{equation}
where $\upgamma_\rL$ can be bounded above and below independently of $\rL$, namely
\begin{equation}
\label{upgammar}
    \rbC_{\rm  sink}       \leq \upgamma_\rL \leq \rbC_{\rm pass}.
\end{equation}

\subsection{Excluding homoclinic solutions }
\begingroup\emergencystretch=1em We consider a situation where we are given  a collection  of $\mathfrak p_\star$ non constant solutions $\{\rbu_k\}_{k=1, \ldots, \mathfrak p_\star}$  on $\Lambda_{\rL}$ satisfying the relations \eqref{mainidentity2}, \eqref{castorama} and \eqref{upgammar}.
Some of these solutions might be homoclinic, that is such that $\upsigma_k^-=\upsigma_k^+$ (solutions which are not {\color{black}homoclinic} are termed heteroclinic).  Let $\mathfrak m_\star$ denote the number of heteroclinic  solutions and $\mathfrak n_\star$ be the number of homoclinic solutions, so that $\mathfrak m_\star+ \mathfrak n_\star =\mathfrak p_\star$. {\color{black}Relabeling the solutions, keeping the heteroclinic ones in their original order}, we may assume that: \\
- The solutions $\rbu_k$, for $k=1,\ldots  \mathfrak m_\star$ are heteroclinic, i.e. 
$\displaystyle{\upsigma_k^-\not =\upsigma_k^+}$,  and we have moreover 
 \begin{equation}
\upsigma_k^+=\upsigma_{k+1}^-, {\rm \ for  \ }  k=1, \ldots, \mathfrak m_\star-1, \, 
\upsigma_1^-=\upsigma_0^-, {\rm \ and \  } \upsigma_{{\color{black}\mathfrak m_\star}}^+ =\upsigma_0^+.
 \end{equation}
- The solutions $\rbu_k$, for $k=\mathfrak m_\star+1, \ldots, \mathfrak p_\star$ are homoclinic.\par\endgroup

\smallskip
We start with the observation that, since $\bm(\upsigma_k^-, {\color{black}\upsigma_k^+})$ is the minimal energy for maps in 
$\mathbbmss X(\upsigma_k^-, \upsigma_k^+)$, we deduce that 
 $\E(\rbu_k)\geq \bm(\upsigma_k^-, \upsigma_k^+)\rL$,  and hence, we have
 \begin{equation}
 \label{ripolin}
 \underset{k=1}{\overset {\mathfrak m_\star} \sum } \E(\rbu_k)
 \geq\rL \left(\underset{k=1}{\overset {\mathfrak m_\star} \sum }\bm(\upsigma_k^-, \upsigma_k^+) \right)\geq \rL \bm(\upsigma_0^-, \upsigma_0^+), 
 \end{equation} 
where the last inequality is a consequence of \eqref{v3}.  On the other hand, we have, by \eqref{mainidentity}
\begin{equation*}
\bm(\upsigma_0^-, \upsigma_0^+)\rL + \rbC_{\rm  pass} \geq  \underset{k=1}{\overset {\mathfrak m_\star} \sum } \E(\rbu_k)+ 
\underset{k=\mathfrak m_\star+1}{\overset {\mathfrak p_\star} \sum }\E(\rbu_k) \geq \rL \bm(\upsigma_0^-, \upsigma_0^+), 
\end{equation*}
so that 
\begin{equation}
0\leq \underset{k=\mathfrak m_\star+1}{\overset {\mathfrak p_\star} \sum }\E(\rbu_k)\leq \rbC_{\rm pass}.
\end{equation}

\begin{lemma}   Consider a collection  of $\mathfrak p_\star$ non constant solutions $\{\rbu_k\}_{k=1, \ldots, \mathfrak p_\star}$  on $\Lambda_{\rL}$ satisfying the relations \eqref{mainidentity2}, \eqref{castorama} and \eqref{upgammar}.
There exists a constant $\rL_0>0$, such that, for $\rL\geq \rL_0$, we have $\mathfrak n_\star=0$,  that is, there are  no non constant {\color{black}homoclinic} solutions  in the collection  $\{\rbu_k\}_{k=1, \ldots, \mathfrak p_\star}$.
\end{lemma}

\begin{proof}  In view of Proposition \ref{mono}, we have, if $\rbu_k$ is homoclinic, that is,  if $k \geq \mathfrak m^\star +1$,  
\begin{equation}
\label{grokki}
\E(\rbu_k, \Lambda_\rL) \geq   \upgamma_{\rm hc} \log \rL+ \rbC_{\rm hc}, 
\end{equation}
where $\upgamma_{\rm hc}>0$. Assume  next that  $\mathfrak  n_\star\not =0$.  It follows from \eqref{grokki} that 
   $$
  \mathfrak n_\star  \upgamma_{\rm hc} \log \rL+\mathfrak n_\star  \rbC_{\rm hc} \leq  
  \underset{k=\mathfrak m_\star+1}{\overset {\mathfrak p_\star} \sum }\E(\rbu_k) \leq  \rbC_{\rm pass}, 
  $$
 and   hence
$$\displaystyle{
\log \rL \leq \frac{\rbC_{\rm pass}}{\mathfrak n_\star\upgamma_{\rm hc}}-  \frac{\rbC_{\rm hc}}{\upgamma_{\rm hc}} 
\leq 
 \frac{\rbC_{\rm pass}-\rbC_{\rm hc}}{\upgamma_{\rm hc}}.
 }$$
If 
$\displaystyle{\rL \geq \rL_0 \equiv  2\exp \left(  \frac{\rbC_{\rm pass}-\rbC_{\rm hc}}{\upgamma_{\rm hc}}\right)}$, we obtain a contradiction, so that $\mathfrak n_\star=0$, for $\rL \geq \rL_0$.
\end{proof}

\subsection{Proof of Theorem \ref{maintheo}}
\label{final}
 In order to complete the proof of Theorem \ref{maintheo}, it suffices to prove that at least one of the solutions $\rbu_k$ is not a profile, i.e. has not   the form $\rbu_k(x_1, x_2)=w_k(x_2)$, for $x_1\in \rbI_\rL$, $x_2\in \R$, where $w_k$ is a solution to \eqref{heterocline}.
 
 \smallskip
 We divide the completion of the proof of Theorem \ref{maintheo} into steps.
 
 \medskip
 \noindent
  {\it Step 1.
  There exists a constant $\rL_1\geq \rL_0$ such that, if $\rL\geq \rL_1$, then we necessarily have  $\mathfrak p_\star=1$, } 
  \begin{proof} {\color{black}Since $\rL\geq\rL_0$, the preceding subsection gives $\upsigma_k^-\not=\upsigma_k^+$ for every $k$; hence the chain is admissible in the definition of $\updelta_\star$.} We assume by contradiction that $\mathfrak p_\star>1$.   Combining \eqref{ripolin} with \eqref{castorama} and  \eqref{v3}, we then obtain
 \begin{equation*}
 \bm(\upsigma_0^-, \upsigma_0^+) \leq \underset{k=1}{\overset {\mathfrak p_\star} \sum }\bm(\upsigma_k^-, \upsigma_k^+)\leq 
 \bm (\upsigma_0^-, \upsigma_{0}^+)+\frac{{\color{black}\upgamma_\rL}}{\rL}, 
 \end{equation*}
 so that, $\updelta_\star>0$ being the constant defined in Remark \ref{remimportante} (see \eqref{deltastar} and \eqref{deltastarplus}), we have
 $$
 0<\updelta_\star\leq \underset{k=1}{\overset {\mathfrak p_\star} \sum }\bm(\upsigma_k^-, \upsigma_k^+)-\bm(\upsigma_0^-, \upsigma_0^+) \leq 
 \frac{{\color{black}\upgamma_\rL}}{\rL}\leq \frac{\rbC_{\rm pass}}{\rL}{\rm \ and \ hence \ } \rL \leq \frac{\rbC_{\rm pass}}{\updelta_\star}.
$$
We obtain  therefore a contradiction when  
$\displaystyle{\rL\geq \rL_1\equiv {\color{black}\max\left\{2,\rL_0,\frac{2\rbC_{\rm pass}}{\updelta_\star}\right\}}}$. 
 \end{proof}
  We  assume throughout the remainder of the proof that $\rL\geq \rL_1$.  Hence, it follows from Step 1 that 
$\mathfrak p_\star=1$, so that   the solution $\rbu_1\in \mathbbmss W(\upsigma_0^-, \upsigma_0^+)$ and 
\begin{equation}
\label{fantomas}
 \E(\rbu_1)=\bm(\upsigma_0^-, \upsigma_0^+){\color{black}\rL}+ \upgamma_\rL, {\rm \ with \ }  0< \upgamma_\rL \leq \rbC_{\rm pass}.
 \end{equation}
As a matter of fact, $\rbu_1$ depends on $\rL$, so that in places we use the notation $\rbu_{1}=\rbu_{1, \rL}$ to stress the dependence on $\rL$.

\medskip
\noindent
{\it   Step 2.  There exists $\rL_\star \geq \rL_1$ such that, for $\rL\geq \rL_\star$,  the map $\rbu_{1, \rL}$  is not a one-dimensional profile}.   We argue by contradiction and assume that  $\rbu_{1, \rL}$ is a one-dimensional profile, that is there exists some solution
$\rbq_\rL \in \mathbbmss X(\upsigma_0^-, \upsigma_0^+)$ to \eqref{heterocline}  such that 
\begin{equation}
\label{monomaths}
\rbu_{1, \rL}(x_1, x_2)=\rbq_\rL(x_2), {\rm \ for \ } x=(x_1, x_2) \in \Lambda_\rL.
\end{equation}
Since $\E(\rbu_{1, \rL}, \Lambda_\rL)=\rL. \E(\rbq_\rL, \R)$, we deduce from \eqref{fantomas} that
\begin{equation}
\label{brianconnet}
\E(\rbq_\rL, \R)\leq \bm (\upsigma_0^-, \upsigma_0^+)+ \frac{\upgamma_\rL}{\rL} \leq\bm (\upsigma_0^-, \upsigma_0^+)+ \frac{\rbC_{\rm pass}}{\rL}.
\end{equation}
 Let $\displaystyle{\rL_\star=\sup \{\frac{\rbC_{\rm pass}}{\uprho_\star}, \rL_1\}}$, where $\uprho_\star$ is the constant provided by Proposition \ref{gapounet}.  If $\rL\geq \rL_\star$, then, it follows from \eqref{brianconnet} that 
 \begin{equation}
 \E(\rbq_\rL, \R)\leq \bm ({\color{black}\upsigma_0^-}, \upsigma_0^+)+\uprho_\star, 
 \end{equation}
 We are therefore in a position to apply Proposition \ref{gapounet} to assert that, for some $i\in \{1, 2\}$ and $a_\rL\in \R$  we have
 $$
 \rbq_\rL=w_i(\cdot-a_\rL), {\rm \ and \ hence \  }   \E(\rbq_\rL)=\bm (\upsigma_0^-, \upsigma_0^+), 
 $$
 This contradicts \eqref{fantomas} {\color{black}and} therefore completes the proof of Theorem \ref{maintheo}.
 \qed
 

\vspace{1cm}
\begin{center}
{\bf \Large Appendix}
\end{center}

\appendix
\section{A potential fulfilling the assumptions of Theorem \ref{maintheo}}
\subsection{Explicit definition of the potential}
The potential $V$ we consider has the following  explicit form, for $\rby=(y_1, y_2) \in \R^2$,
\begin{equation}
\label{explicit}
\rbV(y_1, y_2)=\frac{1}{8}\left( 3-3y_1^2-y_2^2\right)^2+\frac12y_1^2y_2^2+
 \left(y_1^2 +y_2^2-1\right)^2(1+y_2+y_2^2).
\end{equation}
It is convenient to introduce the following  functions on $\R^2$ and $\R$ respectively defined by
\begin{equation}
\label{castex}
\left\{
\begin{aligned}
\Phi(\rby)&=\frac12 y_1(3-y_1^2-y_2^2), {\rm \ for \ } \rby \in \R^2 \\
\rbK(s)&= 1+s+s^2, {\rm \ for \ } s \in \R, 
\end{aligned}
\right.
\end{equation}
Since 
\begin{equation}
\label{nablaphi}
\nabla \Phi (\rby)=\left(\frac12 (3-3y_1^2-y_2^2), -y_1y_2 \right),{ \rm \ for \ } \rby\in \R^2,\\
\end{equation}
we may  recast  the expression of  the potential as
\begin{equation}
\label{recast}
\left\{
\begin{aligned}
\rbV(\rby)&=\frac{1}{2}\vert \nabla \Phi (\rby)\vert^2+ (\vert \rby\vert^2-1)^2\rbK(y_2)=\rbV_1(\rby)+\rbV_2(\rby), {\rm \ for \ } \rby \in \R^2,
{\rm \ where \ } \\
\rbV_1(\rby)&=\frac{1}{2} \vert \nabla \Phi (\rby) \vert^2=\frac{1}{8}(3-3y_1^2-y_2^2)^2+\frac12 y_1^2y_2^2, {\rm \ and \ } \\
  \rbV_2(\rby)&=\left(\vert \rby \vert^2-1\right)^2 \rbK(y_2).
\end{aligned}
\right.
\end{equation}
We notice moreover that 
\begin{equation}
\label{notice0}
\rbK(s) \geq \frac34,{\rm \ for \ } s \in \R, 
\end{equation}
so that 
\begin{equation}
\label{posrbV}
\rbV_2(\rby) \geq 0, \rbV_1(\rby) \geq 0, {\rm \ and \ hence \ } \rbV(\rby) \geq 0 {\rm \ for \ } \rby \in \R^2.
\end{equation}
We notice also that $\rbV_2(\rby)=0$ if and only if $\rby \in \S^1=\{\rby \in \R^2, \vert \rby\vert=1\}$, so that $\Sigma \subset \S^1$.
\subsection{The wells and their properties}
In view of the form \eqref{recast} and \eqref{notice0}, we  observe that  $\rbV(y)=0$, if and only if 
$$ 
\vert \rby \vert=1,\  y_1y_2=0, {\rm \ and \ } 3y_1^2+y_2^2=3, 
$$
so that there are only two solutions 
\begin{equation}
\Sigma=\{\upsigma_0^+, \upsigma_0^-\}, {\rm \ with \ } \upsigma_0^-=(-1, 0) {\rm \ and \ } \upsigma_0^+=(1, 0).
\end{equation}
A direct computation shows that 
\begin{equation}
\label{directcomp}
D^2\rbV (\upsigma_0^-)=D^2\rbV (\upsigma_0^+)=
\begin{pmatrix}
 17&0 \\
  0&1
\end{pmatrix}.
\end{equation}
Hence, assertions about $\rbV$  in Assumption ${\color{black}(\mathbf{H}_{\mathbf{1}})}$ as well as  property ${\color{black}(\mathbf{H}_{\mathbf{3}})}$ are satisfied.
\subsection{Behavior off the wells}
   It follows from     \eqref{recast} and \eqref{notice0} that
   \begin{equation}
   \rbV(\rby) \geq \frac34\left (\vert \rby\vert^2-1\right)^2, {\rm \ so  \ that \ }
    \rbV(\rby) \geq \frac{27}{64} \vert \rby \vert^4, {\rm \ for \ }
   \vert y \vert\geq 2.
   \end{equation}
   Hence $\rbV(\rby) \to +\infty$ as $\vert \rby \vert \to +\infty$, and assumption ${\color{black}(\mathbf{H}_{\mathbf{1}})}$ is fulfilled.  
   
   \smallskip
   We turn next to the first statement in assumption ${\color{black}(\mathbf{H}_{\mathbf{2}})}$, that is we merely have to find a suitable lower bound for 
   $\nabla V(\rby).\rby$.  Going back to the decomposition \eqref{recast},   we compute
      \begin{equation*}
   \left\{
   \begin{aligned}
   \frac{\partial \rbV_1}{\partial y_1}(\rby){\color{black}\cdot}y_1&=\frac{1}{4} \left( 18 y_1^4+10y_1^2y_2^2-18y_1^2\right),  {\rm \ and \ } \\
   \frac{\partial \rbV_1}{\partial y_2}(\rby){\color{black}\cdot}y_2 &=\frac{1}{4}\left( 2y_2^4+10y_1^2y_2^2-6y_2^2\right), 
   \end{aligned}
   \right.
   \end{equation*}
   so that 
   \begin{equation*}
   \begin{aligned}
   \nabla \rbV_1(\rby)\cdot \rby&= \frac{1}{4} \left( 18 y_1^4+2y_2^4+20 y_1^2y_2^2-18y_1^2-6y_2^2\right) \\
  & \geq \frac{1}{2} \vert \rby \vert^4-  \frac{9}{2} \vert \rby \vert^2, 
   \end{aligned}
   \end{equation*}
   whereas 
   \begin{equation*}
   \begin{aligned}
   \nabla \rbV_2(\rby)\cdot \rby&=4\vert \rby \vert^2 (\vert \rby \vert^2-1) \rbK(y_2)+ 
  \left(\vert \rby \vert^2-1\right)^2\rbK'(y_2){y_2} \\
  &=\vert \rby \vert^4 (6y_2^2+5y_2+4)- \vert \rby \vert^2 (8y_2^2+6y_2+4)+ (2y_2^2+y_2).
  \end{aligned}
  \end{equation*}
  \begingroup\color{black}
  Set $r=\vert\rby\vert$ and $s=y_2$. In addition to \eqref{notice0}, we have
  \begin{equation}
  \label{sviluppo}
  s\rbK'(s)=2s^2+s=2\left(s+\frac14\right)^2-\frac18\geq-\frac18.
  \end{equation}
  Hence, for $r\geq1$, the preceding identity gives
  \begin{equation}
  \label{polynomio}
  \nabla\rbV_2(\rby)\cdot\rby\geq3r^2(r^2-1)-\frac18(r^2-1)^2.
  \end{equation}
  Combining this with the lower bound for $\nabla\rbV_1(\rby)\cdot\rby$, we obtain
  \begin{equation}
  \label{sviluppo2}
  \nabla\rbV(\rby)\cdot\rby\geq\frac{27}{8}r^4-\frac{29}{4}r^2-\frac18,
  \qquad r\geq1.
  \end{equation}
  Since $11r^4-58r^2-1\geq41r^2-1>0$ for $r\geq3$, it follows that
  \begin{equation}
  \label{sviluppo3}
  \nabla\rbV(\rby)\cdot\rby\geq2\vert\rby\vert^4,
  \qquad\vert\rby\vert\geq3.
  \end{equation}
  This proves the first condition in ${\color{black}(\mathbf{H}_{\mathbf{2}})}$; its polynomial growth condition follows from the explicit polynomial expression of $\rbV$.
  \endgroup
   \subsection{Heteroclinic minimizing profiles}
   Consider, for $i=1, 2$,   the two   maps $w_i: \R \to \R^2$ defined by 
   \begin{equation}
   \label{tanhsech}
  w_1(s) =(\tanh s, \sech s) {\rm \ and \ } w_2(s)=(\tanh s, -\sech s), {\rm \ for \ } s \in \R, 
   \end{equation} 
   where we recall 
   $$\tanh s=\frac{\exp (s)-\exp(-s)}{\exp (s) +\exp(-s)}{\rm \ and \ }
   \sech s=\frac{2}{\exp (s)+ \exp (-s)}, $$
   so that 
   $\tanh(s) \to 1 $ as $s\to  +\infty$, $\tanh (s) \to -1$ as $s \to -\infty$, and $\sech  s \to 0$ as $\vert s \vert \to +\infty$, so that, for 
   $i=1, 2$, $w_i(s) \to \upsigma_0^\pm$, as $s \to \pm \infty$  and hence $w_i \in \mathbb X(\upsigma_0^-, \upsigma_0^+)$.  We notice also that
   $$ \
   \tanh^2 s+ \sech^2 s=1, {\rm \ and \ hence \  } w_i(s) \in \S^1, {\rm \ for \ } s\in \R.
   $$
   \begin{lemma}
   \label{kigali}
    For $i=1, 2$, the map $w_i$ is a solution to \eqref{heterocline}, that is 
   \begin{equation}
   \label{wisolution}
   w_i''=\nabla_\rby \rbV(w_i) {\rm \ on \ } \R.
   \end{equation}
   \end{lemma}
  
  For the proof, we rely on the following observation:
  
  \begin{lemma} 
  \label{kirikou}
  Let $\Phi$ be a scalar function in  $C^2(\R^k)$, $k \in \N$, let $I\subset \R$ be a nonempty interval, and 
  $v=(v_1, \ldots, v_k): I \to \R^k$ be a solution to
  \begin{equation}
  \label{separa}
  v'(s)=\nabla_\rby \Phi (v(s)) {\rm  \ for  \ }  s \in I, {\rm \ i.e. \ }  
  v'_\ell=\frac{\partial \Phi}{\partial  y_\ell}(v), {\rm \ for \  } \ell=1, \ldots, k.
  \end{equation}
  Then, we have 
  \begin{equation}
  \label{chaotique}
  v''=\frac12 \nabla_\rby  (\vert \nabla_\rby \Phi (v)\vert^2).
  \end{equation}
  \end{lemma}
  \begin{proof}[Proof of Lemma \ref{kirikou}] Differentiating  \eqref{separa}, we obtain
  \begin{equation}
  \label{alinea}
  \begin{aligned}
  v''_\ell(s)&=\frac{\rd} {\rd s}\left [\frac{\partial \Phi} {\partial y_\ell }(v(s))\right] \\
  &=\underset {i=1}{\overset{k}\sum} \frac{\partial  \Phi}{\partial y_i}(v(s))\frac{\partial^2 \Phi }{\partial y_\ell \partial y_i}(v(s)) 
  ={\color{black}\frac 12 \underset {i=1}{\overset{k}\sum} \frac{\partial}{\partial y_\ell}\left[\left(\frac{\partial\Phi}{\partial y_i}\right)^2\right](v(s))}\\
&  =\frac 12 \frac{ \partial}{\partial y_\ell}(\vert \nabla \Phi\vert^2) (v(s)), 
  \end{aligned}
  \end{equation}
  hence the result.
  \end{proof}
  
  \begin{proof}[Proof of Lemma \ref{kigali}]  Since $\rbV_2(\rby)=0$ and $\nabla \rbV_2(\rby)=0$  for $\rby \in \S^1$, it follows that {\color{black}\eqref{wisolution}} is equivalent to
  \begin{equation}
  \label{wisolution-vone}
    w''_i=\nabla_\rby \rbV_1(w_i)=\frac{1}{2} \nabla_\rby (\vert \nabla_\rby \Phi (w_i)\vert^2).
   \end{equation}
Hence, in view of Lemma \ref{kirikou}, it suffices to  check that $w_i=(w_{i, 1}, w_{i, 2})$ satisfies  the first-order equation
\begin{equation}
\label{bamako}
w_i'(s)= \nabla_\rby \Phi (w_i(s)),{\rm \ for \ }  i=1, 2, s \in \R.
\end{equation}
Indeed, since  we  have $\tanh'(s)=1-\tanh (s)^2$ and $\sech'(s)=-\tanh(s)\sech(s)$, we compute 
\begin{equation}
\begin{aligned}
w_1'(s)&=(1-\tanh^2(s), -\tanh(s)\sech(s))=(\frac{1}{2} (3-1-2w_{1, 1} ^2(s)), -w_{1, 1}(s)w_{1, 2}(s))\\
&={\color{black}\left( \frac{1}{2} \left(3-3w_{1, 1} ^2(s)-w_{1, 2}^2(s)\right), -w_{1, 1}(s)w_{1, 2}(s)\right)}, {\rm \ since \ } w_{1, 1}^2+w_{1, 2}^2=1, \\
&=\nabla_\rby \Phi (w_1(s))
\end{aligned}
\end{equation}
so that {\color{black}\eqref{wisolution}} holds for $i=1$. The computations are similar for $i=2$.
  \end{proof}
 
 \begin{remark}{\rm   In the scalar case $k=1$, the statement of Lemma \ref{kirikou} is related to the method of separation of variables as in \eqref{first-order}. In particular any non constant finite energy solution on $\R$  
  {\color{black}satisfies the first-order balance law \eqref{first-order} (see \eqref{conserve}).}
  }
 \end{remark}  
   
Next we have:
\begin{lemma} The maps $w_i$ are the unique  minimizers (up to  translation in the variable $s$) for the energy in $\mathbb X(\upsigma_0^-, \upsigma_0^+)$. 
 \end{lemma} 
 \begin{proof}  Let $v \in \mathbb X (\upsigma_0^-, \upsigma_0^+)$. In view of the definition of the energy, we have, completing the squares
 \begin{equation}
 \label{piccolo}
 \begin{aligned}
 \E(v, \R)&=\int_\R\left[\frac12 \vert v'(s) \vert^2+\frac{1}{2}\vert \nabla_\rby \Phi (v(s))\vert^2+
(1-\vert v \vert^2)^2{\color{black}\rbK (v_2(s))} \right]\, \rd s \\
 &=\frac{1}{2}\int_\R  {\color{black}\left [\vert v'(s)-\nabla_\rby \Phi (v(s))\vert^2+ 2v'(s){\color{black}\cdot} \nabla_\rby \Phi (v(s))\right]}\rd s \\
&+ \int_\R   (1-\vert v \vert^2)^2{\color{black}\rbK (v_2(s))}\, \rd s.
 \end{aligned}
 \end{equation}
 We observe that $\displaystyle{v'(s){\color{black}\cdot} \nabla_\rby \Phi(v(s))= \frac{\rd}{\rd s} {\color{black}[\Phi(v(s))]}}$, so that 
 \begin{equation}
 \label{piccolissimo}
 \int_\R v'(s){\color{black}\cdot} \nabla_\rby \Phi (v(s)) \rd s=\left[\Phi(v(s)) \right]_{-\infty}^{+\infty}={\color{black}\Phi(\upsigma_0^+)-\Phi(\upsigma_0^-)}
 =1-(-1)=2.
 \end{equation}
 Combining \eqref{piccolo} and \eqref{piccolissimo}, we are led to 
 \begin{equation}
 \label{grammaire}
 \E(v, \R)=2+\frac{1}{2}{\color{black}\int_\R \vert v'(s)-\nabla_\rby \Phi (v(s))\vert^2\rd s} 
+ \int_\R  (1-\vert v \vert^2)^2\rbK (v_2(s))\, \rd s.
 \end{equation}
By \eqref{bamako} or a direct computation, we obtain that $\E(w_i)=2$, for $i=1, 2$, and 
\begin{equation}
\label{dakar}
\E(v, \R)\geq \E(w_i), {\rm \ for \  any \ } v \in \mathbb X(\upsigma_0^-, \upsigma_0^+), 
\end{equation}
with equality in \eqref{dakar} if and only if 
\begin{equation}
\label{togo}
v'(s)=\nabla_\rby \Phi (v(s)), {\rm \ and \ } \vert v(s) \vert =1, {\rm \ for \ } s \in \R.
\end{equation}
Writing $v=(v_1, v_2)$, we obtain $v_1^2(s)+v_2^2(s)=1$, $\forall s \in \R$. Equation \eqref{togo} yields  the differential equation 
  $$
  v'_1=1-v_1^2, {\rm \ with   \ } v_1(s) \underset{s \to +\infty} \to  1, {\rm \ and \ }  v_1(s) \underset{{\color{black}s \to -\infty}} \to  -1, 
  $$
  so that $v_1(s)=\tanh (s-\tau)$, for any $s \in{\color{black}\R}$. The relation $v_1^2(s)+v_2^2(s)=1$ then yields
  $v_2^2(s)=1-\tanh^2(s-\tau)= \sech (s-\tau)^2$.  Since $v_2$ is continuous and never vanishes, its sign is constant. Hence $v(s)=w_1(s-\tau)$ or $v(s)=w_2(s-\tau)$, for $s \in \R$, as announced.
 \end{proof}
 
 \subsection{The second derivative of the energy $\E$ at $w_i$, $i=1, 2$}
 \label{seconder}
 As seen in {\color{black}\eqref{derdederdeux}}, the second derivative of the energy at a map $v \in \mathbb X(\sigma_0^-, \sigma_0^+)$ appears in its second-order expansion. In particular, for $\rbh \in H^1(\R, \R^2)$, we have
 \begin{equation}
 \label{appendix-second-variation}
 \rd^2E(w_i)(\rbh, \rbh)=\int_\R \vert \rbh'(s)\vert^2\rd s +\int_\R \rbh(s) \nabla^2 \rbV({\color{black}w_i(s)}) \rbh(s) \rd s, 
 \end{equation}
where the matrix ${\color{black}D^2\Phi} (w_i(s))$ is given, for $i=1, 2$ by
\begin{equation} 
D^2\Phi (w_i(s))=
\begin{pmatrix}
-3w_{i, 1}(s)& -w_{i, 2}(s) \\
-w_{i, 2}(s) & -w_{i, 1}(s)
\end{pmatrix}
=
\begin{pmatrix}
-3 \tanh (s)& (-1)^{i}\sech(s) \\
(-1)^{i}\sech (s) & -\tanh(s)
\end{pmatrix}.
\end{equation}
 If we perform the expansion using the expression \eqref{grammaire} of the energy, we obtain
\begin{equation}
\label{bangladesh}
 \rd^2\E(w_i) (h, h)=\int_\R \vert h'(s)-{\color{black}D^2\Phi} (w_i(s))h(s) \vert^2+8\rbK(w_{i, 2}(s))  (w_i(s){\color{black}\cdot}h(s))^2 \rd s \geq 0.
\end{equation}
  We have:
    \begin{lemma}
    \label{speculus}
    Let $i=1, 2$.  We have $\rd^2\E(w_i)(h, h)=0$ if and only if $h=\rbc w'_i$ , for some $\rbc\in \R$. Hence $\ker \rbA(w_i)=\R {\color{black}w'_i}$, so that 
    $\rbV$ satisfies assumption ${\color{black}(\mathbf{H}_{\mathbf{5}})}$.
        \end{lemma}
 
 \begin{proof} 
 We have already seen in \eqref{coree} that $\rd^2 \E(w_i)(w'_i, w'_i)=0$. {\color{black}Conversely, let $h\in H^1(\R,\R^2)$ be such that $\rd^2\E(w_i)(h,h)=0$.} It follows from formula \eqref{bangladesh}  and \eqref{notice0} that
  \begin{equation}
  \label{oklahoma}
  h'(s)=D^2\Phi (w_i(s))(h(s)), \forall s\in \R, {\rm \ and \ }  w_i(s)\cdot h(s)=0,  \forall s \in \R.
  \end{equation}
  Since $\vert w_i (s) \vert=1$ and $w_i'(s)\not =0$, for any $s \in \R$,  we deduce from the second statement in \eqref{oklahoma} that $w'_i(s)$ and $h(s)$ are colinear for any $s \in \R$, hence there exists $\rbc(s) \in \R$ such that 
  $h(s)=\rbc(s)w_i'(s)$.  It remains to show that $\rbc$ is a constant function. It follows from \eqref{alinea} that we have
  $\displaystyle{w_i''(s)=D^2\Phi (w_i(s))w_i'(s)}$.  Combining with \eqref{oklahoma}, we deduce that 
  \begin{equation}
  \label{mexico}
  \begin{aligned}
  D^2\Phi(w_i(s)) h(s)&=h'(s)=\rbc'(s)w_i'(s)+\rbc(s)w_i''(s) \\
  &=\rbc'(s)w_i'(s)+ \rbc(s)D^2\Phi (w_i(s))w_i'(s)
  \end{aligned}
  \end{equation}
 On the other {\color{black}hand}, since $h(s)=\rbc(s)w_i'(s)$, we have also
\[
{\color{black}D^2\Phi(w_i(s))h(s)=\rbc(s)D^2\Phi(w_i(s))w_i'(s).}
\]
 Combining this with \eqref{mexico}, we deduce that ${\color{black}\rbc'(s)=0}$, for  all $s\in \R$, hence ${\color{black}\rbc}$ is constant, so that the first statement is proved.
 \end{proof}
  We have hence checked that $\rbV$ satisfies all our assumptions 
  ${\color{black}(\mathbf{H}_{\mathbf{1}})}$, ${\color{black}(\mathbf{H}_{\mathbf{2}})}$, ${\color{black}(\mathbf{H}_{\mathbf{3}})}$, ${\color{black}(\mathbf{H}_{\mathbf{4}})}$, and ${\color{black}(\mathbf{H}_{\mathbf{5}})}$.
\bibliographystyle{plain}

\bibliography{biblio-BOL.bib}



\end{document}